\documentclass[a4paper,12pt,twoside]{article}
\usepackage{amssymb,amsmath,latexsym,geometry,fancyhdr,lineno,hyperref,titletoc}

\usepackage{graphicx}

\usepackage{import}
\usepackage{xifthen}
\usepackage{pdfpages}
\usepackage{transparent}
\usepackage{amsmath}
\usepackage{mathrsfs}
\contentsmargin{0pt}

\dottedcontents{section}[2em]{\vspace{-1mm}\small}{2em}{0pt}
\dottedcontents{subsection}[5em]{\vspace{-1mm}\small}{3em}{5pt}

\hypersetup{colorlinks,linkcolor= blue,citecolor=blue}
\newtheorem{theorem}{Theorem}[section]
\newtheorem{theoremx}{Theorem}[section]

\newtheorem{lemma}[theorem]{Lemma}
\newtheorem{proposition}[theorem]{Proposition}
\newtheorem{remark}[theorem]{Remark}
\newtheorem{claim}[theorem]{Claim}
\newtheorem{corollary}[theorem]{Corollary}
\newtheorem{definition}[theoremx]{Definition}

\numberwithin{equation}{section}
\allowdisplaybreaks[0]

\def\vs{\vskip4pt}

\allowdisplaybreaks[4]

\begin{document}

\begin{center}
{\bf\Large Construction of blowup solutions for Liouville systems}
\end{center}

~~

\vs\centerline{Zetao Cheng}  
\begin{center}
{\footnotesize
{Analysis and PDE Unit, Okinawa Institute of Science and Technology\\
Okinawa, Japan\\
{\em E-mail}: chengzt20@mails.tsinghua.edu.cn}}
\end{center}

\vs\centerline{Haoyu Li}  
\begin{center}
{\footnotesize
{Departamento de Matem\'atica, Universidade Federal de S\~ao Carlos\\
S\~ao Carlos-SP, 13565-905, Brazil\\
{\em E-mail}:  hyli1994@hotmail.com}}
\end{center}

\vs\centerline{Lei Zhang}
\begin{center}
{\footnotesize
{Department of Mathematics, University of Florida\\
1400 Stadium Rd, Gainesville FL 32611\\
{\em E-mail}:  leizhang@ufl.edu}}
\end{center}

{\bf\normalsize Abstract.} {\small
We study the following Liouville system defined on a flat torus
\begin{equation}
    \left\{
   \begin{array}{lr}
     -\Delta u_i=\sum_{j=1}^n a_{ij}\rho_j\Big(\frac{h_j e^{u_j}}{\int_\Omega h_j e^{u_j}}-1\Big),\nonumber\\
     u_j\in H_{per}^1(\Omega)\mbox{ for }i\in I=\{1,\cdots,n\}\nonumber,
   \end{array}
   \right.
\end{equation}
where $h_j\in C^3(\Omega)$, $h_j>0$, $\rho_j>0$ and $u=(u_1,..,u_n)$ is doubly periodic on $\partial\Omega$. The matrix $A=(a_{ij})_{n\times n}$ satisfies certain properties. One central problem about Liouville systems is whether multi-bubble solutions do exist. In this work we present a comprehensive construction of multi-bubble solutions in the most general setting.}

~~

\noindent{\bf MCS: }35J47, 35B44

\tableofcontents

\section{Introduction}

In this article we consider the following Liouville system defined on a flat torus: let $\Omega$ be a parallelogram in $\mathbb{R}^2$, $u=(u_1,...,u_n)$ be a solution of
\begin{equation}\label{e:001}
    \left\{
   \begin{array}{lr}
     -\Delta u_i=\sum_{j=1}^n a_{ij}\rho_j\Big(\frac{h_j e^{u_j}}{\int_\Omega h_j e^{u_j}}-1\Big),\\
     u_j\in H_{per}^1(\Omega)\mbox{ for }i\in I=\{1,\cdots,n\},
   \end{array}
   \right.
\end{equation}
where $h_j\in C^3(\Omega)$, $h_j>0$, $\rho_j$ are non-negative constants and $u_j$ is doubly periodic on $\partial\Omega$ for $j\in I$,  without loss of generality  the area of $\Omega$ is $1$, $H_{per}^1(\Omega)$ is the subspace of $H^1(\Omega)$ consisting of doubly periodic functions, the coefficient matrix $A=(a_{ij})_{n\times n}$ is a  matrix that has non-negative entries. Usually for Liouville system we postulate the following standard assumption:
\[\mathcal{H}_1:\quad A\mbox{ is non-negative, invertible, symmetric and irreducible.}\]
$A$ being irreducible means (\ref{e:001}) is not a union of two separate systems: there is no partition of $I=\{1,...,n\}$ such that $I=I_1\cup I_2$, $I_1\cap I_2=\emptyset$ and $a_{ij}=0$ for all $i\in I_1$ and $j\in I_2$.

There are two striking features of Liouville systems. First, Liouville systems are enormously inclusive. It is well known that Liouville systems are deeply rooted and closely related to many fields of Mathematics and Physics. With the minimal assumption
$\mathcal{H}_1$
on the non-negative coefficient matrix $A$, the general Liouville systems represent numerous models under different contexts. In geometry, when the system reduces to a single equation ($n=1$), it generalizes the renowned Nirenberg problem, which has been extensively researched over the past few decades (see \cite{barto2,barto3,bart-taran-jde,bart-taran-jde-2,ccl,chen-li-duke, dwz-ajm, kuo-lin-jdg,li-cmp,li-shafrir,wei-zhang-adv,wei-zhang-plms,wei-zhang-jems,zhangcmp,zhangccm}). In physics, Liouville systems emerge from the mean field limit of point vortices in the Euler flow (see \cite{biler,wolansky1,wolansky2,yang}) and are intricately linked to self-dual condensate solutions of the Abelian Chern-Simons model with $N$ Higgs particles \cite{kimleelee,Wil}. In biology, they appear in stationary solutions of the multi-species Patlak-Keller-Segel system \cite{wolansky3} and are important for studying chemotaxis \cite{childress}.

The second special feature of Liouville system is the structure of solutions. From the classification theorem of \cite{ChipotShafrirWolansky1997,LinZhang2010}, for a global Liouville system of $n$ components, the total integrations of $n$ components form a $n-1$ dimensional hyper-surface. This continuum of energy brings tremendous difficulty in blowup analysis as well as construction of bubbling solutions. Recently there have been breakthroughs in the blowup analysis for Liouville systems \cite{HuangZhang2022,GuZhangArxiv}. However, because of the difficulty of structure of solutions, references on the construction of blow-up solutions have been few and insufficient.  In this article we completely settle the construction of bubbling solutions in the most general setting. Namely we postulate the least assumptions and construct bubbling solutions for multi-bubble situations.

Since our construction is closely related to the breakthroughs in blowup analysis, we list them here as reference.

\subsection{Profile of Blowup solutions} In \cite{LinZhang2011} Lin and Zhang postulated this assumption, which is called \emph{strong interaction assumption}:
\[\mathcal{H}_2: \quad a^{ii}\leq 0,\,\, \forall i\in I,\quad  a^{ij}\geq 0 \,\, \forall i\neq j\in I,  \quad  \sum_{j=1}^na^{ij}\geq 0 \,\,\forall i\in I.\]
here $(a^{ij})$ is to denote $A^{-1}$. Under $\mathcal{H}_1$ and $\mathcal{H}_2$,
Lin and Zhang \cite{LinZhang2010,LinZhang2011,LinZhang2013} prove the following:
For $\rho=(\rho_1,...,\rho_n)$ satisfying
\begin{equation}\label{rho}
8\pi (N-1) \sum_{i\in I} \rho_i<\sum_{i,j\in I}a_{ij}\rho_i\rho_j<8\pi N \sum_{i\in I}\rho_i,\quad N\in \mathbb N,
\end{equation}
there is a priori estimate for all solution $u$ to (\ref{e:001}). See \cite{LinZhang2010,LinZhang2011,LinZhang2013}.
Let
\begin{equation*}
\Gamma_{N}=
\left\{
\rho\,\big|\, \rho_i>0, i\in I;\,\,  \Lambda_{I,N}(\rho)=0.
\right\}
\end{equation*}
be the $N$-th hyper-surface, where
\begin{equation}\label{lambdaIN}\Lambda_{I,N}(\rho)=4\sum_{i=1}^n\frac{\rho_i}{2\pi N}-\sum_{i=1}^n\sum_{j=1}^na_{ij}\frac{\rho_i}{2\pi N}\frac{\rho_j}{2\pi N}.
\end{equation}
We use $\mathscr{O}_{N}$ to denote the region between $\Gamma_{N-1}$ and $\Gamma_{N}$.
The possible blowup only occurs when there exists a sequence  $\rho^k\in\mathscr{O}_N$ such that $\rho^k\to\Gamma_N$ for $N\geq2$.
\begin{center}
\begin{figure}[htbp]
\centering
\includegraphics[width=0.5\textwidth]{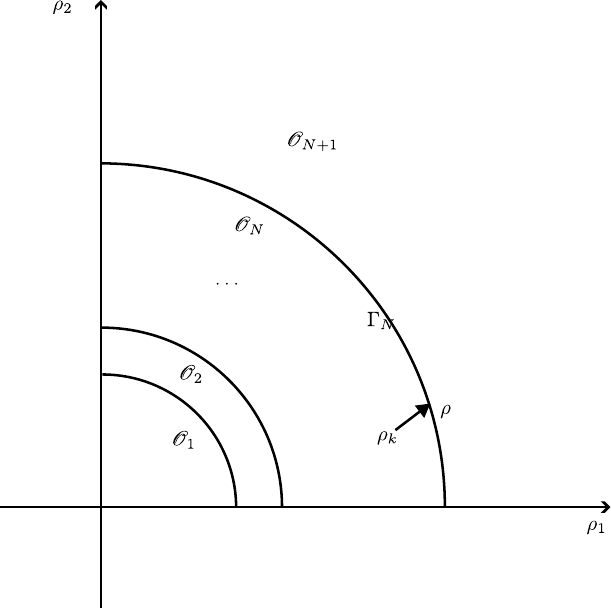}
\caption{Example of profile of blowup solutions}
\end{figure}
\end{center}
Blowup solutions are defined as follows.
\begin{definition}
A family of solutions $u^k=(u^k_1,\cdots,u_{\color{red}{n}}^k)$ to
\begin{equation}
    \left\{
   \begin{array}{lr}
     -\Delta u_i=\sum_{j=1}^n a_{ij}\rho_j^k\Big(\frac{h_j e^{u_j}}{\int_\Omega h_j e^{u_j}}-1\Big),\nonumber\\
     u_j\in H_{per}^1(\Omega)\mbox{ for }i\in I=\{1,\cdots,n\},\nonumber
   \end{array}
   \right.
\end{equation}
is called \emph{blowup solutions} if there exists  $N$ disjoint points $p^*_1,\cdots,p^*_N\in\Omega$ such that
\begin{itemize}
    \item [$(1)$] $u_i^k\to-\infty$ on any compact subset of $\Omega\backslash\{p^*_1,\cdots,p^*_N\}$ as $k\to\infty$;
    \item [$(2)$] for $i=1,\cdots,n$, $\frac{h_i e^{u_i^k}}{\int_\Omega h_i e^{u_i^k}}\to\sum_{t=1}^N \alpha_t\delta_{p^*_t}$ as $k\to\infty$. Here, the constant $\alpha_t>0$ and $\delta_{p^*_t}$ is the Dirac measure supported by $\{p^*_t\}$ for $t=1,\cdots,N$.
\end{itemize}
\end{definition}

Now we recall the works in \cite{LinZhang2010,LinZhang2011,LinZhang2013,HuangZhang2022,GuZhang2020,GuZhangArxiv} and call these results {\em sharp estimates}. To be precise, for any blowup solution $u^k=(u_1^k,\cdots,u_N^k)$ of a regular Liouville system, we denote
\begin{align}
M_{k,t}:=\max_{i\in I}\max_{x\in B_{\delta}(p_t^*)}\Big\{u_i^k(x)-\ln\int_\Omega h_i e^{u_i^k(x)}dx\Big\}\quad \mbox{ and }\quad
\varepsilon_{k,t}:=e^{-\frac{M_{k,t}}{2}}.\nonumber
\end{align}
where $\delta>0$ is small. In the following for simplicity we omit $k$ in the notations and use $\varepsilon$ to denote $\varepsilon_1,\cdots,\varepsilon_N$, and use $\varepsilon_t$, $u_{i,\varepsilon}$, $u_\varepsilon$, $\rho_{i,\varepsilon}$, $\rho_\varepsilon$ to denote $\varepsilon_{k,t}$, $u_i^k$, $u^k$, $\rho_i^k$, $\rho^k$  respectively.
Let $\rho^*\in\Gamma_N$ be the limit of $\rho_\varepsilon$ with
$\Lambda_{I,N}(\rho^*)=0$.
Moreover, we define $p_{t,\varepsilon}\in B_\delta(p_t^*)$ to be the point where \[\max_{i\in I}\max_{x\in B_{\delta}(p_t^*)}\Big\{u_i^k(x)-\ln\int_\Omega h_i e^{u_i^k(x)}dx\Big\}\]
is attained and $p_\varepsilon=(p_{1,\varepsilon},\cdots,p_{N,\varepsilon})$. Let
\[m_i^*=\sum_{j=1}^na_{ij}\frac{\rho_j^*}{2\pi N}
\quad\mbox{for}\quad  i=1,\cdots,n,
\quad \mbox{and}\quad
m^*=\min\{m^*_i|i=1,\cdots,n\}.\]
It is evident that for any $\rho^*\in\Gamma_N$, we get either $2<m^*<4$ or $m^*=4$.
The profile of blowup solutions is significantly different for $m^*<4$ and $m^*=4$. Furthermore, we denote a generic quantity by
\begin{align}\label{def:Om}
O_{m^*}=\left\{
\begin{aligned}
& O(\varepsilon^{m^*-2})&\mbox{ if } m^*<4;\\
& O(\varepsilon^{2}\ln\frac{1}{\varepsilon})&\mbox{ if }\hat m^*= 4.
\end{aligned}
\right.
\end{align}

For any $p\in\Omega$, let $G(x,p)$ be the Green's function at $p$:
\begin{equation}\label{e:GreensFunction}
    \left\{
   \begin{array}{lr}
     -\Delta_x G(x,p)=\delta_p -1\mbox{, }\int_\Omega G(x,p)dx=0,\\
     G(x,p)\mbox{ is doubly periodic on }\partial\Omega.
   \end{array}
   \right.
\end{equation}
where  $\gamma$ denotes the regular part of $G$:
\begin{align}
\gamma(x,p)=G(x,p)+\frac{1}{2\pi}\ln|x-p|,\nonumber
\end{align}
and we set for any $t=1,\cdots,N$
\begin{align}\label{def:Gt}
 {G^*}(x;p_t)=\gamma(x,p_t)+\sum_{s\neq t}G(x,p_s).
\end{align}
The first result in the \emph{sharp estimates} is on the location of blowup points:
\begin{align}\tag{$SE_1$}\label{SE1}
\sum_{i=1}^n\Big[\nabla \ln h_i(p_{t,\varepsilon})+2\pi m_i^* \nabla \tilde{G}_t(p_{t,\varepsilon};p_\varepsilon)\Big]\rho_i^*=O_{m^*}, \quad {\color{red}{t=1,...,N.}}
\end{align}
This is proved by Lin-Zhang \cite{LinZhang2013} for one-bubble case and by Huang-Zhang \cite{HuangZhang2022} and Gu-Zhang \cite{GuZhangArxiv} for multi-bubble case.
 For any $p=(p_1,\cdots,p_N)\in\Omega^N$, we define
  {
\begin{align}\label{def:H}
H_{it}(p)=2\pi m_i^* {G^*}(p_t;p_t) +\ln h_i(p_t),
\end{align}}
and let $\Omega_t$ ($t=1,\cdots,N$) be $N$ sub-domains of $\Omega$ such that $p_t^*\in\Omega_t$, $\overline{\cup_{t=1}^N\Omega_t}=\overline{\Omega}$ and $\Omega_{t_1}\cap\Omega_{t_2}=\emptyset$ for $t_1,t_2=1,\cdots,N$ with $t_1\neq t_2$.
Then we set
\begin{align}\label{NumberD}
D_{i,t}:=\Bigg[\int_{\Omega_t}\frac{\frac{h_i(x)}{h_i(p_{t}^*)}e^{2\pi m^*( {G^*}(x;p_t^*)- {G^*}(p_t^*;p_t^*))} -1}{|x-p_{t}^*|^{m^*}}dx
-\int_{\Omega_t^c}\frac{dx}{|x-p_{t}^*|^{m^*}}\Bigg]e^{I_i}
\end{align}
for $i=1,\cdots,n$ and $t=1,\cdots,N$.
Here, $I_i$ is a constant in the expansion of a global solution (see (\ref{e:LimitSystemAsy})),
 the index set $\hat {I}$ is defined as
 \[\hat{I}=\{i\in I|m^*_i=m^*\}\]
Moreover, we define
\begin{align}\label{NumeberL}
L_{i,t}:=\Bigg|\frac{\nabla h_i(p_{t}^*)}{h_i(p_{t}^*)}+8\pi\nabla_1  {G^*}(p_t^*;p_t^*)\Bigg|^2 e^{I_i}+\Bigg(\frac{\Delta h_i(p_{t}^*)}{h_i(p_{t}^*)}+8N\pi\Bigg)e^{I_i}
\end{align}
for $i=1,\cdots,n$ and $t=1,\cdots,N$.
The second result of the \emph{sharp estimates} states that
\begin{align}\tag{$SE_2$}\label{SE2}
\Lambda_{I,N}(\rho_\varepsilon)=\left\{
\begin{aligned}
&-\sum_{i\in\hat{I}}\frac{m^*-2}{N}\Big(D_{i,t}+ {\sum_{s\neq t}\frac{e^{H_{i,s}(p^*)}}{e^{H_{i,t}(p^*)}} }D_{i,s}+o(1)\Big)\varepsilon_t^{m^*-2} &\mbox{ if } m^*<4\\
&\quad\quad\mbox{ and for any }t=1,\cdots,N;\\
&-\sum_{i=1}^n\frac{2}{N}\Big(L_{i,t}+ {\sum_{s\neq t}\frac{e^{H_{i,s}(p^*)}}{e^{H_{i,t}(p^*)}} }L_{i,s}+o(1)\Big) {\varepsilon_t^{2}}\ln\frac{1}{\varepsilon_t} &\mbox{ if }\hat m^*= 4\\
&\quad\quad\mbox{ and for any }t=1,\cdots,N.
\end{aligned}
\right.
\end{align}
The third result of the \emph{sharp estimates} is:
For any $i,j=1,\cdots,n$ and any $t,s=1,\cdots,N$,
 {
\begin{align}\tag{$SE_3$}\label{SE3}
\frac{H_{i,t}(p_\varepsilon)-H_{i,s}(p_\varepsilon)}{m_i^*-2}=\frac{H_{j,t}(p_\varepsilon)-H_{j,s}(p_\varepsilon)}{m_j^*-2}+O_{m^*}.
\end{align}}
(\ref{SE2}) is proved by Huang and Zhang \cite{HuangZhang2022} and (\ref{SE3}) is proved by Gu and Zhang \cite{GuZhangArxiv}. (\ref{SE3}) is surprising since the one-bubble solutions has no relevant property .

\subsection{Statement of Main Results}

 This article aims to prove the existence of blowup solutions and demonstrate that \emph{sharp estimates} are also sufficient.
Now we list several assumptions essential to our main result.
Select any $\rho^*=(\rho^*_1,\cdots,\rho^*_n)\in(\mathbb{R}_+)^n$ such that
\begin{align}
\Lambda_{I,N}(\rho^*)=4\sum_{i=1}^n\frac{\rho^*_i}{2\pi N}-\sum_{i=1}^n\sum_{j=1}^na_{ij}\frac{\rho^*_i}{2\pi N}\frac{\rho^*_j}{2\pi N}=0.\nonumber
\end{align}
\begin{itemize}
\item [$(A_1)$] Suppose that $(p^*_1,\cdots,p^*_N)\in\Omega^N$ satisfying $p^*_t\neq p^*_s$ is a non-degenerate critical point of
\begin{align}
\sum_{t=1}^N\Big\{\sum_{i=1}^n\Big[\ln h_i(x_t)+2\pi m_i^*  {G^*}(x;p_t^*)\Big]\rho_i^*\Big\}.\nonumber
\end{align}
Moreover, for any $i,j=1,\cdots,n$ and any $t,s=1,\cdots,N$, we assume that
 {
\begin{equation}\label{new-c}
\frac{H_{it}(p^*)-H_{is}(p^*)}{m_i^*-2}=\frac{H_{jt}(p^*)-H_{js}(p^*)}{m_j^*-2}.
\end{equation}}
\end{itemize}

Based on the above notions, we get the existence of a $N$-bubble solution.
\begin{theorem}\label{t:MAIN}
Suppose that  {($\mathcal{H}_1$)} and $(A_{1})$ holds and either of the following holds
\begin{itemize}
    \item [$(1)$] $ m^*<4$ and
$\sum_{i\in\hat{I}}\sum_{t=1}^N e^{ {H_{i,t}(p^*)}}D_{i,t}\neq 0$;
    \item [$(2)$] $ m^*=4$ and
$\sum_{i=1}^n \sum_{t=1}^N e^{ {H_{i,t}(p^*)}}L_{i,t}\neq 0$;
    \item [$(3)$] $ m^*=4$,
$\sum_{i=1}^n \sum_{t=1}^N e^{ {H_{i,t}(p^*)}} L_{i,t}=0$ and
$\sum_{i=1}^N \sum_{t=1}^N e^{ {H_{i,t}(p^*)}}h_i(p_{t}^*)D_{i,t}\neq 0$.
\end{itemize}
Then, there exist $\rho_{\varepsilon}:=(\rho_{1,\varepsilon},\cdots,\rho_{n,\varepsilon})$ and  a family of solutions $\{u_{\varepsilon}:=(u_{1,\varepsilon},\cdots,u_{n,\varepsilon})\}_\varepsilon$ to Problem (\ref{e:001}) corresponding to $\rho_{\epsilon}$ such that
$\rho_{i,\varepsilon}\to\rho_i^*$ ($i=1,...,n$)
and
\begin{align}
\frac{h_i(x)e^{u_{i,\varepsilon}}}{\int_\Omega h_i(x)e^{u_{i,\varepsilon}}}\to\sum_{t=1}^N 2\pi m^*_i\delta_{p^*_t}\mbox{ in the sense of measure}\nonumber
\end{align}
as $\varepsilon\to0+$. Here, $\delta_p$ is the Dirac measure supported at $p\in\Omega$ and the numbers $D_{i,s}$ and $L_{i,s}$ are defined as in (\ref{NumberD}) and (\ref{NumeberL}).
Moreover, the solution $u_\varepsilon$ satisfies the sharp estimates (\ref{SE1}) and (\ref{SE3}). Furthermore, (\ref{SE2}) holds if Cases $(1)$ or $(2)$ occur, we get
\begin{align}
\Lambda_{I,N}(\rho_\varepsilon)=-\sum_{i=1}^n\frac{2}{N}\Big(D_{i,t}+\sum_{s\neq t}\frac{e^{ {H_{i,s}(p^*)}}}{e^{ {H_{i,t}(p^*)}}} D_{i,s}+o(1)\Big)\varepsilon_t^{2}
\end{align}
if Case $(3)$ occurs.
\end{theorem}

\begin{remark}\label{relate}
$(A1)$ and Assumptions $(1-3)$ in Theorem \ref{t:MAIN} are based on the blowup analysis of Liouville systems in
\cite{LinZhang2013,HuangZhang2022,GuZhangArxiv}. In particular, the surprising (\ref{new-c}) plays a crucial role.
\end{remark}

\begin{remark} The subtlety of our construction can be observed from
Assumption $(3)$ in Theorem \ref{t:MAIN} because it is new even for the blowup analysis in \cite{GuZhangArxiv}.
\end{remark}

\begin{remark}
$D_{i,t}$ is a \emph{global} quantity because it is involved with global integration, $L_{i,t}$ is a \emph{local} quantity because its definition only needs information at blowup point. This difference is consistent with the blowup analysis not only for Liouville systems \cite{LinZhang2013}, but also for the singular Liouville equations \cite{BartolucciYangZhang2024}
\end{remark}

\begin{remark} Because of the difficulty on the structure of solutions, there has not been too many works on the construction of blowup solutions of Liouville systems. In \cite{Huang2019} Huang constructed a sequence of one-bubble blowup solutions under further restrictions. In this article we removed all the extra assumptions and constructed the multi-bubble blowup solutions without additional restrictions.
\end{remark}

\begin{remark}
Strictly speaking, the family of solutions $\{u_\varepsilon=(u_{1,\varepsilon},\cdots,u_{n,\varepsilon})\}_\varepsilon$ depends on $N$ parameters $\varepsilon_1,\cdots,\varepsilon_N$, which also depend on each other. For the sake of simplicity, we will still use $\varepsilon$ to parameterize the solutions.
\end{remark}

The main tool we use to construct the blowup solutions is Lyapunov-Schmidt reduction, which is widely used in many planar critical problems \cite{EGP2005,delPinoMusso2005,delPinoMusso2006,LinYan2013,LinYan20131,DPR2015,Huang2019}. However our construction requires us to overcome
at least three major obstacles: First, in our construction, the approximating solutions need to be extremely precise to reduce the error. Not only do we have to apply all the subtle results in \cite{HuangZhang2022,GuZhangArxiv} about bubble interactions, in some situations we have to go further in order for the approximation to be close enough.

Second, we have to prove that certain linearized operator is invertible, which is involved with inner and outer gluing around each blowup point. When the number of equations in a Liouville system is large, the analysis is significantly harder. The number of equations in a Liouville system brings extraordinary difficulty to our construction. The higher the number, the harder the construction. In our main result we have no restriction on the size of a Liouville system.

The third major difficulty comes from the role of the Pohozaev identity. For single equation, the Pohozaev identity is a powerful tool. However, for Liouville systems, the limited information the Pohozaev identity provides is not enough for us to deal with the first derivatives of the coefficient functions.  We have to rely on Fourier analysis to handle the lack of information.

\subsection{Notations}
Throughout this article, we uniformly use the following notation.
\begin{itemize}
\item $n$ denotes the number of components of Problem (\ref{e:001}). We generally use $i$ or $j$ to index them;
\item We denote $S(u_1,\cdots,u_n):=(S_1(u_1,\cdots,u_n),\cdots,S_n(u_1,\cdots,u_n))$ with
$$S_i(U_{1,\varepsilon},\cdots,U_{n,\varepsilon})=-\Delta \Big(\sum_{j=1}^n a^{ij}U_{j,\varepsilon}\Big)-\rho_{i,\varepsilon}\Bigg(\frac{h_i e^{U_{i,\varepsilon}}}{\int_\Omega h_i e^{U_{i,\varepsilon}}}-1\Bigg)$$
for $i=1,\cdots,n$;
\item $N$ denotes the number of bubbles. We generally use $s$, $t$ or $l$ to index them;
\item $\varepsilon_t$ and $\delta_t$ denotes the scale and radii in the construction of the $t$-bubble for $t=1,\cdots,n$, respectively. See (\ref{def:ApproximateSolution1}) and (\ref{def:ApproximateSolution2});
\item $O_{m^*}$ denotes a generic quantity defined as in (\ref{def:Om});
\item $\mathbb N$ is the set of positive integers.
\end{itemize}

\section{Lyapunov-Schmidt Reduction}

In this part, we present a refined version of Lyapunov-Schmidt reduction by combining it with Fourier analysis. Such modification is necessary  because the error terms of Liouville system are more complicated than the higher dimensional semilinear elliptic equations or even mean field equations (see \cite{LinZhang2010,LinZhang2011,LinZhang2013}). We first construct the approximate solution and estimate certain integrations. A preliminary result due to the classical Lyapunov-Schmidt reduction is in Subsection \ref{Subsection:reduction} and the refined version is in Subsections \ref{subsection:RefindeReduction1} and \ref{Subsection:RefindeReduction2}.

\subsection{The Limit System}
The limit system of Problem (\ref{e:001}) is
\begin{equation}\label{e:LimitSystem}
    \left\{
   \begin{array}{lr}
     -\Delta v_i=\sum_{j=1}^n a_{ij}e^{v_j}\mbox{ for }x\in\mathbb{R}^2,\\     \frac{1}{2\pi}\int_{\mathbb{R}^2}e^{v_i}=\sigma_i \mbox{ for }i\in I=\{1,\cdots,n\}.
   \end{array}
   \right.
\end{equation}

By the results in \cite{ChipotShafrirWolansky1997,LinYan2013}, under the Assumption $\mathcal{H}_1$, entire solution $v=(v_1,\cdots,v_n)$ of (\ref{e:LimitSystem}) is radially symmetric with respect to some point $p\in\mathbb{R}^2$. Moreover, by letting $\sigma_i=\frac{1}{2\pi}\int_{\mathbb{R}^2}e^{v_i}dx$ and $\sigma=(\sigma_1,\cdots,\sigma_n)$, we get
\begin{align}
\Lambda_I(\sigma)=0\mbox{ and }\Lambda_J(\sigma)>0\mbox{ for any }\emptyset\subsetneqq J \subsetneqq I.\nonumber
\end{align}
where $\Lambda_I(\sigma)=4\sum_{i\in I}\sigma_i-\sum_{i,j\in I}a_{ij}\sigma_i\sigma_j$. It is proved in \cite{LinZhang2013} that
\[\sigma_i=\frac{\rho_i^*}{2\pi N},\quad i\in I.\]
The entire solution $v=(v_1,\cdots,v_n)$ satisfies
\begin{equation}\label{e:LimitSystemAsy}
    \left\{
   \begin{array}{lr}
     v_i(x)=-m_i^* \ln|x|+I_i-\sum_{j=1}^n\frac{a_{ij}e^{I_j}}{m_j-2}|x|^{2-m_j^*}+O(|x|^{2-m^*-\delta}),\\
     Dv_i(x)=-m_i^*\frac{x}{|x|^2}+O(|x|^{1-\hat{m}_*})\mbox{ for }i\in I=\{1,\cdots,n\}
   \end{array}
   \right.
\end{equation}
for $|x|$ large. Here, we use \cite[Lemma 4.1]{GuZhangArxiv}.
Moreover, Lin and Zhang \cite{LinZhang2013} prove the following property of the entire solution.

\begin{theorem}\label{t:KernelLocal}
Suppose $\phi=(\phi_1,\cdots,\phi_n)$ satisfies
\begin{equation}
\left\{
\begin{array}{lr}
-\Delta(\sum_{j=1}^n a^{ij}\phi_j)=e^{v_i(y)}\phi_i\mbox{ in }\mathbb{R}^2,\nonumber\\
|\phi_i(x)|\leq C(1+|x|)^\tau\mbox{ for }i=1,\cdots,n\nonumber
\end{array}
\right.
\end{equation}
for small $\tau>0$. Then, $\phi$ is a linear combination of
\begin{align}
(\partial_{x_j}v_1,\cdots,\partial_{x_j}v_n)\mbox{ for }j=1,2\nonumber
\end{align}
and
\begin{align}
(|x|\dot{v}_1(|x|)+2,\cdots,|x|\dot{v}_n(|x|)+2).\nonumber
\end{align}
\end{theorem}

\begin{corollary}\label{coro:KernelLocal}
Under the assumptions of Theorem \ref{t:KernelLocal}, if
\begin{align}\label{Condition1}
\sum_{i=1}^n\int_{\mathbb{R}^2} e^{v_i}\partial_{x_j}v_i\phi_idx=0\mbox{ for }j=1,2
\end{align}
and
\begin{align}\label{Condition2}
\sum_{i=1}^n\int_{\mathbb{R}^n}e^{v_i}(|x|\dot{v}_i(|x|)+2)\phi_idx=0,
\end{align}
then $\phi_i\equiv0$ for $i=1,\cdots,n$.
\end{corollary}

\subsection{The Approximate Solutions and The Corresponding Parameters}

Let us define a smooth radial cut-off function
\begin{equation}\label{def:CUTOFF}
\chi(x)=
\left\{
\begin{aligned}
& 1 &\mbox{ if }|x|\leq\frac{1}{2};\nonumber\\
& 0 &\mbox{ if }|x|\geq 1.\nonumber
\end{aligned}
\right.
\nonumber
\end{equation}
For any $i=1,\cdots,n$, $t=1,\cdots,N$ and $x\in B_{\delta_t}(p_{t,\varepsilon_t})$,
\begin{align}\label{def:ApproximateSolution1}
W^*_{i,t,\varepsilon}(x)=v_i\Big(\frac{x-p_{t,\varepsilon}}{\varepsilon_t}\Big)+2\ln\frac{1}{\varepsilon_t}+ 2\pi m_i^*\Big[\gamma(x,p_{t,\varepsilon})-\gamma(p_{t,\varepsilon},p_{t,\varepsilon})\Big]
\end{align}
and
\begin{align}\label{def:ApproximateSolution2}
W^{**}_{i,t,\varepsilon}(x)=v_i(\frac{\delta_t}{\varepsilon_t})+m^*_i\ln\delta_t +m^*_i\ln\frac{1}{|x-p_{t,\varepsilon}|}+2\ln\frac{1}{\varepsilon_t} + 2\pi m_i^*\Big[\gamma(x,p_{t,\varepsilon})-\gamma(p_{t,\varepsilon},p_{t,\varepsilon})\Big]
\end{align}
for $x\notin B_\frac{\delta_t}{2}(p_{t,\varepsilon})$.
We define
\begin{align}\label{def:ApproximateSolutionW}
W_{i,t,\varepsilon}(x)=W_{i,t,\varepsilon}^*(x)\chi_t(x)+W_{i,t,\varepsilon}^{**}(x)(1-\chi_t(x)).
\end{align}
Here,
\begin{align}\label{def:chi}
\chi_t(x)=\chi\big(\frac{x-p_{t,\varepsilon}}{\delta_t}\big).
\end{align}
Moreover, let
\begin{align}\label{e:Scales}
\varepsilon_1=e^{H_{1,1}(p_\varepsilon)}\varepsilon
\end{align}
and for $t=1,\cdots,N$
\begin{align}\label{def:deltat}
\delta_t=e^{H_{1,t}(p^*)-H_{1,1}(p^*)}\delta
\end{align}
for some small constant $\delta>0$.
Here, the function $H_{i,t}$ is defined as in (\ref{def:H}).
For $\varepsilon_t$ with $t=2,\cdots,N$, we will study in Hypothesis (\ref{ASSUMPTION1}).
The approximate solution is defined as
\begin{align}
U_{i,\varepsilon}(x)=\sum_{t=1}^N W_{i,t,\varepsilon}(x)\nonumber
\end{align}
for $i=1,\cdots,n$.

\begin{remark}
Notice that the footnotes $i$ and $t$ of $W_{i,t,\varepsilon}$ represent the $i$-th component and the $t$-th bubble.
\end{remark}

\begin{remark}
Throughout this article, with an obvious abuse of notations, we will denote $(\varepsilon_1,\cdots,\varepsilon)$ as $\varepsilon$ in short.
\end{remark}

We will look for the solution in the form of
\begin{align}(U_{1,\varepsilon}+w_{1,\varepsilon},\cdots, U_{n,\varepsilon}+w_{n,\varepsilon}).\nonumber
\end{align}
In a priori, we assume that
\begin{equation}\tag{$H_1$}\label{ASSUMPTION0}
|\Lambda(\rho_\varepsilon)|\leq
\left\{
\begin{aligned}
& \Lambda_0 \varepsilon^{m^*-2}&\mbox{ if } m^*<4;\\
& \Lambda_0 \varepsilon^2\ln\frac{1}{\varepsilon}&\mbox{ if }m^*= 4
\end{aligned}
\right.
\end{equation}
for a sufficiently large positive number $\Lambda_0$ and that for any $t=1,\cdots,N$
\begin{equation}\tag{$H_2$}\label{ASSUMPTION}
\sum_{i=1}^n\Big[\nabla(\ln h_i)(p_{t,\varepsilon})+2\pi m_i^*\nabla_1  {G^*}(p_{t,\varepsilon};p_{t,\varepsilon})\Big]\rho_{i}^*=
O_{m^*}\ln\frac{1}{\varepsilon}
\end{equation}
and that for
any $t,s=1,\cdots,N$ with $t\neq s$, we get
\begin{equation}\tag{$H_3$}\label{ASSUMPTION1}
(\varepsilon_{t}/\varepsilon_{s})^{m_i^*-2}-e^{ {H_{i,t}(p_\varepsilon)-H_{i,s}(p_\varepsilon)}} =\left\{
\begin{aligned}
& O(\delta^{2-\frac{m^*}{2}}\varepsilon^{m^*-2})&\mbox{ if } m^*<4;\\
& O(\delta\varepsilon^{2})&\mbox{ if }\hat m^*= 4.
\end{aligned}
\right.
\end{equation}
as $\varepsilon\to0+$.
Here, $H_{i,t}$ is defined in  (\ref{def:H}).
We remark that if $\hat m^*=4$, all $m_i^*=4$.
\begin{remark}\label{r:Equiv}
An immediate corollary of Hypothesis (\ref{ASSUMPTION1}) is that $\frac{\varepsilon_t}{\varepsilon_s}\sim 1$ for any $t,s=1,\cdots,N$.
\end{remark}

We will find a solution to Problem (\ref{e:001}) \emph{inside the interior} of these regions defined as in Hypotheses (\ref{ASSUMPTION0}), (\ref{ASSUMPTION}) and (\ref{ASSUMPTION1}).
Before we go to the next subsection, we have a discussion on $\rho_{j,\varepsilon}-\rho_j^*$ and on the error term of $\frac{\varepsilon_{t}^{m_i^*-2}}{\varepsilon_{s}^{m_i^*-2} }$.

\begin{lemma}\label{l:GammaRho}
Under Hypothesis (\ref{ASSUMPTION0}), suppose we have
\begin{align}
\frac{\rho_{i,\varepsilon}-\rho_i^*}{\rho_{j,\varepsilon}-\rho_j^*}\sim 1\quad \forall i,j\in I,\nonumber
\end{align}
then
\begin{align}
|\rho_{i,\varepsilon}-\rho_i^*|\leq
\left\{
\begin{aligned}
& \Lambda'_0 \varepsilon^{m^*-2}&\mbox{ if } m^*<4;\nonumber\\
& \Lambda'_0 \varepsilon^2\ln\frac{1}{\varepsilon}&\mbox{ if }m^*= 4\nonumber
\end{aligned}
\right.
\end{align}
for any $i\in I$. Here, $\Lambda'_0$ is a sufficiently large constant.
\end{lemma}
{\bf Proof.}
We observe that $\Lambda_{I,N}(\rho^*)=0$. If we use $s_i=\rho_i-\rho_i^*$ we see that all $s_i$ have the same sign and they are comparable to one another. Writing $\rho_i=\rho_i^*+s_i$ we have
\begin{align*}\Lambda_{I,N}(\rho)&=\Lambda_{I,N}(\rho^*)+\sum_i(2-m_i^*)\frac{s_i}{2\pi N}-\sum_{i,j}a_{ij}\frac{s_is_j}{4\pi^2N^2}\\
&=\sum_i(2-m_i^*-\sum_ja_{ij}\frac{s_j}{2\pi N})\frac{s_i}{2\pi N}.
\end{align*}Since $m_i^*>2$ and all $s_i$ have the same sign, it is easy to complete the proof by the bound of $\Lambda_{I,N}(\rho^{\varepsilon})$. $\Box$

A by-product of the above computation is that
\begin{lemma}\label{l:LambdaRho}
Assume that
\begin{align}
\frac{\rho_{i,\varepsilon}-\rho_i^*}{\rho_{j,\varepsilon}-\rho_j^*}\sim 1\quad \forall i,j\in I,\nonumber
\end{align}
we get
$\Lambda_{I,N}(\rho)-\Lambda_{I,N}(\rho^*)=\sum_{i=1}^n(2-m_i^*)\frac{\rho_i -\rho_i^*}{2\pi N}+O_{m^*}^2$.
\end{lemma}
\begin{remark}
It is worth to be pointed out that we can always find a sequence of $\rho_\varepsilon$ such that
\begin{align}
\frac{\rho_{i,\varepsilon}-\rho_i^*}{\rho_{j,\varepsilon}-\rho_j^*}\sim 1\quad \forall i,j\in I.\nonumber
\end{align}
This is due to Sard's theorem. See, for instance, \cite{HirschBook1976}.
\end{remark}

Recall that we denote
\begin{align}
S_i(U_{1,\varepsilon},\cdots,U_{n,\varepsilon})=-\Delta U_\varepsilon^i-\rho_{i,\varepsilon}\Bigg(\frac{h_i e^{U_{i,\varepsilon}}}{\int_\Omega h_i e^{U_{i,\varepsilon}}}-1\Bigg)\nonumber
\end{align}
for $i=1,\cdots,n$
and
\begin{align}
    S(u_1,\cdots,u_n)=(S_1(u_1,\cdots,u_n),\cdots,S_n(u_1,\cdots,u_n)).\nonumber
\end{align}
Here, $U_\varepsilon^i=\sum_{j=1}^n a^{ij}U_{j,\varepsilon}$ for any $i=1,\cdots,n$.

\subsection{Estimates Related with the Approximate Solution}

\begin{lemma}\label{l:Sexpansion}
It holds that
\begin{itemize}
  \item [$(1).$] For $x\in B_{\delta_t}(p_{t,\varepsilon})$,
  \begin{align}
  -\Delta U_\varepsilon^i =e^{v_i\big(\frac{x-p_{t,\varepsilon}}{\varepsilon_t}\big)+2\ln\frac{1}{\varepsilon_t}}-2\pi \sum_{t=1}^N \sum_{j=1}^n a^{ij}m_j^*+ O(\varepsilon^{m_i^*-2})\nonumber
  \end{align}
  and
\begin{align}
\frac{h_i(x)e^{U_{i,\varepsilon}}}{\int_\Omega h_i(x)e^{U_{i,\varepsilon}}}=\frac{h_i(x)}{\rho_i^* h_i(p_{t,\varepsilon})}e^{v_i\big(\frac{x-p_{t,\varepsilon}}{\varepsilon_t}\big)+2\ln\frac{1}{\varepsilon_t}+2\pi m_i^*[ {G^*}(x;p_{t,\varepsilon})- {G^*}(p_{t,\varepsilon};p_{t,\varepsilon})]}(1+O_{m^*}).\nonumber
\end{align}
Here, $ {G^*}(x;p_t)$ is defined in (\ref{def:Gt}), we use $  {G^*}(x;p_{t,\varepsilon})$ to denote
\begin{align}\label{def:G}
 {G^*}(x;p_{t,\varepsilon})=\gamma(x,p_{t,\varepsilon_t})+\sum_{s\neq t}G(x,p_{s,\varepsilon_s})
\end{align}
  \item [$(2).$] For $x\in\Omega \backslash\cup_{s=1}^N B_{\delta_s}(p_{s,\varepsilon})$,
  \begin{align}
  -\Delta U_\varepsilon^i=-2\pi \sum_{t=1}^N \sum_{j=1}^n a^{ij}m_j^*\nonumber
  \end{align}
  and
  \begin{align}
\frac{h_i(x)e^{U_{i,\varepsilon}}}{\int_\Omega h_i(x)e^{U_{i,\varepsilon}}}=\frac{N e^{I_i} h_i(x)e^{\sum_{s=1}^N 2\pi m_i^*[G(x,p_{s,\varepsilon})-\gamma(p_{s,\varepsilon},p_{s,\varepsilon})]}\big(\Pi_{s=1}^N\varepsilon_s\big)^{m_i^*-2}}{\rho_i^*\sum_{l=1}^N h_i(p_l)e^{\sum_{s\neq l}2\pi m_i^*[G(p_{l,\varepsilon},p_{s,\varepsilon})-\gamma(p_{s,\varepsilon},p_{s,\varepsilon})]}\big(\Pi_{s\neq l}\varepsilon_s\big)^{m_i^*-2}}\cdot\big(1+o_\varepsilon(1)\big).\nonumber
\end{align}
\end{itemize}
\end{lemma}
\noindent{\bf Proof.}
First of all, we notice that for any $x\in B_{\delta_t}(p_{t,\varepsilon})$ and some $t=1,\cdots,N$,
\begin{align}
-\Delta U_\varepsilon^i&=-\Delta\Big(\sum_{j=1}^n a^{ij} U_{j,\varepsilon}\Big)\nonumber\\
&=-\Delta\Big[\sum_{j=1}^n a^{ij}\chi_t\cdot W^*_{j,t,\varepsilon}+(1-\chi_t)\cdot W^{**}_{j,t,\varepsilon}\Big]\nonumber\\
&=-\Delta\Big(\sum_{j=1}^n a^{ij}W^*_{j,t,\varepsilon}\Big)\chi_t+2\nabla\Big(\sum_{j=1}^n a^{ij}W^*_{j,t,\varepsilon}\Big)\cdot\nabla\chi_t-\Big(\sum_{j=1}^n a^{ij} W^*_{j,t,\varepsilon}\Big)\Delta\chi_t\nonumber\\
&\quad-\Delta\Big(\sum_{j=1}^n a^{ij}W^{**}_{j,t,\varepsilon}\Big)+2\nabla\Big(\sum_{j=1}^n a^{ij} W^{**}_{j,t,\varepsilon}\Big)\cdot\nabla\chi_t +\Big(\sum_{j=1}^n a^{ij}W^{**}_{j,t,\varepsilon}\Big)\Delta\chi_t\nonumber\\
&=\chi_t e^{v_i\big(\frac{x-p_{t,\varepsilon}}{\varepsilon_t}\big)+2\ln\frac{1}{\varepsilon}}-2\pi\sum_{t=1}^N\sum_{j=1}^n a^{ij}m_j^*+2\nabla\chi_t\cdot\nabla\Big[\sum_{j=1}^n a^{ij}( W^*_{j,t,\varepsilon}-W^{**}_{j,t,\varepsilon})\Big]\nonumber\\
&\quad+\sum_{j=1}^n a^{ij}( W^*_{j,t,\varepsilon}-W^{**}_{j,t,\varepsilon})\Delta\chi_t.\nonumber
\end{align}
Here, for any $x\in B_{\delta_t}(p_{t,\varepsilon})\backslash B_\frac{\delta_t}{2}(p_{t,\varepsilon})$, we get
\begin{align}\label{Ineq:W-W}
W_{j,t,\varepsilon}^{**}(x)-W_{j,t,\varepsilon}^{*}(x)&=v_i(\frac{\delta_t}{\varepsilon_t})+m_i^*\ln\delta_t +m_i^*\ln\frac{1}{|x-p_{t,\varepsilon}|}-v_i\big(\frac{x-p_{t,\varepsilon}}{\varepsilon_t}\big)\nonumber\\
&=-m_i^*\ln(\frac{\delta_t}{\varepsilon_t})+I_i+O(\varepsilon^{m_i^*-2})+m_i^*\ln\delta_t+m_i^* \ln\frac{1}{|x-p_{t,\varepsilon}|}+m_i^*\ln\Big|\frac{x-p_{t,\varepsilon}}{\varepsilon_t}\Big|\nonumber\\
&\quad-I_i+O(\varepsilon^{m_i^*-2})\nonumber\\
&=O(\varepsilon^{m_i^*-2}).
\end{align}
Here, we use (\ref{def:ApproximateSolution1}), (\ref{def:ApproximateSolution2}) and (\ref{e:LimitSystemAsy}).
By a similar approach, we get
\begin{align}
\nabla \big(W_{j,t,\varepsilon}^{**}(x)-W_{j,t,\varepsilon}^{*}(x)\big)\cdot\nabla\chi_t=O(\varepsilon^{m_i^*-2}).\nonumber
\end{align}
The above implies that
\begin{align}
-\Delta U_\varepsilon^i=e^{v_i\big(\frac{x-p_{t,\varepsilon}}{\varepsilon_t}\big)+2\ln\frac{1}{\varepsilon_t}}-2\pi \sum_{t=1}^N \sum_{j=1}^n a^{ij}m_j^*+ O(\varepsilon^{m_i^*-2})\nonumber
\end{align}

By a direct computation, for $x\in \Omega\backslash\cup_{t=1}^N B_{\delta_t}(p_{t,\varepsilon})$, we get
\begin{align}
  -\Delta U_\varepsilon^i=-2\pi \sum_{t=1}^N \sum_{j=1}^n a^{ij}m_j^*\nonumber
  \end{align}

Now we estimate
$\int_\Omega h_i(x)e^{U_{i,\varepsilon}}dx$.
For $x\in B_{\delta_t}(p_{t,\varepsilon})$, we get
\begin{align}
h_i(x)e^{U_{i,\varepsilon}}&=h_i(x)e^{\sum_{s=1}^N W_{i,s,\varepsilon_s}(x)}\nonumber\\
&=h_i(x)e^{W_{i,t,\varepsilon_t}(x)} e^{\sum_{s\neq t}W_{i,s,\varepsilon_s}(x)}\nonumber\\
&=h_i(x)e^{v_i\big(\frac{x-p_{t,\varepsilon_t}}{\varepsilon_t}\big)+2\ln\frac{1}{\varepsilon_t}}e^{2\pi m_i^*[\gamma(x,p_{t,\varepsilon})-\gamma(p_{t,\varepsilon_t},p_{t,\varepsilon_t})]}\times\nonumber\\
&\quad\times e^{\sum_{s\neq t}v_i\big(\frac{\delta}{\varepsilon_s}\big)+2\ln\frac{1}{\varepsilon_s}}e^{\sum_{s\neq t}2\pi m_i^*[G(x,p_{s,\varepsilon_s})+\frac{1}{2\pi}\ln\delta -\gamma(p_{t,\varepsilon_t},p_{s,\varepsilon_s})]}e^{ (1-\chi_t)(W_{j,t,\varepsilon}^{**}-W_{j,t,\varepsilon}^*)}\nonumber\\
&=e^{v_i\big(\frac{x-p_{t,\varepsilon_t}}{\varepsilon_t}\big)+2\ln\frac{1}{\varepsilon_t}}\Big(\Pi_{s\neq t}\varepsilon_s\Big)^{m_i^*-2} e^{2\pi m_i^*\sum_{s\neq t}[G(p_{t,\varepsilon_t},p_{s,\varepsilon_s})-\gamma(p_{s,\varepsilon_s},p_{s,\varepsilon_s})]}\times\nonumber\\
&\quad\times\Big[h_i(p_{t,\varepsilon})+\nabla_y\Big(h_i(y)e^{2\pi m^*[ {G^*}(y;p_{t,\varepsilon})- {G^*}(p_{t,\varepsilon};p_{t,\varepsilon})]}\Big)\Big|_{y=p_{t,\varepsilon}}\cdot(x-p_{t,\varepsilon})+O(|x-p_{t,\varepsilon}|)^2\Big]\times\nonumber\\
&\quad\times\Big[1+(1-\chi_t)\cdot(W_{j,t,\varepsilon}^{**}-W_{j,t,\varepsilon}^*)+O(\varepsilon^{2m_i^*-4})\Big].\nonumber
\end{align}

Notice that
\begin{align}
\int_{B_\delta(p_{t,\varepsilon_t})} e^{v_i(\frac{x-p_{t,\varepsilon}}{\varepsilon_t})+2\ln\frac{1}{\varepsilon_t}}=\int_{B_{\frac{\delta}{\varepsilon}}(0)}e^{v_i(x)}dx=\frac{\rho_i^*}{N}+O(\int_\frac{\delta}{\varepsilon}^\infty \frac{rdr}{r^{m_i^*}})=\frac{\rho_i^*}{N}+O(\varepsilon^{m^*-2}).\nonumber
\end{align}
On the other hand, by a direct computation, we get
\begin{align}
\int_{B_\delta(p_{t,\varepsilon_t})}e^{v_i(\frac{x-p_{t,\varepsilon}}{\varepsilon_t})+2\ln\frac{1}{\varepsilon_t}}|x-p_{t,\varepsilon}|^2dx
&=\varepsilon_t^2 \int_{B_\frac{\delta}{\varepsilon_t}(0)} e^{v_i(x)}|x|^2dx\nonumber\\
&\leq C\varepsilon_t^2 \Big(C+\int_1^\frac{\delta}{\varepsilon}\frac{r^3 dr}{r^{m^*_i}}\Big)=
O_{m^*}\nonumber
\end{align}
Combining (\ref{Ineq:W-W}), we get for any $t=1,\cdots,N$
\begin{align}
\int_{B_\delta(p_{t,\varepsilon})}h_i(x)e^{U_{i,\varepsilon}}dx&=\frac{\rho_i^* h_i(p_{t,p_{t,\varepsilon}})}{N}e^{\sum_{s\neq t}2\pi m_i^*[G(p_{t,\varepsilon},p_{s,\varepsilon}) -\gamma(p_{s,\varepsilon},p_{s,\varepsilon})]} \Big(\Pi_{s\neq t}\varepsilon_s\Big)^{m_i^* -2}(1+O_{m^*}).\nonumber
\end{align}
Here, $O_{m^*}$ is defined in (\ref{def:Om}).

For $x\in\Omega\backslash\cup_{t=1}^N B_{\delta_t}(p_{t,\varepsilon})$, we get
\begin{align}
e^{U_{i,\varepsilon}}=e^{\sum_{t=1}^N v_i(\frac{\delta}{\varepsilon_t})+2\ln\frac{1}{\varepsilon_t}+2\pi m_i^*[G(x,p_{t,\varepsilon})+\frac{1}{2\pi}\ln\frac{1}{\delta}-\gamma(p_{t,\varepsilon},p_{t,\varepsilon})](1-\theta_{i,t})}\leq C\Pi_{t=1}^N e^{v_i(\frac{\delta}{\varepsilon_t})+2\ln\frac{1}{\varepsilon_t}+m_i\ln\delta}.\nonumber
\end{align}
Recalling Hypothesis (\ref{ASSUMPTION1}) and Remark \ref{r:Equiv}, we get
\begin{align}
\int_{\Omega\backslash\cup_{t=1}^N B_{\delta_t}(p_{t,\varepsilon})}h_i(x)e^{U_{i,\varepsilon}}dx=O(\varepsilon^{N(m^*-2)}).\nonumber
\end{align}
The above computation gives that
\begin{align}
\int_{\Omega}h_i(x)e^{U_{i,\varepsilon}}=
\frac{\rho_i^*}{N}\sum_{t=1}^N h_i(p_{t,\varepsilon}) e^{\sum_{s\neq t}2\pi m_i^*[G(p_{t,\varepsilon},p_{s,\varepsilon})-\gamma(p_{s,\varepsilon},p_{s,\varepsilon})]}\Big(\Pi_{s\neq t}\varepsilon_s\Big)^{m_i^*-2}(1+O_{m^*}).\nonumber
\end{align}

Therefore, if $x\in B_{\delta_t}(p_{t,\varepsilon})$, we get
\begin{align}\label{e:FracInBall}
&\frac{h_i(x)e^{U_{i,\varepsilon}}}{\int_\Omega h_i(x)e^{U_{i,\varepsilon}}}\nonumber\\
&=\frac{h_i(x) e^{v_i\big(\frac{x-p_{t,\varepsilon}}{\varepsilon_t}\big)+2\ln\frac{1}{\varepsilon_t}+2\pi m_i^* [ {G^*}(x;p_{t,\varepsilon})-\gamma(p_{t,\varepsilon},p_{t,\varepsilon})]-\sum_{s\neq t}2\pi m_i^*\gamma(p_{s,\varepsilon},p_{s,\varepsilon})}\Big(\Pi_{s\neq t}\varepsilon_s\Big)^{m_i^*-2}}{\frac{\rho_i^*}{N}\sum_{l=1}^N h_i(p_{l,\varepsilon}) e^{\sum_{s\neq t}2\pi m_i^*[G(p_{l,\varepsilon},p_{s,\varepsilon})-\gamma(p_{s,\varepsilon},p_{s,\varepsilon})]}\Big(\Pi_{s\neq t}\varepsilon_s\Big)^{m_i^*-2}}(1+O_{m^*})\nonumber\\
&=\frac{h_i(x) e^{v_i\big(\frac{x-p_{t,\varepsilon}}{\varepsilon_t}\big)+2\ln\frac{1}{\varepsilon_t}+2\pi m_i^* [ {G^*}(x;p_{t,\varepsilon})-\gamma(p_{t,\varepsilon},p_{t,\varepsilon})]-\sum_{s\neq t}2\pi m_i^*\gamma(p_{s,\varepsilon},p_{s,\varepsilon})}}{\frac{\rho_i^*}{N}\sum_{l=1}^N h_i(p_{l,\varepsilon}) e^{\sum_{s\neq t}2\pi m_i^*[G(p_{l,\varepsilon},p_{s,\varepsilon})-\gamma(p_{s,\varepsilon},p_{s,\varepsilon})]}\big(\frac{\varepsilon_t}{\varepsilon_l}\big)^{m_i^*-2}}(1+O_{m^*}).
\end{align}

Now we compute its denominator.
By Hypothesis (\ref{ASSUMPTION1}), we get
\begin{align}
\mbox{Denominator of (\ref{e:FracInBall})}&= {\frac{\rho_i^*}{N}}h_i(p_{t,\varepsilon})e^{\sum_{s\neq t}2\pi m_i^* [G(p_{t,\varepsilon},p_{s,\varepsilon})-\gamma(p_{s,\varepsilon},p_{s,\varepsilon})]}\nonumber\\
&\quad+\frac{\rho_i^*}{N}\sum_{l\neq t}h_i(p_{l,\varepsilon})e^{\sum_{s\neq t}2\pi m_i^* [G(p_{l,\varepsilon},p_{s,\varepsilon})-\gamma(p_{s,\varepsilon},p_{s,\varepsilon})]}\big(\frac{\varepsilon_t}{\varepsilon_l}\big)^{m_i^*-2}\nonumber\\
&=\frac{\rho_i^*}{N}h_i(p_{t,\varepsilon})\Big\{e^{\sum_{s\neq t}2\pi m_i^*[G(p_{t,\varepsilon},p_{s,\varepsilon})-\gamma(p_{s,\varepsilon},p_{s,\varepsilon})]}\nonumber\\
&\quad+\sum_{l\neq t} e^{\sum_{s\neq l}2\pi m_i^*[G(p_{l,\varepsilon},p_{s,\varepsilon})-\gamma(p_{s,\varepsilon},p_{s,\varepsilon})]}\times\frac{e^{(m_i^*-2)H_{i,t}(p_\varepsilon)}}{e^{(m_i^*-2)H_{i,s}(p_\varepsilon)}} \Big\}\cdot (1+o_\varepsilon(1))\nonumber\\
&=\rho_i^* h_i(p_{t,\varepsilon}) e^{\sum_{s\neq t}2\pi m_i^*[G(p_{t,\varepsilon},p_{s,\varepsilon})-\gamma(p_{s,\varepsilon},p_{s,\varepsilon})]}(1+o_\varepsilon(1)).\nonumber
\end{align}
Here, $H_{i,t}$ is defined as in (\ref{def:H}).
This implies that
\begin{align}
\frac{h_i(x)e^{U_{i,\varepsilon}}}{\int_\Omega h_i(x)e^{U_{i,\varepsilon}}}
=\frac{h_i(x)}{\rho_i^* h_i(p_{t,\varepsilon})}e^{v_i\big(\frac{x-p_{t,\varepsilon}}{\varepsilon_t}\big)+2\ln\frac{1}{\varepsilon_t}+2\pi m_i^*[ {G^*}(x;p_{t,\varepsilon})- {G^*}(p_{t,\varepsilon};p_{t,\varepsilon})]}(1+o_\varepsilon(1))\nonumber
\end{align}
for $x\in B_{\delta_t}(p_{t,\varepsilon_t})$.
Here, $\tilde{G}_t$ is defined as in (\ref{def:Gt}).

By a similar computation, if $x\in \Omega\backslash\cup_{s=1}^N B_{\delta_s}(p_{s,\varepsilon_s})$, we get
\begin{align}
\frac{h_i(x)e^{U_{i,\varepsilon}}}{\int_\Omega h_i(x)e^{U_{i,\varepsilon}}}=\frac{N h_i(x)e^{\sum_{s=1}^N 2\pi m_i^*[G(x,p_{s,\varepsilon_s})-\gamma(p_{s,\varepsilon_s},p_{s,\varepsilon_s})]}\big(\Pi_{s=1}^N\varepsilon_s\big)^{m_i^*-2}}{\rho_i^*\sum_{l=1}^N h_i(p_l)e^{\sum_{s\neq l}2\pi m_i^*[G(p_{l,\varepsilon_l},p_{s,\varepsilon_s})-\gamma(p_{s,\varepsilon_s},p_{s,\varepsilon_s})]}\big(\Pi_{s\neq l}\varepsilon_s\big)^{m_i^*-2}}\cdot(1+o_\varepsilon(1)).\nonumber
\end{align}
The proof is completed.
$\Box$

\begin{corollary}\label{coro:OutSmall}
It holds that
\begin{equation}
\frac{h_i(x)e^{U_{i,\varepsilon}}}{\int_\Omega h_i(x)e^{U_{i,\varepsilon}}}=
\left\{
\begin{aligned}
& O(\varepsilon^{m^*-2})&\mbox{ if } m^*<4;\\
& O(\varepsilon^{2})&\mbox{ if } {m^*= 4.}
\end{aligned}
\right.
\nonumber
\end{equation}
for $x\in\Omega\backslash\cup_{s=1}^N B_{\delta_s}(p_{s,\varepsilon})$.
\end{corollary}

\subsection{The Reduction Scheme: The Preliminary Results}\label{Subsection:reduction}
Let us denote
\begin{align}
K_{i,\varepsilon}(x)=\sum_{t=1}^N{\bf{1}}_{B_{\delta_t}(p_{t,\varepsilon})} e^{v_i(\frac{x-p_{t,\varepsilon}}{\varepsilon_t})+2\ln\frac{1}{\varepsilon_t}},\nonumber
\end{align}

\begin{align}\label{OperatorL}
L_{i,\varepsilon}(w_\varepsilon^1,\cdots,w_\varepsilon^n) = -\Delta w_\varepsilon^i-K_{i,\varepsilon}(x)w_{i,\varepsilon} +\frac{K_{i,\varepsilon}(x)}{\rho_i^*}\int_\Omega K_{i,\varepsilon}(x) w_{i,\varepsilon},
\end{align}
\begin{align}\label{Reminderg}
g_{i,\varepsilon}(x,w)&=
\Delta U_\varepsilon^i-\rho_{i,\varepsilon}+\rho_{i,\varepsilon}\frac{h_i(x)e^{ U_{i,\varepsilon}}}{\int_\Omega h_i(x)e^{U_{i,\varepsilon}}}-\Bigg(K_{i,\varepsilon}(x)-\rho_{i,\varepsilon}\frac{h_i(x)e^{ U_{i,\varepsilon}}}{\int_\Omega h_i(x)e^{ U_{i,\varepsilon}}}\Bigg)w_{i,\varepsilon}\nonumber\\
&\quad\quad\quad+\frac{K_{i,\varepsilon}(x)}{\rho_i^*}\int_\Omega K_{i,\varepsilon}(x)w_{i,\varepsilon}-\frac{\rho_{i,\varepsilon}h_i(x)e^{ U_{i,\varepsilon}}}{\Big(\int_\Omega h_i(x)e^{ U_{i,\varepsilon}}\Big)^2}\int_\Omega\Big(h_i(x)e^{ U_{i,\varepsilon}}w_{i,\varepsilon}\Big)+N_{i,\varepsilon}
\end{align}
and
\begin{align}\label{ReminderN}
N_{i,\varepsilon}&=\rho_{i,\varepsilon}\frac{h_i(x)e^{\ U_{i,\varepsilon}+w_{i,\varepsilon}}}{\int_\Omega h_i(x)e^{ U_{i,\varepsilon}+w_{i,\varepsilon}}} - \rho_{i,\varepsilon}\frac{h_i(x)e^{ U_{i,\varepsilon}}}{\int_\Omega h_i(x)e^{ U_{i,\varepsilon}}} - \rho_{i,\varepsilon}\Bigg(\frac{h_i(x)e^{ U_{i,\varepsilon}}}{\int_\Omega h_i(x)e^{ U_{i,\varepsilon}}}\Bigg)w_{i,\varepsilon}\nonumber\\
&\quad\quad\quad+\frac{\rho_{i,\varepsilon}h_i(x)e^{ U_{i,\varepsilon}}}{\Big(\int_\Omega h_i(x)e^{ U_{i,\varepsilon}}\Big)^2}\int_\Omega\Big(h_i(x)e^{ U_{i,\varepsilon}}w_{i,\varepsilon}\Big).
\end{align}
Here, ${\bf{1}}_A$ is the characteristic function on set $A$ and
\begin{itemize}
    \item $w^i_{\varepsilon}=\sum_{j=1}^n a^{ij}w_{j,\varepsilon}$;
    \item $U^i_{\varepsilon}=\sum_{j=1}^n a^{ij}U_{j,\varepsilon}$
\end{itemize}
for $i=1,\cdots,n$.

\begin{remark}
Notice that all of the coefficients in $L_{i,\varepsilon}$, $g_{i,\varepsilon}$ and $N_{i,\varepsilon}$ are double periodic.
\end{remark}

It is evident that
\begin{align}
S_i(U_{1,\varepsilon}+w_{1,\varepsilon},\cdots, U_{n,\varepsilon}+w_{n,\varepsilon})=L_{i,\varepsilon}(w_\varepsilon^1,\cdots,w_\varepsilon^n)-g_{i,\varepsilon}(x,w)\nonumber
\end{align}
for $i=1,\cdots,n$.

To build the framework for Lyapunov-Schmidt reduction, we need certain weighted Sobolev spaces depending on parameters. To begin with, let us define
\begin{align}
\rho_\beta(x)=(1+|x|)^{1+\frac{\beta}{2}}\mbox{  and  }
\widetilde{\rho}_\beta(x)=\frac{1}{(1+|x|)(\ln(2+|x|))^{1+\frac{\beta}{2}}}\nonumber
\end{align}
with $\beta\in(0,\frac{1}{2})$ fixed. For $\xi=(\xi_1,\cdots,\xi_n)\in(L^1_{loc}(\Omega))^n$, we define
\begin{align}
\|\xi\|_{\mathbb{X}_\varepsilon}^2&=\sum_{i=1}^n\sum_{t=1}^N\Bigg(\|\Delta\widetilde{\xi}_{t,i}\rho_\beta(x-p_{t,\varepsilon})\|^2_{L^2(B_\frac{2\delta_t}{\varepsilon_t}(p_{t,\varepsilon}))}+\|\widetilde{\xi}_{t,i}\widetilde{\rho}_\beta(x-p_{t,\varepsilon})\|^2_{L^2(B_\frac{2\delta_t}{\varepsilon_t}(p_{t,\varepsilon}))}\nonumber\\
&\quad\quad\quad\quad\quad+\|\Delta\xi_i\|^2_{L^2(\Omega\backslash \cup_{t=1}^N B_\frac{2\delta_t}{\varepsilon_t}(p_{t,\varepsilon}))}+\|\xi_i\|^2_{L^2(\Omega \backslash \cup_{t=1}^N B_\frac{2\delta_t}{\varepsilon_t}(p_{t,\varepsilon}))}\Bigg)\nonumber
\end{align}
and
\begin{align}
\|\xi\|_{\mathbb{Y}_\varepsilon}^2=\sum_{i=1}^n\sum_{t=1}^N\Bigg(\varepsilon_t^4\|\widetilde{\xi}_{t,i}\rho_\beta(x-p_{t,\varepsilon})\|^2_{L^2(B_\frac{2\delta_t}{\varepsilon_t}(p_{t,\varepsilon}))}+\|\xi_i\|^2_{L^2(\Omega\backslash \cup_{t=1}^N B_\frac{2\delta_t}{\varepsilon_t}(p_{t,\varepsilon}))}\Bigg).\nonumber
\end{align}
Here, $\widetilde{\xi}_{t,i}(y)=\xi_{i}(\varepsilon_t y+p_{t,\varepsilon})$ for $i=1,\cdots,n$ and $t=1,\cdots,N$.
The weighted Sobolev spaces are defined as
\begin{align}
    \mathbb{X}_\varepsilon=\{\xi\in (L^1_{loc}(\Omega))^n|\|\xi\|_{\mathbb{X}_\varepsilon}<\infty\mbox{ and }\xi_i\mbox{ is double perodic for }i=1,\cdots,n\}\nonumber
\end{align}
and
\begin{align}
    \mathbb{Y}_\varepsilon=\{\xi\in (L^1_{loc}(\Omega))^n|\|\xi\|_{\mathbb{Y}_\varepsilon}<\infty\mbox{ and }\xi_i\mbox{ is double perodic for }i=1,\cdots,n\}.\nonumber
\end{align}

We denote for $t=1,\cdots,N$
\begin{align}
Z_{1,t,\varepsilon}^*(x)&=\Big(\chi_t\partial_{x_1} v_1\Big(\frac{x-p_{t,\varepsilon}}{\varepsilon_t}\Big), \cdots, \chi_t\partial_{x_1} v_n\Big(\frac{x-p_{t,\varepsilon}}{\varepsilon_t}\Big)\Big);\nonumber\\
Z_{2,t,\varepsilon}^*(x)&=\Big(\chi_t\partial_{x_2} v_1\Big(\frac{x-p_{t,\varepsilon}}{\varepsilon_t}\Big), \cdots, \chi\partial_{x_2} v_n\Big(\frac{x-p_{t,\varepsilon}}{\varepsilon_t}\Big)\Big);\nonumber\\
Z_{3,t,\varepsilon}^*(x)&=(Z^*_{3,t,1,\varepsilon},\cdots,Z^*_{3,t,n,\varepsilon}).\nonumber
\end{align}
Here, $\chi_t$ is the cut-off function defined as in (\ref{def:CUTOFF}), $(v_1,\cdots,v_n)$ is the entire solution to Problem (\ref{e:LimitSystem}) and
\begin{align}
Z_{3,t,i,\varepsilon}^*(x)=\chi_t(x)\Big[|x-p_{t,\varepsilon}|v'_i(\frac{x-p_{t,\varepsilon}}{\varepsilon_t})+2\varepsilon_t\Big] + (1-\chi_t(x)) (2-m^*_i)\varepsilon_t\nonumber
\end{align}
for $i=1,\cdots,n$.
In follows, we denote the $i$-th component of $Z^*_{j,t,\varepsilon}$ as $Z^*_{j,t,i,\varepsilon}$ for $i=1,\cdots,n$, $j=1,2,3$ and $t=1,\cdots,N$.

Define
\begin{align}
E_\varepsilon =\Bigg\{w=(w_1,\cdots,w_n)\in\mathbb{X}_\varepsilon&\Bigg|\sum_{t=1}^N\int_\Omega \chi_t e^{v_i(\frac{x-p_{t,\varepsilon}}{\varepsilon_t})}w_i dx=0\mbox{ for }i=1,\cdots,n\mbox{ and }\nonumber\\
&\sum_{i=1}^n\int_\Omega\sum_{s=1}^N\chi_s e^{v_i(\frac{x-p_{s,\varepsilon}}{\varepsilon_s})+2\ln\frac{1}{\varepsilon_s}} Z_{1,t,i,\varepsilon}^* w_idx
\nonumber\\
&=\sum_{i=1}^n\int_\Omega\sum_{s=1}^N\chi_s e^{v_i(\frac{x-p_{s,\varepsilon}}{\varepsilon_s})+2\ln\frac{1}{\varepsilon_s}} Z_{2,t,i,\varepsilon}^* w_idx\nonumber\\
&=\sum_{i=1}^n\int_\Omega\sum_{s=1}^N\chi_s e^{v_i(\frac{x-p_{s,\varepsilon}}{\varepsilon_s})+2\ln\frac{1}{\varepsilon_s}} Z_{3,t,i,\varepsilon}^* w_idx=0\mbox{ for }t=1,\cdots,N.\Bigg\}\nonumber
\end{align}
and
\begin{align}\label{def:F}
F_\varepsilon=\Bigg\{w=(w_1,\cdots,w_n)\in\mathbb{Y}_\varepsilon&\Bigg|\int_\Omega w_i dx=0\mbox{ for }i=1,\cdots,n\mbox{ and }\nonumber\\
&\sum_{i=1}^n\int_\Omega Z_{1,t,i,\varepsilon}^* w_idx
=\sum_{i=1}^n\int_\Omega Z_{2,t,i,\varepsilon}^* w_idx\nonumber\\
&=\sum_{i=1}^n\int_\Omega Z_{3,t,i,\varepsilon}^* w_idx=0\mbox{ for }t=1,\cdots,N.\Bigg\}.
\end{align}
Define the projection operator
$Q_\varepsilon: {\mathbb{Y}_\varepsilon}\to F_\varepsilon$
by
\begin{align}\label{ProjectionQ}
Q_\varepsilon(u)=u&-\sum_{j=1}^2\sum_{t=1}^N r_{j,t}\Big(\chi_t\partial_{x_j}v_1(\frac{x-p_{t,\varepsilon}}{\varepsilon_t})\sum_{s=1}^Ne^{v_1(\frac{x-p_{s,\varepsilon}}{\varepsilon_s})+2\ln\frac{1}{\varepsilon_s}},\cdots,\nonumber\\
&\quad\quad\chi_t\partial_{x_j}v_n(\frac{x-p_{t,\varepsilon}}{\varepsilon_t})\sum_{s=1}^Ne^{v_n(\frac{x-p_{s,\varepsilon}}{\varepsilon_s})+2\ln\frac{1}{\varepsilon_s}}\Big)\nonumber\\
&-\sum_{t=1}^N r_{3,t}\Big(Z_{3,t,1,\varepsilon }^*\sum_{s=1}^N\chi_s e^{v_1(\frac{x-p_{s,\varepsilon}}{\varepsilon_s})+2\ln\frac{1}{\varepsilon_s}},\cdots,Z_{3,t,n,\varepsilon }^*\sum_{s=1}^N\chi_s e^{v_n(\frac{x-p_{s,\varepsilon}}{\varepsilon_s})+2\ln\frac{1}{\varepsilon_s}}\Big)\nonumber\\
&-\Big(s_1\sum_{t=1}^N\chi_t e^{v_1(\frac{x-p_{t,\varepsilon}}{\varepsilon_t})+2\ln\frac{1}{\varepsilon_t}},\cdots,s_n\sum_{t=1}^N\chi_t e^{v_n(\frac{x-p_{t,\varepsilon}}{\varepsilon_t})+2\ln\frac{1}{\varepsilon_t}}\Big).
\end{align}
\begin{proposition}\label{prop:Recution}
Under Hypotheses (\ref{ASSUMPTION0}), (\ref{ASSUMPTION}) and (\ref{ASSUMPTION1}), there exists a constant $\varepsilon_0>0$ such that for any $\varepsilon\in(0,\varepsilon_0)$, there exists a $w_\varepsilon=(w_\varepsilon^1,\cdots,w_\varepsilon^n)\in E_\varepsilon$ satisfies
\begin{align}\label{e:projected}
Q_\varepsilon(L_{1,\varepsilon}(w_\varepsilon)-g_{1,\varepsilon}(x,w_\varepsilon),\cdots,L_{n,\varepsilon}(w_\varepsilon)-g_{n,\varepsilon}(x,w_\varepsilon))=0.
\end{align}
Moreover, $w_\varepsilon$ is a $C^1$ map of $\varepsilon$ in $\mathbb{X}_\varepsilon$ and
\begin{align}\label{ineq:wNorm}
\|w_\varepsilon\|_{L^\infty(\Omega)}+\|w_\varepsilon\|_{\mathbb{X}_\varepsilon}\leq
\left\{
\begin{aligned}
& C\varepsilon^{m^*-2-\frac{\alpha}{2}}(\ln\frac{1}{\varepsilon})^2&\mbox{ if } m^*\leq 3;\\
& C\varepsilon\ln\frac{1}{\varepsilon}&\mbox{ if }m^*_i>3,
\end{aligned}
\right.
\end{align}
where $C$ is a positive constant independent in $\varepsilon$.
\end{proposition}

\noindent{\bf Proof.}
We prove the result in steps.

\vs

\noindent{\bf Step 1. An equivalent problem.}

\vs

By Theorem \ref{t:invertibility}, Problem (\ref{e:projected}) can be rewritten as
\begin{align}\label{e:projection2}
w_\varepsilon=B_\varepsilon(w_\varepsilon):=(Q_\varepsilon\mathbb{L}_\varepsilon)^{-1}Q_\varepsilon(g_{1,\varepsilon}(x,w_\varepsilon),\cdots,g_{n,\varepsilon}(x,w_\varepsilon))
\end{align}
and
\begin{align}\label{ineq:NormofBw}
\|B_\varepsilon w_\varepsilon\|_{L^\infty(\Omega)}+\|B_\varepsilon w_\varepsilon\|_{\mathbb{X}_\varepsilon}\leq C\ln\frac{1}{\varepsilon}\|g_\varepsilon(x,w_\varepsilon)\|_{\mathbb{Y}_\varepsilon}.
\end{align}
Fix a small constant $\beta\in(0,\frac{\alpha}{4})$. Let
\begin{align}\label{SetS1}
\mathcal{S}_\varepsilon=\{w_\varepsilon\in E_\varepsilon|\|w_\varepsilon\|_{L^\infty(\Omega)}+\|w_\varepsilon\|_{\mathbb{X}_\varepsilon}\leq\varepsilon^{1-\beta}\}
\end{align}
if $m^*>3$ and
\begin{align}\label{SetS2}
\mathcal{S}_\varepsilon=\{w_\varepsilon\in E_\varepsilon|\|w_\varepsilon\|_{L^\infty(\Omega)}+\|w_\varepsilon\|_{\mathbb{X}_\varepsilon}\leq\varepsilon^{m^*-2-\beta}\}
\end{align}
for $m^*\leq 3$
We will prove that $B_\varepsilon$ is a contraction map from $\mathcal{S}_\varepsilon$ to $\mathcal{S}_\varepsilon$.

\vs

\noindent{\bf Step 2. $B_\varepsilon$ is a map from $\mathcal{S}_\varepsilon$ to $\mathcal{S}_\varepsilon$.}

\vs

First, we estimate $\|g_\varepsilon(x,w_\varepsilon)\|_{\mathbb{Y}_\varepsilon}$.
Recall
\begin{align}
g_{i,\varepsilon}(x,w)&=
\Delta U_\varepsilon^i-\rho_{i,\varepsilon}+\rho_{i,\varepsilon}\frac{h_i(x)e^{ U_\varepsilon^i}}{\int_\Omega h_i(x)e^{U_\varepsilon^i}}-\Bigg(K_{i,\varepsilon}(x)-\rho_{i,\varepsilon}\frac{h_i(x)e^{ U_\varepsilon^i}}{\int_\Omega h_i(x)e^{ U_\varepsilon^i}}\Bigg)w_{i,\varepsilon}\nonumber\\
&\quad\quad\quad+\frac{K_{i,\varepsilon}(x)}{\rho_i^*}\int_\Omega K_{i,\varepsilon}(x)w_{i,\varepsilon}-\frac{\rho_{i,\varepsilon}h_i(x)e^{ U_{i,\varepsilon}}}{\Big(\int_\Omega h_i(x)e^{ U_{i,\varepsilon}}\Big)^2}\int_\Omega\Big(h_i(x)e^{ U_{i,\varepsilon}}w_{i,\varepsilon}\Big)+N_{i,\varepsilon}\nonumber
\end{align}
and
\begin{align}
N_{i,\varepsilon}&=\rho_{i,\varepsilon}\frac{h_i(x)e^{\ U_{i,\varepsilon}+w_{i,\varepsilon}}}{\int_\Omega h_i(x)e^{ U_{i,\varepsilon}+w_{i,\varepsilon}}} - \rho_{i,\varepsilon}\frac{h_i(x)e^{ U_\varepsilon^i}}{\int_\Omega h_i(x)e^{ U_\varepsilon^i}} - \rho_{i,\varepsilon}\Bigg(\frac{h_i(x)e^{ U_{i,\varepsilon}}}{\int_\Omega h_i(x)e^{ U_{i,\varepsilon}}}\Bigg)w_{i,\varepsilon}\nonumber\\
&\quad\quad\quad+\frac{\rho_{i,\varepsilon}h_i(x)e^{ U_{i,\varepsilon}}}{\Big(\int_\Omega h_i(x)e^{ U_{i,\varepsilon}}\Big)^2}\int_\Omega\Big(h_i(x)e^{ U_{i,\varepsilon}}w_{i,\varepsilon}\Big).\nonumber
\end{align}
By a direct computation, for $x\in\Omega\backslash \cup_{s=1}^N B_{\delta_s}(p_{s,\varepsilon})$, we get
\begin{align}\label{ineq:gYOUT}
g_{i,\varepsilon}(x,w_\varepsilon)=O_{m^*}
\end{align}
for $i=1,\cdots,n$. For any $t=1,\cdots,N$ and any $x\in B_{\delta_t}(p_{t,\varepsilon})$, we get
\begin{align}\label{e:gTOTAL}
g_{i,\varepsilon}(x,w_\varepsilon)&=\Bigg\{\Delta U_\varepsilon^i-\rho_{i,\varepsilon}+ \rho_{i,\varepsilon}\frac{h_i(x)e^{U_{i,\varepsilon}}}{\int_\Omega h_i(x)e^{U_{i,\varepsilon}}}\Bigg\}\nonumber\\
&+\Bigg\{-\Bigg(e^{v_i(\frac{x-p_{t,\varepsilon}}{\varepsilon_t})+2\ln\frac{1}{\varepsilon_t}}
-\rho_{i,\varepsilon}\frac{h_i(x)e^{U_{i,\varepsilon}}}{\int_\Omega h_i(x)e^{U_{i,\varepsilon}}}\Bigg)w_i\Bigg\}
\nonumber\\
&+\Bigg\{ \frac{e^{v_i(\frac{x-p_{t,\varepsilon}}{\varepsilon_t})+2\ln\frac{1}{\varepsilon_t}}}{\rho_i^*}\int_\Omega \Bigg(K_{i,\varepsilon}(x)w_{i,\varepsilon} \Bigg)-\frac{\rho_{i,\varepsilon}h_i(x)e^{U_{i,\varepsilon}}}{(\int_\Omega h_i(x)e^{U_{i,\varepsilon}})^2}\int_\Omega \Bigg(h_i(x)e^{U_{i,\varepsilon}}w_{i,\varepsilon}\Bigg) \Bigg\}
\nonumber\\
&+N_{i,\varepsilon}\nonumber\\
&=:A_{i,1}(x)+A_{i,2}(x)+A_{i,3}(x)+N_{i,\varepsilon}(x).
\end{align}

To begin with, we estimate $A_{i,1}$ in $B_{\delta_t}(p_{t,\varepsilon_t})$. For any $x\in B_{\delta_t}(p_{t,\varepsilon})$, by Lemma  \ref{l:GammaRho} and Lemma \ref{l:Sexpansion}, we get
\begin{align}
A_{i,1}&=\Delta U_\varepsilon^i -\rho_{i,\varepsilon}+\rho_{i,\varepsilon}\frac{h_i(x)e^{U_{i,\varepsilon}}}{\int_\Omega h_i(x)e^{U_{i,\varepsilon}}} \nonumber\\
&=-e^{v_i\big(\frac{x-p_{t,\varepsilon}}{\varepsilon_t}\big)+2\ln\frac{1}{\varepsilon_t}}+2\pi\sum_{t=1}^N \sum_{j=1}^n a^{ij} m_j^*-\rho_{i,\varepsilon}\nonumber\\
&\quad+\rho_{i,\varepsilon}\frac{h_i(x)e^{v_i\big(\frac{x-p_{t,\varepsilon}}{\varepsilon_t}\big)+2\ln\frac{1}{\varepsilon_t}}e^{2\pi m_i^*[ {G^*}(x;p_{t,\varepsilon})- {G^*}(p_{t,\varepsilon};p_{t,\varepsilon})]}}{h_i(p_{t,\varepsilon})\rho_i*}(1+O_{m^*})\nonumber\\
&=O_{m^*} e^{v_i\big(\frac{x-p_{t,\varepsilon}}{\varepsilon_t}\big)+2\ln\frac{1}{\varepsilon_t}}+O_{m^*}+(1+o_\varepsilon(1))e^{v_i\big(\frac{x-p_{t,\varepsilon}}{\varepsilon_t}\big)+2\ln\frac{1}{\varepsilon_t}}\Bigg\{1-\frac{h_i(x)}{h_i(p_{t,\varepsilon})} e^{2\pi m_i^*( {G^*}(x;p_{t,\varepsilon})- {G^*}(p_{t,\varepsilon};p_{t,\varepsilon}))}\Bigg\}.\nonumber
\end{align}
By a direct computation, we get
\begin{align}\label{ineq:A1}
\varepsilon_t^4\|
& A_{i,1}(\varepsilon_t x+\varepsilon_t p_{t,\varepsilon})\rho(x)\|^2_{L^2(B_\frac{\delta}{\varepsilon_t}(0))}=\varepsilon_t^2\|A_{i,1}( x+\varepsilon_t p_{t,\varepsilon_t})\rho(\frac{x}{\varepsilon_t})\|^2_{L^2(B_{\delta}(0))}\nonumber\\
&=O_{m^*}^2\varepsilon^{-\alpha}+\varepsilon^{-2}_t\int_{B_\delta(p_{t,\varepsilon})}e^{2v_i\big(\frac{x-p_{t,\varepsilon}}{\varepsilon_t}\big)}\Bigg\{1-\frac{h_i(x)}{h_i(p_{t,\varepsilon})}e^{2\pi m_i^*( {G^*}(x;p_{t,\varepsilon})- {G^*}(p_{t,\varepsilon};p_{t,\varepsilon}))} \Bigg\}^2\rho(\frac{x}{\varepsilon_t})^2 dx \nonumber\\
&=O_{m^*}^2\varepsilon^{-\alpha}+O\Big(\varepsilon_t^{-2}\int_{B_\delta(p_{t,\varepsilon})}e^{2 v_i\big(\frac{x-p_{t,\varepsilon}}{\varepsilon_t}\big)}|x|^2\rho(\frac{x}{\varepsilon_t})^2 dx\Big)\nonumber\\
&=O_{m^*}^2\varepsilon^{-\alpha}+O(\varepsilon^2).
\end{align}

For $A_{i,2}$, by Lemma \ref{l:Sexpansion}, we get
\begin{align}
|A_{i,2}(x)|\leq(|A_{i,1}(x)|+C\varepsilon^{m^*-2})\|w\|_{L^\infty(\Omega)}.\nonumber
\end{align}
This implies that
\begin{align}\label{ineq:A21}
\varepsilon_t^4 \| A_{i,2}(\varepsilon x+p_{t,\varepsilon})\rho(x)\|^2_{L^2(B_\frac{\delta_t}{\varepsilon_t}(0),dx)}\leq C\varepsilon^2\|w\|^2_{L^\infty(\Omega)}\leq C\varepsilon^{2(m^*-2)-\alpha-\beta}
\end{align}
if $m^*\leq 3$ and
\begin{align}\label{ineq:A22}
\varepsilon_t^4 \| A_{i,2}(\varepsilon x+p_{t,\varepsilon})\rho(x)\|^2_{L^2(B_\frac{\delta_t}{\varepsilon_t}(0))}\leq C\varepsilon^{2-\beta}
\end{align}
if $m^*>3$.
For $A_{i,3}$, we get
\begin{align}
A_{i,3}(x)&=\frac{e^{v_i(\frac{x-p_{t,\varepsilon}}{\varepsilon_t})+2\ln\frac{1}{\varepsilon_t}}}{\rho_i^*}\int_\Omega K_{i,\varepsilon}(x)w_{i,\varepsilon}-\frac{\rho_{i,\varepsilon}h_i(x)e^{U_{i,\varepsilon}}}{\big(\int_\Omega h_i(x)e^{U_{i,\varepsilon}}\big)^2}\int_\Omega h_i(x)e^{U_{i,\varepsilon}}w_{i,\varepsilon}.\nonumber
\end{align}
First, by $w_\varepsilon\in E_\varepsilon$, we get
\begin{align}
\frac{e^{v_i(\frac{x-p_{t,\varepsilon}}{\varepsilon_t})+2\ln\frac{1}{\varepsilon_t}}}{\rho_i^*}\int_\Omega K_{i,\varepsilon}(x)w_{i,\varepsilon}=O_{m^*}\|w_{i,\varepsilon}\|_{L^\infty(\Omega)}e^{v_i(\frac{x-p_{t,\varepsilon}}{\varepsilon_t})+2\ln\frac{1}{\varepsilon_t}}.\nonumber
\end{align}
On the other hand, we get
\begin{align}
&\frac{\rho_{i,\varepsilon}h_i(x)e^{U_{i,\varepsilon}}}{\big(\int_\Omega h_i(x)e^{U_{i,\varepsilon}}\big)^2}\int_\Omega h_i(x)e^{U_{i,\varepsilon}}w_{i,\varepsilon}
=\frac{\rho_{i,\varepsilon}h_i(x)e^{U_{i,\varepsilon}}}{\int_\Omega h_i(x)e^{U_{i,\varepsilon}}}\cdot\frac{\int_\Omega h_i(x)e^{U_{i,\varepsilon}}w_{i,\varepsilon}}{\int_\Omega h_i(x)e^{U_{i,\varepsilon}}}\nonumber\\
&=O(1)e^{v_i\big(\frac{x-p_{t,\varepsilon_t}}{\varepsilon_t}\big)+2\ln\frac{1}{\varepsilon_t}}\frac{\int_{\Omega\backslash\cup_{s=1}^N B_{\delta_s}(p_{s,\varepsilon})}+\sum_{s=1}^N\int_{B_\delta(p_{s,\varepsilon})} h_i(x)e^{U_{i,\varepsilon}}w_{i,\varepsilon}}{\int_\Omega h_i(x)e^{U_{i,\varepsilon}}}\nonumber\\
&=O(1)e^{v_i\big(\frac{x-p_{t,\varepsilon}}{\varepsilon_t}\big)+2\ln\frac{1}{\varepsilon_t}}\Big(O_{m^*}\|w_{i,\varepsilon}\|_{L^\infty}\nonumber\\
&\quad\quad+O(1)\sum_{s=1}^N\int_{B_\delta(p_{s,\varepsilon})}e^{v_i\big(\frac{x-p_{s,\varepsilon}}{\varepsilon_s}\big)+2\ln\frac{1}{\varepsilon_s}}\Big\{\frac{h_i(x)}{h_i(p_{s,\varepsilon})}e^{2\pi m_i^*( {G^*}(x;p_{t,\varepsilon})- {G^*}(p_{t,\varepsilon};p_{t,\varepsilon}))} -1 \Big\}
\|w_{i,\varepsilon}\|_{L^\infty}\Big)\nonumber\\
&=O_{m^*}\|w_{i,\varepsilon}\|e^{v_i\big(\frac{x-p_{t,\varepsilon}}{\varepsilon_t}\big)+2\ln\frac{1}{\varepsilon_t}}.\nonumber
\end{align}
By a direct computation, we get
\begin{equation}\label{ineq:A3}
\varepsilon_t^4\|A_{i,3}(\varepsilon_t x+p_{t,\varepsilon})\rho(x)\|_{L^2(B_\frac{\delta_t}{\varepsilon_t}(0))}=
\left\{
\begin{aligned}
& O_{m^*}^4\varepsilon^{-2\beta}&\mbox{ if } m^*\leq 3;\\
& O_{m^*}^2\varepsilon^{2-2\beta}&\mbox{ if }m^*_i>3.
\end{aligned}
\right.
\end{equation}

For $N_{i,\varepsilon}$, we get
\begin{align}
N_{i,\varepsilon}=\frac{\rho_{i,\varepsilon}h_i(x)e^{U_{i,\varepsilon}}}{\int_\Omega h_i(x)e^{U_{i,\varepsilon}}}\Bigg(-\frac{\int_\Omega h_i(x)e^{U_{i,\varepsilon}}w_{i,\varepsilon}}{\int_\Omega h_i(x)e^{U_{i,\varepsilon}}}w_{i,\varepsilon}+\|w_{i,\varepsilon}\|^2_{L^\infty(\Omega)}\Bigg)\leq  C e^{v_i(\frac{x-p_{t,\varepsilon_t}}{\varepsilon_t})+2\ln\frac{1}{\varepsilon_t}}\|w_\varepsilon\|^2_{L^\infty(\Omega)}.\nonumber
\end{align}
Therefore,

\begin{align}\label{ineq:A4}
\varepsilon_t^4\|N_{i,\varepsilon}(\varepsilon x+p_\varepsilon)\rho(x)\|^2_{L^2(B_\frac{\delta}{\varepsilon}(0),dx)}
=\left\{
\begin{aligned}
& O_{m^*}^4\varepsilon^{-4\beta}&\mbox{ if } m^*\leq 3;\\
& \varepsilon^{4-4\beta}&\mbox{ if }m^*_i>3.
\end{aligned}
\right.
\end{align}

Combining (\ref{ineq:gYOUT}), (\ref{ineq:A1}), (\ref{ineq:A21}), (\ref{ineq:A22}), (\ref{ineq:A3}) and (\ref{ineq:A4}), we get $\|g_{i,\varepsilon}\|_{\mathbb{Y}_\varepsilon}=O_{m^*}$.
The above computation gives
\begin{align}\label{ineq:BwNorm}
\|B_\varepsilon w_\varepsilon\|_{L^\infty(\Omega)} +\|B_\varepsilon w_\varepsilon\|_{\mathbb{X}_\varepsilon}\leq C\ln\frac{1}{\varepsilon}\|g_\varepsilon(x,w_\varepsilon)\|_{\mathbb{Y}_\varepsilon} \leq
\left\{
\begin{aligned}
& C\varepsilon^{m^*-2-\frac{\alpha}{2}}(\ln\frac{1}{\varepsilon})^2&\mbox{ if } m^*\leq 3;\\
& C\varepsilon\ln\frac{1}{\varepsilon}&\mbox{ if }m^*_i>3.
\end{aligned}
\right.
\end{align}

~

\noindent{\bf Step 3. $B_\varepsilon$ is a contraction map from $\mathcal{S}_\varepsilon$ to $\mathcal{S}_\varepsilon$.}

~

For any $w=(w_1,\cdots,w_n)\in\mathcal{S}_\varepsilon$ and any $v=(v_1,\cdots,v_n)\in \mathcal{S}_\varepsilon$, we get
\begin{align}
\|B_\varepsilon w-B_\varepsilon v\|_{L^\infty(\Omega)}+\|B_\varepsilon w-B_\varepsilon v\|_{\mathbb{X}_\varepsilon}
\leq C\ln\frac{1}{\varepsilon}
\|g_{i,\varepsilon}(x,w) - g_{i,\varepsilon}(x,v)\|_{\mathbb{Y}_\varepsilon}\nonumber
\end{align}
Now we estimate $\|g_{i,\varepsilon}(x,w) - g_{i,\varepsilon}(x,v)\|_{\mathbb{Y}_\varepsilon}$. To begin with, observe that
\begin{align}
g_{i,\varepsilon}(x,w) - g_{i,\varepsilon}(x,v) &=-\Bigg(K_{i,\varepsilon}(x)-\rho_{i,\varepsilon}\frac{h_i(x)e^{ U_\varepsilon^i}}{\int_\Omega h_i(x)e^{ U_\varepsilon^i}}\Bigg)(w_{i}-v_i) +\frac{K_{i,\varepsilon}(x)}{\rho_i^*}\int_\Omega K_{i,\varepsilon}(x)(w_{i}-v_i)\nonumber\\
&\quad\quad\quad-\frac{\rho_{i,\varepsilon}h_i(x)e^{ U_{i,\varepsilon}}}{\Big(\int_\Omega h_i(x)e^{ U_{i,\varepsilon}}\Big)^2}\int_\Omega\Big(h_i(x)e^{ U_{1,\varepsilon}}(w_{i}-v_i)\Big)+(N_{i,\varepsilon}(x,w)-N_{i,\varepsilon}(x,v)).\nonumber
\end{align}
A similar idea as in (\ref{ineq:A21}), (\ref{ineq:A22}) and (\ref{ineq:A3}) gives
\begin{align}\label{Ineq:Contraction1}
\Bigg\|&-\Bigg(K_{i,\varepsilon}(x)-\rho_{i,\varepsilon}\frac{h_i(x)e^{ U_\varepsilon^i}}{\int_\Omega h_i(x)e^{ U_\varepsilon^i}}\Bigg)(w_{i}-v_i) +\frac{K_{i,\varepsilon}(x)}{\rho_i^*}\int_\Omega K_{i,\varepsilon}(x)(w_{i}-v_i)\nonumber\\
&\quad\quad\quad-\frac{\rho_{i,\varepsilon}h_i(x)e^{ U_{i,\varepsilon}}}{\Big(\int_\Omega h_i(x)e^{ U_{i,\varepsilon}}\Big)^2}\int_\Omega\Big(h_i(x)e^{ U_{1,\varepsilon}}(w_{i}-v_i)\Big)\Bigg\|_{\mathbb{Y}_\varepsilon}= O_{m^*}\|w-v\|_{L^\infty}.
\end{align}
For $N_{i,\varepsilon}(x,w)-N_{i,\varepsilon}(x,v)$, we find
\begin{align}
|N_{i,\varepsilon}(x,w)-N_{i,\varepsilon}(x,v)|\leq Ce^{U_{i,\varepsilon}(x)}(\|w\|_{L^\infty(\Omega)}+\|v\|_{L^\infty(\Omega
)})\|w-v\|_{L^\infty(\Omega)}.\nonumber
\end{align}
Combining (\ref{Ineq:Contraction1}), (\ref{SetS1}) and (\ref{SetS2}), we get
\begin{align}
\|g_{i,\varepsilon}(x,w) - g_{i,\varepsilon}(x,v)\|_{\mathbb{Y}_\varepsilon}\leq C(O_{m^*}\varepsilon^{-\beta}+\varepsilon^{1-\beta})\ln\frac{1}{\varepsilon}\|w-v\|_{L^\infty(\Omega)}.\nonumber
\end{align}
This implies that
\begin{align}
\|B_\varepsilon w-B_\varepsilon v\|_{L^\infty(\Omega)}+\|B_\varepsilon w-B_\varepsilon v\|_{\mathbb{X}_\varepsilon}
\leq \frac{1}{2}\big(\|w-v\|_{L^\infty(\Omega)}+\|w-v\|_{\mathbb{X}_\varepsilon}\big).\nonumber
\end{align}
Therefore, $B_\varepsilon$ is a map from $\mathcal{S}_\varepsilon$ to $\mathcal{S}_\varepsilon$. By Banach contraction mapping theorem, we find a unique $w_\varepsilon\in \mathcal{S}_\varepsilon$ with
$w_\varepsilon=B_\varepsilon w_\varepsilon$.
Combining (\ref{ineq:BwNorm}), we completes the proof.
$\Box$

\subsection{Refined Estimates on the Error Term I: Parts of $w_\varepsilon$ with zero frequency}\label{subsection:RefindeReduction1}

In this part, we provide a refined estimate on $w_\varepsilon$ in the view point of Fourier analysis  (See \cite{LinZhang2013}). This method gives finer results by splitting the error term into parts with different frequencies and analysing the detailed properties of them. To begin with, let us recall the rescaled problem
\begin{align}\label{e:ScaledEquation}
\widetilde{L}_{i,\varepsilon}\widetilde{w}_\varepsilon&=\varepsilon_t^2\widetilde{g}_{i,\varepsilon}+\sum_{j=1}^3\sum_{t=1}^N r_{j,t}\widetilde{Z}_{j,t,i,\varepsilon}^*\cdot \sum_{l=1}^N\widetilde{\chi}_l e^{v_i(y+\frac{p_{t,\varepsilon}-p_{l,\varepsilon}}{\varepsilon_t})+2\ln\frac{\varepsilon_t}{\varepsilon_l}}+ s_i\sum_{l=1}^N\widetilde{\chi}_l e^{v_i(y+\frac{p_{t,\varepsilon}-p_{l,\varepsilon}}{\varepsilon_t})+2\ln\frac{\varepsilon_t}{\varepsilon_l}}
\end{align}
with
\begin{align}
\widetilde{w}_{t,i,\varepsilon}(y)&:=w_{i,\varepsilon}(\varepsilon_t y+p_{t,\varepsilon}),\nonumber\\
\widetilde{w}^i_{t,\varepsilon}(y)&:=\sum_{j=1}^n a^{ij}w_{j,\varepsilon}(\varepsilon_t y+p_{t,\varepsilon}),\nonumber\\
\widetilde{w}_{t,\varepsilon}(y)&:=(\widetilde{w}_{t,1,\varepsilon}(y), \cdots,\widetilde{w}_{t,n,\varepsilon}(y)),\nonumber\\
\widetilde{g}_{t,i,\varepsilon}(y)&:=g_{i,\varepsilon}(\varepsilon_t y+p_{t,\varepsilon}),\nonumber\\
\widetilde{K}_{t,i,\varepsilon}(y)&:=K_{i,\varepsilon}(\varepsilon_t y+p_{t,\varepsilon}),\nonumber\\
\widetilde{\chi}_{t,s}(y)&:=\chi_s(\varepsilon_t y+p_{t,\varepsilon}),\nonumber\\
\widetilde{Z}_{t,j,s,\varepsilon}^*&:= Z^*_{j,s,\varepsilon}(\varepsilon_t y+p_{t,\varepsilon}),\nonumber\\
\Omega_{t,\varepsilon}&:=\{y|\varepsilon_t y+p_{t,\varepsilon}\in\Omega\}\nonumber
\end{align}
and
\begin{align}
\widetilde{L}_{t,i,\varepsilon}\widetilde{w}_{t,\varepsilon}:=-\Delta_y\Big(\sum_{l=1}^n a^{il}\widetilde{w}_{t,l,\varepsilon}\Big)-\varepsilon_t^2\widetilde{K}_{t,i,\varepsilon}\widetilde{w}_{t,i,\varepsilon}+\frac{\varepsilon_t^2\widetilde{K}_{t,i,\varepsilon}}{\rho_i^*}\int_{\Omega_{t,\varepsilon}} \widetilde{K}_{t,i,\varepsilon}(y)\widetilde{w}_{t,i,\varepsilon}(y)dy.\nonumber
\end{align}
We study the 0-frequency part of $\widetilde{w}_\varepsilon$.
Let us begin with the average of $\widetilde{w}_{t,i,\varepsilon}$ on the circles
\begin{align}
R^0_{t,i,\varepsilon}(r)=\frac{1}{2\pi r}\int_{\partial B_r(p_{t,\varepsilon_t})}\widetilde{w}_{t,i,\varepsilon}(x)dS_x,\nonumber
\end{align}
which is governed by the following ODE
\begin{align}\label{e:0Frequency}
-\Big(\frac{d^2}{dr^2}+\frac{1}{r}\frac{d}{dr}\Big)\Big(\sum_{l=1}^n a^{il} R_{t,l,\varepsilon}^0(r)\Big)-e^{v_i(r)}R_{t,i,\varepsilon}^0(r)&=\varepsilon_t^2\widetilde{g}_{i,\varepsilon}^0+\sum_{t=1}^N r_{3,t}\widetilde{Z}_{3,t,i,\varepsilon}^*\cdot \sum_{l=1}^N\widetilde{\chi}_l e^{v_i(y+\frac{p_{t,\varepsilon}-p_{l,\varepsilon}}{\varepsilon_t})+2\ln\frac{\varepsilon_t}{\varepsilon_l}}\nonumber\\
&\quad+ s_i\sum_{l=1}^N\widetilde{\chi}_l e^{v_i(y+\frac{p_{t,\varepsilon}-p_{l,\varepsilon}}{\varepsilon_t})+2\ln\frac{\varepsilon_t}{\varepsilon_l}}.
\end{align}
Here,
\begin{align}
\widetilde{g}_{i,\varepsilon}^0(r)=\frac{1}{2\pi r}\int_{\partial B_r(p_{t,\varepsilon_t})}\widetilde{g}_{i,\varepsilon}(x)dS_x.\nonumber
\end{align}
By a similar approach as in Proposition \ref{prop:Recution}, we find
\begin{align}
g_{i,\varepsilon}(x)&=O_{m^*}\ln\frac{1}{\varepsilon} e^{v_i\big(\frac{x-p_{t,\varepsilon}}{\varepsilon_t}\big)+2\ln\frac{1}{\varepsilon_t}}+O_{m^*}\ln\frac{1}{\varepsilon_t}\nonumber\\
&\quad+(1+o_\varepsilon(1))e^{v_i\big(\frac{x-p_{t,\varepsilon}}{\varepsilon_t}\big)+2\ln\frac{1}{\varepsilon_t}}\times\nonumber\\
&\quad\quad\times\Bigg\{1-\frac{h_i(x)}{h_i(p_{t,\varepsilon})} e^{2\pi m_i^*[\gamma(x,p_{t,\varepsilon})-\gamma(p_{t,\varepsilon},p_{t,\varepsilon})] +2\pi m_i^*\sum_{s\neq t} [G(x,p_{t,\varepsilon})-G(p_{t,\varepsilon},p_{t,\varepsilon})]} \Bigg\}\nonumber\\
&\quad+\Bigg(O_{m^*}\ln\frac{1}{\varepsilon} e^{v_i\big(\frac{x-p_{t,\varepsilon}}{\varepsilon_t}\big)+2\ln\frac{1}{\varepsilon_t}}+O_{m^*}\ln\frac{1}{\varepsilon_t}+O(|x|)\Bigg)|w_{i,\varepsilon}(x)|\nonumber\\
&\quad+O_{m^*}\|w_{i,\varepsilon}\|_{L^\infty(\Omega)}e^{v_i\big(\frac{x-p_{t,\varepsilon}}{\varepsilon_t}\big)+2\ln\frac{1}{\varepsilon_t}}+O\Big(e^{v_i\big(\frac{x-p_{t,\varepsilon}}{\varepsilon_t}\big)+2\ln\frac{1}{\varepsilon_t}}\|w_{i,\varepsilon}\|_{L^\infty(\Omega)}^2\Big).
\end{align}
This implies that
$\|g_{i,\varepsilon}^0\|_{\mathbb{Y}_\varepsilon}=\varepsilon^{-\alpha}O_{m^*}$
with
$$g_{i,\varepsilon}^0(r)=\frac{1}{2\pi r}\int_{\partial B_r(p_{t,\varepsilon_t})}g_{i,\varepsilon} (x)d S_x.$$


\begin{proposition}\label{prop:w0Pointwise}
There exists constants $\tau,C>0$ such that
\begin{align}\label{Ineq:R0Pointwise}
|R^0_{t,i,\varepsilon}(r)|\leq \varepsilon^{-\alpha} O_{m^*}(1+r)^\tau.
\end{align}
for $|r|< \frac{3\delta_t}{\varepsilon}$. Here, $\alpha$ is a small positive constant.
\end{proposition}

We apply a method \`a la \cite{LinZhang2013,BartolucciYangZhang2024}. A similar computation can be found in Claim \ref{c:wlndecay}.

\noindent{\bf Proof.}
Denote $w^0_{t,i,\varepsilon}(x)=R^0_{t,i,\varepsilon}(|x|)$ and it satisfies the equation
Recall that
\begin{align}
L_{i,\varepsilon}w^0_\varepsilon&=g^0_{i,\varepsilon}+\sum_{t=1}^N r_{3,t}Z_{3,t,i,\varepsilon}^*\cdot \sum_{s=1}^N\chi_s e^{v_i(\frac{x-p_{s,\varepsilon}}{\varepsilon_s})+2\ln\frac{1}{\varepsilon_s}}+ s_i\sum_{s=1}^N\chi_s e^{v_i(\frac{x-p_{s,\varepsilon}}{\varepsilon_s})+2\ln\frac{1}{\varepsilon_s}}.\nonumber
\end{align}

\begin{claim}\label{c:Coefficients}
Under the above assumptions, we get
\begin{align}
r_{3,t}&=O(\frac{1}{\varepsilon}\|g^0_{i,\varepsilon}\|_{\mathbb{Y}_\varepsilon})+O(\varepsilon^{m^*-3}\|w^0_\varepsilon\|_{L^\infty(B_{3\delta_t}(p_{t,\varepsilon}))}),\nonumber\\
s_i&=O(\|g^0_{i,\varepsilon}\|_{\mathbb{Y}_\varepsilon})+O(\varepsilon^{m^*-2}\|w^0_\varepsilon\|_{L^\infty(B_{3\delta_t}(p_{t,\varepsilon}))}).\nonumber
\end{align}
\end{claim}

We omit the proof since we will prove a more general result in Lemma \ref{l:constantsrrs}.
Notice that
\begin{align}\label{Ineq:gwSize}
\|g^0_{i,\varepsilon}\|_{\mathbb{Y}_\varepsilon}=\varepsilon^{-\alpha}O_{m^*},\quad\|w^0_\varepsilon\|_{L^\infty(B_{3\delta_t}(p_{t,\varepsilon}))}\leq \|w_\varepsilon\|_{L^\infty(\Omega)}\leq
\left\{
\begin{aligned}
& C\varepsilon^{m^*-2-\frac{\alpha}{2}}(\ln\frac{1}{\varepsilon})^2&\mbox{ if } m^*\leq 3;\\
& C\varepsilon\ln\frac{1}{\varepsilon}&\mbox{ if }m^*_i>3,
\end{aligned}
\right.
\end{align}
This implies that
\begin{align}\label{Ineq:rsSize}
\varepsilon r_{3,t},\,s_i=\varepsilon^{-\alpha}O_{m^*}
\end{align}
for $t=1,\cdots,N$ and $i=1,\cdots,n$. And $\alpha$ is a generically small positive number.

To prove Proposition \ref{prop:w0Pointwise} is to prove
\begin{align}\label{Ineq:w0Pointwise}
|\widetilde{w}^0_{t,i,\varepsilon}(y)|\leq \varepsilon^{-\alpha}O_{m^*}(1+|y|)^\tau
\end{align}
for $y\in B_{\frac{3\delta_t}{\varepsilon_t}}(p_{t,\varepsilon})$. We argue by contradiction. Assume that
\begin{align}
\Lambda_\varepsilon:=\max_{y\in B_{\frac{3\delta_t}{\varepsilon_t}}(p_{t,\varepsilon})}\frac{|\widetilde{w}^0_{t,i,\varepsilon}(y)|}{\varepsilon^{-\alpha}O_{m^*}(1+|y|)^\tau}\to+\infty.\nonumber
\end{align}
Letting $y_\varepsilon\in B_{\frac{3\delta}{\varepsilon_t}}(p_{t,\varepsilon_t})$ be the point where the function $\frac{|\widetilde{w}_{t,i,\varepsilon}(y)|}{\varepsilon(1+|y|)^\tau}$ achieves its maximum, then we define
\begin{align}
\widetilde{w}^*_{t,i,\varepsilon}(y)= \frac{\widetilde{w}^0_{t,i,\varepsilon}(y)}{\Lambda_\varepsilon\cdot\varepsilon^{-\alpha}O_{m^*}(1+|y_\varepsilon|)^\tau}\nonumber
\end{align}
and
\begin{align}
\widetilde{L}_{t,i,\varepsilon}\widetilde{w}^*_{t,\varepsilon}:=-\Delta_y\Big(\sum_{l=1}^n a^{il}\widetilde{w}^*_{t,l,\varepsilon}\Big)-\varepsilon_t^2\widetilde{K}_{t,i,\varepsilon}\widetilde{w}^*_{t,i,\varepsilon}+\frac{\varepsilon_t^2\widetilde{K}_{t,i,\varepsilon}}{\rho_i^*}\int_{\Omega_{t,\varepsilon}} \widetilde{K}_{t,i,\varepsilon}(y)\widetilde{w}^*_{t,i,\varepsilon}(y)dy.\nonumber
\end{align}
An immediate observation is that
\begin{align}\label{Ineq:wGrowth}
|\widetilde{w}^*_{t,i,\varepsilon}(y)|\leq \Big|\frac{\widetilde{w}^0_{t,i,\varepsilon}(y)}{\Lambda_\varepsilon\cdot\varepsilon^{-\alpha}O_{m^*} (1+|y|)^\tau}\cdot\frac{(1+|y|)^\tau}{(1+|y_\varepsilon|)^\tau}\Big|\leq\frac{(1+|y|)^\tau}{ (1+|y_\varepsilon|)^\tau}.
\end{align}

Then, there holds the equation
\begin{align}\label{e:ScaledExpansion}
\widetilde{L}_{i,\varepsilon}\widetilde{w}^*_\varepsilon&=\frac{\varepsilon_t^2\widetilde{g}^0_{i,\varepsilon}}{\Lambda_\varepsilon\cdot\varepsilon^{-\alpha}O_{m^*}(1+|y_\varepsilon|)^\tau}+\sum_{t=1}^N \frac{r_{3,t}\widetilde{Z}_{3,t,i,\varepsilon}^*}{\Lambda_\varepsilon\cdot\varepsilon^{-\alpha}O_{m^*}(1+|y_\varepsilon|)^\tau}\cdot \sum_{s=1}^N\widetilde{\chi}_s e^{v_i(y+\frac{p_{t,\varepsilon}-p_{l,\varepsilon}}{\varepsilon_t})+2\ln\frac{\varepsilon_t}{\varepsilon_l}}\nonumber\\
&\quad+ s_i\sum_{s=1}^N\frac{\widetilde{\chi}_s e^{v_i(y+\frac{p_{t,\varepsilon}-p_{l,\varepsilon}}{\varepsilon_t})+2\ln\frac{\varepsilon_t}{\varepsilon_l}}}{\Lambda_\varepsilon\cdot\varepsilon^{-\alpha}O_{m^*}(1+|y_\varepsilon|)^\tau}.
\end{align}

Now we apply a skill \`a la \cite{LinZhang2013,BartolucciYangZhang2024}. To do this, we divide the argument into three parts.

~~

\noindent{\bf Case 1. $|y_\varepsilon|$ is bounded.}

~~

By a direct computation, for any $R>0$ and any $p\in(1,2)$, (\ref{Ineq:gwSize}) implies that
\begin{align}\label{ineq:gLp1}
\Big\|\frac{\varepsilon_t^2 \widetilde{g}_{i,\varepsilon}}{\Lambda_\varepsilon\cdot\varepsilon^{-\alpha}O_{m^*}(1+|y_\varepsilon|)^\tau}\Big\|_{L^p(B_R(0))} \leq\frac{C}{\Lambda_\varepsilon\cdot\varepsilon^{-\alpha}O_{m^*}(1+|y_\varepsilon|)^\tau}\|g_{i,\varepsilon}\|_{\mathbb{Y}_\varepsilon}\leq \frac{C}{\Lambda_\varepsilon\cdot(1+|y_\varepsilon|)^\tau}=o_\varepsilon(1)
\end{align}
and
\begin{align}\label{ineq:Z3Lp1}
\Big\|&\sum_{t=1}^N \frac{r_{3,t}\widetilde{Z}_{3,t,i,\varepsilon}^*}{\Lambda_\varepsilon\cdot\varepsilon^{-\alpha}O_{m^*} (1+|y_\varepsilon|)^\tau}\cdot \sum_{s=1}^N\widetilde{\chi}_s e^{v_i(y+\frac{p_{t,\varepsilon}-p_{l,\varepsilon}}{\varepsilon_t})+ 2\ln\frac{\varepsilon_t}{\varepsilon_l}}\Big\|_{L^p(B_R(0))}\nonumber\\
&\leq\frac{C \varepsilon r_{t,3}}{\Lambda_\varepsilon\cdot\varepsilon^{-\alpha}O_{m^*}(1+|y_\varepsilon|)^\tau}\leq \frac{C}{\Lambda_\varepsilon(1+|y_\varepsilon|)^\tau}=o_\varepsilon(1)
\end{align}
and
\begin{align}\label{ineq:esLp1}
\Big\|s_i\sum_{s=1}^N\frac{\widetilde{\chi}_s e^{v_i(y+\frac{p_{t,\varepsilon}-p_{l,\varepsilon}}{\varepsilon_t})+ 2\ln\frac{\varepsilon_t}{\varepsilon_l}}}{\Lambda_\varepsilon\cdot\varepsilon^{-\alpha}O_{m^*} (1+|y_\varepsilon|)^\tau} \Big\|_{L^p(B_R(0))}\leq\frac{C s_i}{\Lambda_\varepsilon\cdot\varepsilon^{-\alpha}O_{m^*}(1+|y_\varepsilon|)^\tau}\leq \frac{C}{\Lambda_\varepsilon(1+|y_\varepsilon|)^\tau}=o_\varepsilon(1).
\end{align}
Here, we use Claim \ref{c:Coefficients}.
Then there exists a radial $C^2(\mathbb{R}^2)$ function $\widetilde{w}_{t,i}^*$ such that for some $\gamma\in(0,1)$
\begin{align}
\widetilde{w}^*_{t,i,\varepsilon}(y)\to \widetilde{w}_{t,i}^*\mbox{ in }C^\gamma(\mathbb{R}^2)\nonumber
\end{align}
with
\begin{equation}
\left\{
\begin{array}{lr}
-\Delta\Big(\sum_{j=1}^n a^{ij}\widetilde{w}^*_{t,j}\Big)-e^{v_i(y)}\widetilde{w}^*_{t,i}=0\mbox{ in }\mathbb{R}^2,\\
|w_{t,i}^*(y)|\leq C(1+|y|)^\tau;\\
\sum_{i=1}^n\int_{\mathbb{R}^2}e^{v_i(y)}(y\cdot\nabla v_i(y)+2)\widetilde{w}_{t,i}^*(y)dy=0\nonumber\\
w_{t,i}^*(y_0)=1\mbox{ with }y_0=\lim_{\varepsilon\to0}y_\varepsilon.
\end{array}
\right.
\end{equation}
This contradicts with Corollary \ref{coro:KernelLocal}.

~

\noindent{\bf Case 2. $|y_\varepsilon|$ is unbounded and $|y_\varepsilon|=o_\varepsilon(1)\varepsilon^{-1}$.}

~

By a similar computation as in {\bf Case 1}, we know that $\hat{\hat{w}}_{t,i,\varepsilon}(0)\to0$ as $\varepsilon\to0$. On the other hand, we always have
$\widetilde{w}^*_{t,i,\varepsilon}(y_\varepsilon)=\pm 1$.
Then, by Green's formula, we get
\begin{align}\label{Ineq:GREEN}
\frac{1}{2}&\leq|\widetilde{w}^*_{t,i,\varepsilon}(y_\varepsilon)-\widetilde{w}^*_{t,i,\varepsilon}(0)|\nonumber\\
&\leq\Bigg\{\int_{B_{\frac{2\delta}{\varepsilon}}(0)}\Big|G_\varepsilon(y_\varepsilon,\eta) -G_\varepsilon(0,\eta)\Big|\cdot\Big|\sum_{j=1}^n a_{ij}e^{v_j(\eta)}\widetilde{w}^*_{t,i,\varepsilon}(\eta)\Big|d\eta\Bigg\}\nonumber\\
&+\Bigg\{ \int_{B_{\frac{2\delta}{\varepsilon}}(0)}\Big|G_\varepsilon(y_\varepsilon,\eta) -G_\varepsilon(0,\eta)\Big|\times\nonumber\\
&\times\Big|\frac{\varepsilon_t^2\widetilde{g}_{i,\varepsilon}(\eta)}{\Lambda_\varepsilon\cdot \varepsilon^{-\alpha}O_{m^*}(1+|y_\varepsilon|)^\tau} +\sum_{t=1}^N \frac{r_{3,t}\widetilde{Z}_{3,t,i,\varepsilon}^*(\eta)}{\Lambda_\varepsilon\cdot\varepsilon^{-\alpha}O_{m^*} (1+|y_\varepsilon|)^\tau}\cdot \sum_{s=1}^N\widetilde{\chi}_s (\eta) e^{v_i(\eta+\frac{p_{t,\varepsilon}-p_{l,\varepsilon}}{\varepsilon_t})+ 2\ln\frac{\varepsilon_t}{\varepsilon_l}}\nonumber\\
&\quad +s_i\sum_{s=1}^N\frac{\widetilde{\chi}_s(\eta) e^{v_i(\eta+\frac{p_{t,\varepsilon}-p_{l,\varepsilon}}{\varepsilon_t})+ 2\ln\frac{\varepsilon_t}{\varepsilon_l}}}{\Lambda_\varepsilon\cdot\varepsilon^{-\alpha}O_{m^*}(1+|y_\varepsilon|)^\tau}\Big|d\eta\Bigg\}\nonumber\\
&=:A+B.
\end{align}
Here, $G_\varepsilon(y,\eta)$ denotes the Green's function with Dirichlet boundary condition on $B_{\frac{2\delta_t}{\varepsilon_t}}(0)$. The boundary terms are canceled out since $\widetilde{w}^*_{t,i,\varepsilon}$ is radial.

We will draw a contradiction by proving that $A+B=o_\varepsilon(1)$.

By (\ref{Ineq:wGrowth}), we have the pointwise estimate
\begin{align}\label{Ineq:Pointwise0}
\Big|e^{v_j(\eta)}\widetilde{w}^*_{t,j,\varepsilon}(\eta)\Big|\leq\frac{C}{(1+|y_\varepsilon|)^\tau(1+|\eta|)^{m^*-\tau}}.
\end{align}

For $p\in(1,2)$, (\ref{Ineq:gwSize}) implies that
\begin{align}\label{Ineq:Lp010}
\Bigg(\int_{B_{\frac{2\delta}{\varepsilon}}(0)} \Bigg|\frac{\varepsilon_t^2\widetilde{g}_{i,\varepsilon}(\eta)}{\Lambda_\varepsilon\cdot\varepsilon^{-\alpha}O_{m^*} (1+|y_\varepsilon|)^\tau}\Bigg|^pd\eta
\Bigg)^\frac{1}{p}
\leq \frac{C\|g_{i,\varepsilon}\|_{\mathbb{Y}_\varepsilon}}{\Lambda_\varepsilon\cdot \varepsilon^{-\alpha}O_{m^*} (1+|y_\varepsilon|)^\tau}=\frac{C}{\Lambda_\varepsilon(1+|y_\varepsilon|)^\tau}.
\end{align}
Moreover, we get
\begin{align}\label{Ineq:Lp020}
\Bigg(\int_{B_\frac{2\delta}{\varepsilon}(0)}\Bigg|\frac{r_{3,t}\widetilde{Z}_{3,t,i,\varepsilon}^*}{\Lambda_\varepsilon \cdot\varepsilon^{-\alpha}O_{m^*}(1+|y_\varepsilon|)^\tau}\cdot \sum_{s=1}^N\widetilde{\chi}_s e^{v_i(y+\frac{p_{t,\varepsilon}-p_{l,\varepsilon}}{\varepsilon_t})+2\ln\frac{\varepsilon_t}{\varepsilon_l}}\Bigg|^p\Bigg)^\frac{1}{p}\leq \frac{C}{\Lambda_\varepsilon(1+|y_\varepsilon|)^\tau}.
\end{align}
and
\begin{align}\label{Ineq:Lp030}
\Bigg(\int_{B_\frac{2\delta}{\varepsilon}(0)}\Bigg|s_i\sum_{s=1}^N\frac{\widetilde{\chi}_s(\eta) e^{v_i(\eta+\frac{p_{t,\varepsilon}-p_{l,\varepsilon}}{\varepsilon_t}) +2\ln\frac{\varepsilon_t}{\varepsilon_l}}}{\Lambda_\varepsilon\cdot\varepsilon^{-\alpha}O_{m^*} (1+|y_\varepsilon|)^\tau}\Bigg|^p d\eta\Bigg)^\frac{1}{p}=\frac{o_\varepsilon(1)}{\Lambda_\varepsilon(1+|y_\varepsilon|)^\tau}.
\end{align}
Here, we use (\ref{Ineq:rsSize}).

In order to estimate $A+B$, we apply the methods in \cite[Theorem 3.2]{LinZhang2013} (see also \cite[Proposition 3.1]{BartolucciYangZhang2024}). By \cite[Lemma 3.2]{LinZhang2013}, we get
\begin{equation}\label{Ineq:GreenEstimates0}
|G_\varepsilon(y_\varepsilon,
    \eta)-G_\varepsilon(0,\eta)| \leq \left\{
    \begin{aligned}
    & C(\ln|y_\varepsilon| +|\ln|\eta||)&\mbox{ if }\eta\in\Sigma_1:=\{\eta\in B_{\frac{2\delta_t}{\varepsilon_t}}(0)||\eta|<\frac{|y_\varepsilon|}{2}\}; \\
&C(\ln|y_\varepsilon|+|\ln|y-\eta||)&\mbox{ if }\eta\in\Sigma_2:=\{\eta\in B_{\frac{2\delta_t}{\varepsilon_t}}(0)||y_\varepsilon-\eta|<\frac{|y_\varepsilon|}{2}\};\\
    &\frac{C|y_\varepsilon|}{|\eta|}&\mbox{ if }\eta\in\Sigma_3:= B_{\frac{2\delta_t}{\varepsilon_t}(0)}\backslash(\Sigma_1\cup\Sigma_2).
    \end{aligned}
    \right.
\end{equation}
In follows, we estimate $A+B$ with the help of (\ref{Ineq:Pointwise0}), (\ref{Ineq:Lp010}), (\ref{Ineq:Lp020}), (\ref{Ineq:Lp030}) and (\ref{Ineq:GreenEstimates0}).

As for $A$, we get
\begin{align}\label{ineq:ATOTAL1}
A&\leq C\int_{B_\frac{2\delta}{\varepsilon}(0)}\frac{|G_\varepsilon(y_\varepsilon,\eta)-G_\varepsilon(0,\eta)|d\eta}{(1+|y_\varepsilon|)^\tau(1+|\eta|)^{m^*-\tau}}\nonumber\\
&=C\Bigg(\int_{\Sigma_1}+\int_{\Sigma_2}+\int_{\Sigma_3}\Bigg)\frac{|G_\varepsilon(y_\varepsilon,\eta)-G_\varepsilon(0,\eta)|d\eta}{(1+|y_\varepsilon|)^\tau(1+|\eta|)^{m^*-\tau}}=:C(A_1+A_2+A_3).
\end{align}
Here,
\begin{align}\label{ineq:A11}
A_1&\leq C\int_{|\eta|<\frac{y_\varepsilon}{2}}\frac{(\ln|y_\varepsilon|+|\ln|\eta||)d\eta}{(1+|y_\varepsilon|)^\tau(1+|\eta|)^{m^*-\tau}}\nonumber\\
&\leq C\int_0^\frac{|y_\varepsilon|}{2}\frac{rdr}{(1+r)^{m^*-\tau}}\cdot\frac{\ln|y_\varepsilon|}{(1+|y_\varepsilon|)^\tau}+C\int_0^\frac{|y_\varepsilon|}{2}\frac{|\ln r|rdr}{(1+r)^{m^*-\tau}}\frac{1}{(1+|y_\varepsilon|)^\tau}=o_\varepsilon(1)
\end{align}
and
\begin{align*}
A_2&\leq C\int_{|\eta|<\frac{y_\varepsilon}{2}}\frac{(\ln|y_\varepsilon|+|\ln|\eta-y_\varepsilon||)d\eta}{(1+|y_\varepsilon|)^\tau(1+|\eta|)^{m^*-\tau}}\nonumber\\
&=o_\varepsilon(1)+ C\Big\{\int_{|\eta-y_\varepsilon|<1,|\eta|<\frac{y_\varepsilon}{2}}+\int_{|\eta-y_\varepsilon|>1}\Big\}\frac{|\ln|\eta-y_\varepsilon||d\eta}{(1+|y_\varepsilon|)^\tau(1+|\eta|)^{m^*-\tau}}\nonumber\\
&\leq o_\varepsilon(1)+\frac{C}{(1+|y_\varepsilon|)^\tau}\int_{|\eta|<1}|\ln|\eta||d\eta+\frac{\ln(2|y_\varepsilon|)}{(1+|y_\varepsilon|)^\tau}\int_{\mathbb{R}^2}\frac{d\eta}{(1+|\eta|)^{m^*-\tau}}=o_\varepsilon(1),
\end{align*}
which can be deduced by a similar method. Moreover,
\begin{align}\label{ineq:A31}
A_3=C\int_{\Sigma_3}\frac{d\eta}{(1+|\eta|)^{m^*-\tau}|\eta|}\cdot\frac{|y_\varepsilon|}{(1+|y_\varepsilon|)^\tau}=
C\int_{\Sigma_3}\frac{d\eta}{(1+|\eta|)^{m^*-\tau}}\cdot\frac{1}{(1+|y_\varepsilon|)^\tau}=o_\varepsilon(1).
\end{align}
Together, they give
\begin{align}\label{Ineq:A}
A=o_\varepsilon(1).
\end{align}
On the other hand, we get
\begin{align}\label{ineq:BTOTAL1}
B&\leq \frac{o_\varepsilon(1)}{\Lambda_\varepsilon(1+|y_\varepsilon|)^\tau}\Big(\int_{B_\frac{2\delta}{\varepsilon}}|G_\varepsilon(y_\varepsilon,\eta)-G_\varepsilon(0,\eta)|^{p'}d\eta\Big)^\frac{1}{p'}\nonumber\\
&\leq \frac{o_\varepsilon(1)}{\Lambda_\varepsilon(1+|y_\varepsilon|)^\tau}\Big[\Big(\int_{\Sigma_1}|G_\varepsilon(y_\varepsilon,\eta)-G_\varepsilon(0,\eta)|^{p'}d\eta\Big)^\frac{1}{p'}+\Big(\int_{\Sigma_2}|G_\varepsilon(y_\varepsilon,\eta)-G_\varepsilon(0,\eta)|^{p'}d\eta\Big)^\frac{1}{p'}\nonumber\\
&\quad\quad\quad+\Big(\int_{\Sigma_3}|G_\varepsilon(y_\varepsilon,\eta)-G_\varepsilon(0,\eta)|^{p'}d\eta\Big)^\frac{1}{p'}\Big]=:B_1+B_2+B_3
\end{align}
Here, $\frac{1}{p}+\frac{1}{p'}=1$. Notice that
\begin{align}
B_1\leq \frac{o_\varepsilon(1)}{\Lambda_\varepsilon(1+|y_\varepsilon|)^\tau}\Big(\int_{\Sigma_1}(\ln|y_\varepsilon|)^{p'}+(\ln|\eta|)^{p'}\Big)^\frac{1}{p'}\leq\frac{C(\ln|y_\varepsilon|)\cdot|y_\varepsilon|^\frac{2+\theta}{p'}}{\Lambda_\varepsilon(1+|y_\varepsilon|)^\tau}.\nonumber
\end{align}
Here, $\theta$ is a small positive constant.
Recall that $\Lambda_\varepsilon\geq c_0>0$ for a uniform $c_0$ due to the assumption. Let $p$ sufficiently close to $1$, we get $p'$ sufficiently large. Under this assumption, we get
\begin{align}\label{ineq:B11}
B_1\leq \frac{C(\ln|y_\varepsilon|)\cdot|y_\varepsilon|^\frac{2+\theta}{p'}}{\Lambda_\varepsilon(1+|y_\varepsilon|)^\tau}=o_\varepsilon(1).
\end{align}
By a similar argument, we get
\begin{align}\label{ineq:B21}
B_2=o_\varepsilon(1).
\end{align}
For $B_3$, we get
\begin{align}\label{ineq:B31}
B_3\leq\frac{o_\varepsilon(1)}{\Lambda_\varepsilon(1+|y_\varepsilon|)^\tau}\Big(\frac{|y_\varepsilon|^{p'}d\eta}{|\eta|^{p'}}\Big)^\frac{1}{p'}\leq C\frac{|y_\varepsilon|^{\frac{2-p'}{p'}}}{(1+|y_\varepsilon|)^\tau}=o_\varepsilon(1).
\end{align}
Together they imply that
\begin{align}\label{Ineq:B}
B=o_\varepsilon(1).
\end{align}
(\ref{Ineq:A}) and (\ref{Ineq:B}) imply that
$A+B=o_\varepsilon(1)$.
We get a contradiction.

~

\noindent{\bf Case 3. $|y_\varepsilon|$ is unbounded and $|y_\varepsilon|\simeq\varepsilon^{-1}$.}

~

We can draw the contradiction in a similar way. Since the large part of the argument is the same, we only point out the difference.
As in \cite[Proposition 3.1]{BartolucciYangZhang2024}, we get
\begin{align}
|G_\varepsilon(y_\varepsilon,\eta)-G_\varepsilon(0,\eta)|\leq C\Big(|\ln|\eta||+|\ln|y_\varepsilon -\eta||+\ln\frac{1}{\varepsilon}\Big).\nonumber
\end{align}
Notice that
\begin{align}
\int_{B_{\frac{2\delta}{\varepsilon}}(0)}& \Big|G_\varepsilon(y_\varepsilon,\eta) -G_\varepsilon(0,\eta)\Big|\times\Big\{\Big|\sum_{j=1}^n a_{ij}e^{v_j(\eta)}\widetilde{w}^*_{t,i,\varepsilon}(\eta)\Big|+ \Big|\frac{\varepsilon_t^2\widetilde{g}^0_{i,\varepsilon}(\eta)}{\Lambda_\varepsilon\cdot \varepsilon^{-\alpha}O_{m^*}(1+|y_\varepsilon|)^\tau}\nonumber\\
&\quad+\sum_{t=1}^N \frac{r_{3,t}\widetilde{Z}_{3,t,i,\varepsilon}^*(\eta)}{\Lambda_\varepsilon\cdot \varepsilon^{-\alpha}O_{m^*}(1+|y_\varepsilon|)^\tau}\cdot \sum_{s=1}^N\widetilde{\chi}_s (\eta) e^{v_i(\eta+\frac{p_{t,\varepsilon}-p_{l,\varepsilon}}{\varepsilon_t})+ 2\ln\frac{\varepsilon_t}{\varepsilon_l}}\nonumber\\
&\quad +s_i\sum_{s=1}^N\frac{\widetilde{\chi}_s(\eta) e^{v_i(\eta+\frac{p_{t,\varepsilon}-p_{l,\varepsilon}}{\varepsilon_t})+ 2\ln\frac{\varepsilon_t}{\varepsilon_l}}}{\Lambda_\varepsilon\cdot\varepsilon^{-\alpha}O_{m^*} (1+|y_\varepsilon|)^\tau}\Big|\Bigg\}d\eta\nonumber\\
&\leq \Bigg\{ C\int_{B_{\frac{2\delta}{\varepsilon}(0)}}\Big(|\ln|\eta||+|\ln|y_\varepsilon -\eta||+\ln\frac{1}{\varepsilon}\Big) \Big|\sum_{j=1}^n a_{ij}e^{v_j(\eta)}\widetilde{w}^*_{t,i,\varepsilon}(\eta)\Big|d\eta\Bigg\}\nonumber\\
&\quad+\Bigg\{ C\int_{B_{\frac{2\delta}{\varepsilon}(0)}}\Big(|\ln|\eta||+|\ln|y_\varepsilon -\eta||+\ln\frac{1}{\varepsilon}\Big) \times\nonumber\\
&\quad\times\Big|\frac{\varepsilon_t^2\widetilde{g}^0_{i,\varepsilon}(\eta)}{\Lambda_\varepsilon\cdot \varepsilon^{-\alpha}O_{m^*}(1+|y_\varepsilon|)^\tau}+\sum_{t=1}^N \frac{r_{3,t}\widetilde{Z}_{3,t,i,\varepsilon}^*(\eta)}{\Lambda_\varepsilon\cdot \varepsilon^{-\alpha}O_{m^*}(1+|y_\varepsilon|)^\tau}\cdot \sum_{s=1}^N\widetilde{\chi}_s (\eta) e^{v_i(\eta+\frac{p_{t,\varepsilon}-p_{l,\varepsilon}}{\varepsilon_t})+2\ln\frac{\varepsilon_t}{\varepsilon_l}}\nonumber\\
&\quad +s_i\sum_{s=1}^N\frac{\widetilde{\chi}_s(\eta) e^{v_i(\eta+\frac{p_{t,\varepsilon}-p_{l,\varepsilon}}{\varepsilon_t})+2\ln\frac{\varepsilon_t}{\varepsilon_l}}}{\Lambda_\varepsilon \cdot\varepsilon^{-\alpha}O_{m^*}(1+|y_\varepsilon|)^\tau}\Big|d\eta\Bigg\}\nonumber\\
&=:A'+B'.\nonumber
\end{align}
As in {\bf Step 2}, by (\ref{Ineq:gwSize}) and (\ref{Ineq:rsSize}), we get
\begin{align}
A' &\leq C\int_{B_\frac{2\delta}{\varepsilon}(0)}\frac{|\ln|\eta||d\eta}{(1+|y_\varepsilon|)^\tau (1+|\eta|)^{m^*-\tau}}+ C\int_{B_\frac{2\delta}{\varepsilon}(0)}\frac{|\ln|y_\varepsilon -\eta||d\eta}{(1+|y_\varepsilon|)^\tau (1+|\eta|)^{m^*-\tau}}\nonumber\\
&\quad\quad+C\ln\frac{1}{\varepsilon}\cdot\int_{B_\frac{2\delta}{\varepsilon}(0)}\frac{d\eta}{(1+|y_\varepsilon|)^\tau (1+|\eta|)^{m^*-\tau}}\nonumber
\end{align}
Here,
\begin{align}
\int_{B_\frac{2\delta}{\varepsilon}(0)}\frac{|\ln|\eta||d\eta}{(1+|y_\varepsilon|)^\tau (1+|\eta|)^{m^*-\tau}} \leq
\frac{C}{(1+|y_\varepsilon|)^\tau} \int_0^{\frac{2\delta}{\varepsilon}}\frac{r|\ln r|dr}{(1+r)^{m^*-\tau}}=o_\varepsilon(1),\nonumber
\end{align}
\begin{align}
\int_{B_\frac{2\delta}{\varepsilon}(0)}\frac{|\ln|y_\varepsilon -\eta||d\eta}{(1+|y_\varepsilon|)^\tau (1+|\eta|)^{m^*-\tau}}&=\int_{B_1(y_\varepsilon)}\frac{|\ln|y_\varepsilon -\eta||d\eta}{(1+|y_\varepsilon|)^\tau (1+|\eta|)^{m^*-\tau}}\nonumber\\
&\quad\quad+\int_{B_\frac{2\delta}{\varepsilon}(0)\backslash B_1(y_\varepsilon)}\frac{|\ln|y_\varepsilon -\eta||d\eta}{(1+|y_\varepsilon|)^\tau (1+|\eta|)^{m^*-\tau}}\nonumber\\
&\leq \frac{C}{|y_\varepsilon|^{m^*}}+\frac{\ln\frac{C}{\varepsilon}}{(1+|y_\varepsilon|)^\tau}\int_{\mathbb{R}^2}\frac{d\eta}{(1+|\eta|)^{m^*-\tau}}\nonumber\\
&=o_\varepsilon(1)\nonumber
\end{align}
and
\begin{align}
\ln\frac{1}{\varepsilon}\cdot\int_{B_\frac{2\delta}{\varepsilon}(0)}\frac{d\eta}{(1+|y_\varepsilon|)^\tau (1+|\eta|)^{m^*-\tau}}\leq\frac{\ln\frac{1}{\varepsilon}}{(1+|y_\varepsilon|)^\tau}\int_{\mathbb{R}^2}\frac{d\eta}{(1+|\eta|)^{m^* -\tau}}.\nonumber
\end{align}
Together they imply
\begin{align}
A'=o_\varepsilon(1).\nonumber
\end{align}

On the other hand, we get
\begin{align}
B'\leq C\Big(\int_{B_\frac{2\delta}{\varepsilon}(0)}|\ln|\eta||^{p'}+|\ln|y_\varepsilon-\eta||^{p'}+|\ln\frac{1}{\varepsilon}|^{p'}\Big)^\frac{1}{p'}\cdot\frac{o_\varepsilon(1)}{(1+|y_\varepsilon|)^\tau}\leq\frac{C}{\varepsilon^\frac{2+\theta}{p'}}\frac{1}{(1+|y_\varepsilon|)^\tau}=o_\varepsilon(1)\nonumber
\end{align}
if we let $\theta$ small enough and $p$ sufficiently closed to $1$.

Therefore, for any case, we get a contradiction. Then, $\Lambda_\varepsilon\to\infty$ and Proposition
\ref{prop:w0Pointwise} is proved.
$\Box$

In such a way, we obtain a function defined in $\cup_{t=1}^N B_{3\delta_t}(p_{t,\varepsilon_t})$. Then, we get a  $C^2$ 0-extension $w_{i,\varepsilon}^0$ defined on $\Omega$ with
\begin{itemize}
    \item [$(1)$] $\|w^0_{i,\varepsilon}\|_{L^\infty(\Omega)}=O_{m^*}\varepsilon^{-\alpha-\tau}$;
    \item [$(2)$] $w^0_{i,\varepsilon}|_{B_{2\delta_t}(p_{t,\varepsilon})}=R^0_{t,i,\varepsilon}(|x-p_{t,\varepsilon}|)$ for any $t=1,\cdots,N$;
    \item [$(3)$] $w^0_{i,\varepsilon}|_{\Omega\backslash\cup_{t=1}^N B_{3\delta_t}(p_{t,\varepsilon})}=0$.
\end{itemize}

\subsection{Refined Estimates on the Error Term II: Parts of $w_\varepsilon$ with non-zero frequency}\label{Subsection:RefindeReduction2}

\subsubsection{Parts of $w_{\varepsilon}$ with Frequency no less than 2}

Now we consider the parts of $w_\varepsilon$ with frequency no less than 2. By the Fourier expansion, the $k$-frequency part of $w_{i,\varepsilon}$ is defined by
\begin{align}
w_{t,i,\varepsilon}^k:=\sin(k\theta)\cdot R_{1,t,i,\varepsilon}^k(r)+ \cos(k\theta)\cdot R_{2,t,i,\varepsilon}^k(r)\nonumber
\end{align}
with
\begin{align}\label{e:EquationR1k}
-\Big(\frac{d^2}{dr^2}+\frac{1}{r}\frac{d}{dr}-\frac{k^2}{r^2}\Big)\Big(\sum_{l=1}^n a^{li} R^k_{1,t,l,\varepsilon}\Big)-e^{v_i(r)}R_{1,t,i,\varepsilon}^k=O(\varepsilon^{-\alpha}O_{m^*}(1+r)^{\tau-2})
\end{align}
and
\begin{align}\label{e:EquationR2k}
-\Big(\frac{d^2}{dr^2}+\frac{1}{r}\frac{d}{dr}-\frac{k^2}{r^2}\Big)\Big(\sum_{l=1}^n a^{li} R^k_{2,t,l,\varepsilon}\Big)-e^{v_i(r)}R_{2,t,i,\varepsilon}^k=O(\varepsilon^{-\alpha}O_{m^*}(1+r)^{\tau-2}).
\end{align}
Notice that the right hand sides of (\ref{e:EquationR1k}) and of (\ref{e:EquationR2k}) are projections of  $g_{i,\varepsilon}$ on nodes higher than $1$.

Then, the function
\begin{align}\label{e:SummationW2}
\widetilde{w}^{II}_{t,i,\varepsilon}(x):=\sum_{k\geq 2}\Big(\sin(k\theta)\cdot R_{1,t,i,\varepsilon}^k(r)+ \cos(k\theta)\cdot R_{2,t,i,\varepsilon}^k(r)\Big)\mbox{ in }B_{\frac{4\delta_t}{\varepsilon_t}}(p_{t,\varepsilon})
\end{align}
satisfies that
\begin{align}
\widetilde{L}_{t,i,\varepsilon}\widetilde{w}^{II}_{t,\varepsilon}=
\left\{
\begin{aligned}
& \varepsilon_t^2\widetilde{g}_{t,i,\varepsilon}^{II} &\mbox{ for }y\in \cup_{s=1}^N B_{\frac{5\delta_s}{\varepsilon_s}}(p_{s,\varepsilon}) ;\\
& 0&\mbox{ for }\Omega\backslash \cup_{s=1}^N B_{\frac{5\delta_s}{\varepsilon_s}}(p_{s,\varepsilon}).
\end{aligned}
\right.
\end{align}
Here, $\widetilde{g}_{t,i,\varepsilon}^{II}$ is the projection of $\widetilde{g}_{t,i,\varepsilon}$ over nodes higher than $1$ in $\cup_{s=1}^N B_{\frac{5\delta_s}{\varepsilon_s}}(p_{s,\varepsilon})$.  By a similar computation as in Proposition \ref{prop:Recution} and in Subsection \ref{subsection:RefindeReduction1}, we get $\|\widetilde{g}_{t,i,\varepsilon}^{II}\|_{\mathbb{Y}_\varepsilon}=\varepsilon^{-\alpha}O_{m^*}$. Denoting $w^{II}_{i,\varepsilon}(x)=\widetilde{w}^{II}_{t,\varepsilon}(\frac{x-p_{t,\varepsilon}}{\varepsilon_t})$. Notice that $(w^{II}_{1,\varepsilon},\cdots,w^{II}_{n,\varepsilon})\in E_\varepsilon$ due to its frequency. By Theorem
\ref{t:invertibility}, we get
\begin{align}\label{Ineq:WIIPointwise}
\|w^{II}_{t,\varepsilon}\|_{L^\infty(\Omega)}=\varepsilon^{-\alpha-\tau}O_{m^*}.
\end{align}
Here, $\alpha$ and $\tau$ are two small positive constants. Then, we get
\begin{align}
-\Big(\frac{d^2}{dr^2}+\frac{1}{r}\frac{d}{dr}-\frac{k^2}{r^2}\Big) R^k_{1,t,i,\varepsilon}=O(\varepsilon^{-\alpha-\tau}O_{m^*}(1+r)^{\tau-2}).\nonumber
\end{align}
An elementary ODE method (see \cite[pp.2610-2616]{LinZhang2013}) implies that
$$|R_{1,i,\varepsilon}^k(r)|\leq C(\varepsilon r)^k+\frac{C\varepsilon^{m^*-2-\alpha-\tau}}{k^2}r^2(1+r)^{\tau-2}.$$ A similar computation can be given for $R^k_{2,t,i,\varepsilon}$. Then, we get for any $i=1,\cdots,n$ and any $t=1,\cdots,N$,
$|w^{II}_{i,\varepsilon}(x-p_{t,\varepsilon})|\leq C\varepsilon^2|x-p_{t,\varepsilon}|^2$ for
$x\in B_{\delta_t}(p_{t,\varepsilon})$.

\subsubsection{The 1-Frequency Part of $w_{\varepsilon}$}\label{susbection:1Frequence}

Now we study the 1-frequency part of $w_\varepsilon$.
To begin with, we define the function
\begin{align}
w_{i,\varepsilon}^I:=w_{i,\varepsilon}-w^0_{i,\varepsilon}-w_{i,\varepsilon}^{II}.\nonumber
\end{align}

By the above deduction,
\begin{itemize}
    \item [$(1)$] $w^I_\varepsilon:=(w^I_{1,\varepsilon},\cdots,w^I_{n,\varepsilon})\in E_\varepsilon$;
\item [$(2)$] It holds that
\begin{align}
\|w^I_{i,\varepsilon}\|_{L^\infty(\Omega)}\leq
\left\{
\begin{aligned}
& C\varepsilon^{m^*-2-\frac{\alpha}{2}}(\ln\frac{1}{\varepsilon})^2&\mbox{ if } m^*\leq 3;\nonumber\\
& C\varepsilon\ln\frac{1}{\varepsilon}&\mbox{ if }m^*_i>3.\nonumber
\end{aligned}
\right.
\end{align}
\end{itemize}

Assertion $(1)$ holds because that $w_\varepsilon^0,w_\varepsilon^{II}, w_\varepsilon\in E_\varepsilon$. Notice that $w_\varepsilon\in E_\varepsilon$ due to Proposition \ref{prop:Recution} and $w_\varepsilon^{II}\in E_\varepsilon$ due to its frequency. For $w_\varepsilon^0$, we get
\begin{align}
\sum_{i=1}^n\int_\Omega\chi_t\partial_{x_1}v_i(\frac{x-p_{t,\varepsilon}}{\varepsilon_t})w^0_{i,\varepsilon} dx=\sum_{i=1}^n\int_\Omega\chi_t\partial_{x_2}v_i(\frac{x-p_{t,\varepsilon}}{\varepsilon_t})w^0_{i,\varepsilon} dx=0\nonumber
\end{align}
due to its symmetry and
\begin{align}
0&=\sum_{t=1}^N\int_\Omega \chi_t e^{v_i(\frac{x-p_{t,\varepsilon}}{\varepsilon_t})+2\ln\frac{1}{\varepsilon_t}}w_i dx=\sum_{t=1}^N\int_{B_{5\delta_t}(p_{t,\varepsilon})} e^{v_i(\frac{x-p_{t,\varepsilon}}{\varepsilon_t})+2\ln\frac{1}{\varepsilon_t}}w_i dx\nonumber\\
&=2\pi\sum_{t=1}^N\int_{0}^{5\delta_t}dr \int_{\partial B_r(p_{t,\varepsilon})}e^{v_i(\frac{x-p_{t,\varepsilon}}{\varepsilon_t})+2\ln\frac{1}{\varepsilon_t}}w_i dx=2\pi\sum_{t=1}^N\int_{0}^{5\delta_t}dr \int_{\partial B_r(p_{t,\varepsilon})}e^{v_i(\frac{x-p_{t,\varepsilon}}{\varepsilon_t})+2\ln\frac{1}{\varepsilon_t}}w_{i,\varepsilon}^0 dx.\nonumber
\end{align}
And $\sum_{i=1}^n\int_\Omega\sum_{s=1}^N\chi_s e^{v_i(\frac{x-p_{t,\varepsilon}}{\varepsilon_t})+2\ln\frac{1}{\varepsilon_t}} Z_{3,t,i,\varepsilon}^* w_{\varepsilon,i}^0dx=0$ for the same reason.
Assertion $(2)$ holds due to Proposition \ref{prop:w0Pointwise}, (\ref{Ineq:R0Pointwise}) and (\ref{Ineq:WIIPointwise}).

Moreover, due to the local expansion of $w^0_\varepsilon$ and of $w_\varepsilon^{II}$, we know that for any $t=1,\cdots,N$
\begin{align}\label{e:W1Local}
w_{i,\varepsilon}^I\big(\frac{x-p_{t,\varepsilon}}{\varepsilon_t}\big)=\sin\theta\cdot R_{1,t,i,\varepsilon}^1(r)+ \cos\theta\cdot  {R_{2,t,i,\varepsilon}^1(r)}
\end{align}
with
\begin{align}
&-\Big(\frac{d^2}{dr^2}+\frac{1}{r}\frac{d}{dr}-\frac{1}{r^2}\Big)\Big(\sum_{j=1}^n a^{ji}  {R^1_{l,t,j,\varepsilon}}\Big)-e^{v_i(r)}R_{l,t,i,\varepsilon}^1\nonumber\\
&=\left\{
\begin{aligned}
& \frac{\varepsilon_t^2}{2\pi}\int_0^{2\pi}\widetilde{g}_{t,i,\varepsilon}(r\cos\theta,r\sin\theta)\cos\theta d\theta+\sum_{j=1}^2\sum_{t=1}^N \frac{r_{j,t}}{2\pi}\int_0^{2\pi}\widetilde{Z}_{j,t,i,\varepsilon}^*(r\cos\theta,r\sin\theta)\widetilde{\chi}_t e^{v_i(r\cos\theta,r\sin\theta)} \cos\theta d\theta\nonumber\\
&\quad\quad\quad\quad\quad\quad\quad\mbox{ if } l=1;\nonumber\\
& \frac{\varepsilon_t^2}{2\pi}\int_0^{2\pi}\widetilde{g}_{t,i,\varepsilon}(r\cos\theta,r\sin\theta)\sin\theta d\theta+\sum_{j=1}^2\sum_{t=1}^N \frac{r_{j,t}}{2\pi}\int_0^{2\pi}\widetilde{Z}_{j,t,i,\varepsilon}^*(r\cos\theta,r\sin\theta)\widetilde{\chi}_t e^{v_i(r\cos\theta,r\sin\theta)} \sin\theta d\theta\nonumber\\
&\quad\quad\quad\quad\quad\quad\quad\mbox{ if } l=2.\nonumber
\end{aligned}
\right.
\end{align}

By a similar method as in \cite[pp.2610]{LinZhang2013}, we find a sequence $\vec{C}_\varepsilon\in\mathbb{R}^2$ bounded in $\varepsilon$ and
$w_{i,\varepsilon}^I(x)\sim\varepsilon_t \vec{C}_\varepsilon\cdot(x-p_{t,\varepsilon})$ as $x\to p_{t,\varepsilon}$.

\subsubsection{Conclusion: A Refined Lyapunov-Schmidt Reduction}
In this subsection, we provide a refined version of Proposition \ref{prop:Recution} based on Fourier analysis.

\begin{proposition}\label{prop:RecutionFiner}
Under Hypotheses (\ref{ASSUMPTION0}), (\ref{ASSUMPTION}) and (\ref{ASSUMPTION1}), there exists a constant $\varepsilon_0>0$ and for any $\varepsilon\in(0,\varepsilon_0)$, there exists a $w_\varepsilon=(w_\varepsilon^1,\cdots,w_\varepsilon^n)\in E_\varepsilon$ satisfies
\begin{align}
Q_\varepsilon(L_{1,\varepsilon}(w_\varepsilon)-g_{1,\varepsilon}(x,w_\varepsilon),\cdots,L_{n,\varepsilon}(w_\varepsilon)-g_{n,\varepsilon}(x,w_\varepsilon))=0.\nonumber
\end{align}
Moreover, there exists three $C^1$ maps $w_\varepsilon^0$, $w_\varepsilon^I$, $w_\varepsilon^{II}$ of $\varepsilon$ in $E_\varepsilon$ with $w_\varepsilon=w_\varepsilon^0+w_\varepsilon^I+w_\varepsilon^{II}$. Here, $w_\varepsilon^0$ satisfies
\begin{itemize}
    \item [$(W^0_{1})$] $\|w^0_{i,\varepsilon}\|_{L^\infty(\Omega)}=O_{m^*}\varepsilon^{-\alpha-\tau}$;
    \item [$(W^0_{2})$] $w^0_{i,\varepsilon}|_{B_{2\delta_t}(p_{t,\varepsilon})}=R^0_{t,i,\varepsilon}(|x-p_{t,\varepsilon}|)$ for any $t=1,\cdots,N$. Here, $R_{t,i,\varepsilon}^0$ are radial functions satisfying (\ref{e:0Frequency});
    \item [$(W^0_{3})$] $w^0_{i,\varepsilon}|_{\Omega\backslash\cup_{t=1}^N B_{3\delta_t}(p_{t,\varepsilon})}=0$;
    \item [$(W^0_{4})$] $\Big|w^0_{i,\varepsilon}\big(\frac{x-p_{t,\varepsilon}}{\varepsilon_t}\big)\Big|\leq \varepsilon^{-\alpha} O_{m^*}\big(1+\Big|\frac{y-p_{t,\varepsilon}}{\varepsilon_t}\Big|\big)^\tau$ for $x\in B_{3\delta_t}(p_{t,\varepsilon})$ and for any $t=1,\cdots,N$.
\end{itemize}
$w_\varepsilon^{I}$ satisfies
\begin{itemize}
\item [$(W^I_{1})$] It holds that
\begin{align}
\|w^I_{i,\varepsilon}\|_{L^\infty(\Omega)}\leq
\left\{
\begin{aligned}
& C\varepsilon^{m^*-2-\frac{\alpha}{2}}(\ln\frac{1}{\varepsilon})^2&\mbox{ if } m^*\leq 3;\nonumber\\
& C\varepsilon\ln\frac{1}{\varepsilon}&\mbox{ if }m^*_i>3.\nonumber
\end{aligned}
\right.
\end{align}
   \item [$(W^I_{2})$] for any $t=1,\cdots,N$, $w^I_{i,\varepsilon}$ has an expansion in the form of (\ref{e:W1Local}) in $B_{2\delta_t}(p_{t,\varepsilon})$;
   \item [$(W^I_3)$] there exists a sequence $\vec{C}_\varepsilon\in\mathbb{R}^2$ bounded in $\varepsilon$ and
$w_{i,\varepsilon}^I(x)\sim\varepsilon_t \vec{C}_\varepsilon\cdot(x-p_{t,\varepsilon})$ as $x\to p_{t,\varepsilon}$.
\end{itemize}
And $w_\varepsilon^{II}$ satisfies
\begin{itemize}
    \item [$(W^{II}_{1})$] $\|w^{II}_{i,\varepsilon}\|_{L^\infty(\Omega)}=O_{m^*}\varepsilon^{-\alpha-\tau}$;
    \item [$(W^{II}_{2})$] $w^{II}_{i,\varepsilon}|_{B_{2\delta}(p_{t,\varepsilon})}=\widetilde{w}^{II}_{i,\varepsilon}(x-p_{t,\varepsilon})$ for any $t=1,\cdots,N$ and $\widetilde{w}^{II}_{i,\varepsilon}$ satisfies (\ref{e:SummationW2});
    \item [$(W^{II}_{3})$] for any $i=1,\cdots,n$ and any $t=1,\cdots,N$,
$|w^{II}_{i,\varepsilon}(x-p_{t,\varepsilon})|\leq C\varepsilon^2|x-p_{t,\varepsilon}|^2$ for
$x\in B_{\delta_t}(p_{t,\varepsilon})$.
\end{itemize}
Here, $\alpha$ and $\tau$ are generically small positive numbers.
\end{proposition}

\section{Finite-dimensional Problems and Proof of Theorem \ref{t:MAIN}}
In this section, we complete the proof by analyzing the reduced finite-dimensional problem with a standard approach used in  \cite{delPinoWei,Huang2019}. The central issue arises when the projection of $S(U_{1,\varepsilon}+w_{1,\varepsilon},\cdots,U_{n,\varepsilon}+w_{n,\varepsilon})$ onto $F_\varepsilon$ (see \eqref{def:F}) vanishes. Consequently, the primary focus is on the terms involving the projections of $S(U_{1,\varepsilon}+w_{1,\varepsilon}^I,\cdots,U_{n,\varepsilon}+w_{n,\varepsilon}^I)$.

In Subsections \ref{Subsection:Projection1}, \ref{Subsection:Projection2}, and \ref{Subsection:Projection3}, we explore the projections onto $Z^*_{j,t,\varepsilon}$ for $j=1,2,3$ and $t=1,\cdots,N$. Importantly, the projections onto $e_i=(0,\cdots,1,\cdots,0)$ for $i=1,\cdots,n$ clearly vanish due to the double periodicity of the arguments. Based on the computations presented above, we complete the proof of Theorem \ref{t:MAIN} in Subsection \ref{Subsection:proof}.

\subsection{A Finite-dimensional Problem}
Denoting
\begin{align}
L_\varepsilon w_\varepsilon =(L_{1,\varepsilon}w_\varepsilon,\cdots, L_{n,\varepsilon}w_\varepsilon)\nonumber
\end{align}
and
\begin{align}
g_\varepsilon(x,w_\varepsilon)=(g_{1,\varepsilon}(x,w_\varepsilon),\cdots, g_{n,\varepsilon}(x,w_\varepsilon)),\nonumber
\end{align}
the reduction scheme ensures the existence of $w_\varepsilon=(w_{1,\varepsilon},\cdots,w_{n,\varepsilon})$ such that
\begin{align}
L_\varepsilon w_\varepsilon -g_\varepsilon(x,w_\varepsilon)=&\sum_{j=1}^2\sum_{t=1}^N r_{j,t}\Big(\chi_t\partial_{x_j}v_1(\frac{x-p_{t,\varepsilon}}{\varepsilon_t})\sum_{s=1}^Ne^{v_1(\frac{x-p_{s,\varepsilon}}{\varepsilon_s})+2\ln\frac{1}{\varepsilon_s}},\cdots,\nonumber\\
&\quad\quad\chi_t\partial_{x_j}v_n(\frac{x-p_{t,\varepsilon_t}}{\varepsilon_t})\sum_{s=1}^Ne^{v_n(\frac{x-p_{s,\varepsilon}}{\varepsilon_s})+2\ln\frac{1}{\varepsilon_s}}\Big)\nonumber\\
&+\sum_{t=1}^N r_{3,t}\Big(Z_{3,t,1,\varepsilon }^*\sum_{s=1}^N\chi_s e^{v_1(\frac{x-p_{s,\varepsilon}}{\varepsilon_s})+2\ln\frac{1}{\varepsilon_s}},\cdots,Z_{3,t,n,\varepsilon }^*\sum_{s=1}^N\chi_s e^{v_n(\frac{x-p_{s,\varepsilon}}{\varepsilon_s})+2\ln\frac{1}{\varepsilon_s}}\Big)\nonumber\\
&+\Big(s_1\sum_{t=1}^N\chi_s e^{v_1(\frac{x-p_{s,\varepsilon}}{\varepsilon_s})+2\ln\frac{1}{\varepsilon_s}},\cdots,s_n\sum_{s=1}^N\chi_t e^{v_n(\frac{x-p_{s,\varepsilon}}{\varepsilon_s})+2\ln\frac{1}{\varepsilon_s}}\Big)\nonumber
\end{align}
with
\begin{align}\label{e:FiniteDimension}
L_{i,\varepsilon}w_\varepsilon - g_{i,\varepsilon}(x,w_\varepsilon)=-\Delta (U_{\varepsilon}^i +w_\varepsilon^i)-\rho_{i,\varepsilon}\Bigg(\frac{h_i(x)e^{U_{i,\varepsilon}+w_{i,\varepsilon}}}{\int_\Omega h_i(x)e^{U_{i,\varepsilon}+w_{i,\varepsilon}}}-1\Bigg)
\end{align}
for $i=1,\cdots,n$ and for some constant $r_{1,1},\cdots,r_{1,N},r_{2,1},\cdots,r_{2,N},r_{3,1},\cdots,r_{3,N},s_1,\cdots,s_n\in\mathbb{R}$. The following lemma provides a condition to find a solution to Problem (\ref{e:001}).

\begin{lemma}\label{l:Coefficients}
If
\begin{align}
\sum_{i=1}^n\int_\Omega\Big[L_{i,\varepsilon}w_\varepsilon -g_{i,\varepsilon}(x,w_\varepsilon)\Big]Z_{j,t,i,\varepsilon}^*=0
\end{align}
for $t=1,\cdots,N$ and  $j=1,2,3$, then $r_{1,1}=\cdots=r_{1,N}=r_{2,1}=\cdots=r_{2,N}=r_{3,1}=\cdots=r_{3,N}=s_1=\cdots=s_n=0$.
\end{lemma}
\noindent{\bf Proof.}
For $j'=1,2$ and $t'=1,\cdots,N$, multiplying (\ref{e:FiniteDimension}) by $Z_{j',t',i,\varepsilon}^*$, integrating over $\Omega$ and summing with respect to $i$, we get
\begin{align}
0=&\sum_{i=1}^n\int_\Omega \Big[L_{i,\varepsilon}w_\varepsilon -g_{i,\varepsilon}(x,w_\varepsilon)\Big]
Z_{j,i,\varepsilon}^*\nonumber\\
=&\sum_{i=1}^n \Bigg(\sum_{j=1}^3 r_{j,t} \int_\Omega \sum_{s=1}^N\chi_s e^{v_i(\frac{x-p_{s,\varepsilon}}{\varepsilon_s})+2\ln\frac{1}{\varepsilon_s}} Z_{j,t,i,\varepsilon}^* Z_{j',t',i,\varepsilon}^* +s_i \int_\Omega \chi_t e^{v_i(\frac{x-p_{t,\varepsilon}}{\varepsilon_t})+2\ln\frac{1}{\varepsilon_s}} Z_{j',t',i,\varepsilon}\Bigg).\nonumber
\end{align}
By the symmetry, we get $r_{1,1}=\cdots=r_{1,N}=r_{1,i}=r_{2,N}=0$.

Now we study $r_{3,1},\cdots,r_{3,N},s_1,\cdots,s_n$. To do this, notice that
\begin{align}\label{e:LINEAR01}
\int_\Omega\Big(\sum_{t=1}^N  r_{3,t}\sum_{s=1}^N\chi_s e^{v_i(\frac{x-p_{s,\varepsilon}}{\varepsilon_s})+2\ln\frac{1}{\varepsilon_s}}Z_{3,t,i,\varepsilon}^* +s_i \sum_{s=1}^N\chi_s e^{v_i(\frac{x-p_{s,\varepsilon}}{\varepsilon_s})+2\ln\frac{1}{\varepsilon_s}}\Big)=0.
\end{align}
This is due to (\ref{e:FiniteDimension}) and the fact that $U_{\varepsilon}^i+w_\varepsilon^i$ are double periodic.
On the other hand, multiplying $Z_{3,i,\varepsilon}^*$ on both sides of (\ref{e:FiniteDimension}), integrating over $\Omega$ and summing with respect to $i$, we get
\begin{align}\label{e:LINEAR02}
\sum_{i=1}^n \Bigg(\sum_{t=1}^N r_{3,t} \int_\Omega \sum_{s=1}^N\chi_s e^{v_i(\frac{x-p_{t,\varepsilon}}{\varepsilon_s})+2\ln\frac{1}{\varepsilon_s}} Z_{3,t,i,\varepsilon}^* Z_{3,t',i,\varepsilon}^* +s_i \int_\Omega \sum_{s=1}^N\chi_s e^{v_i(\frac{x-p_{s,\varepsilon}}{\varepsilon_s})+2\ln\frac{1}{\varepsilon_s}} Z_{3,t',i,\varepsilon}^*\Bigg)=0.
\end{align}

By a direct computation, we get for $t=t'$
\begin{align}
\sum_{i=1}^n\int_\Omega\sum_{s=1}^N\chi_s e^{v_i(\frac{x-p_{s,\varepsilon}}{\varepsilon_s})+2\ln\frac{1}{\varepsilon_s}}Z_{3,t,i,\varepsilon}^* Z_{3,t',i,\varepsilon}^* =\sum_{i=1}^n&\Big[(N-1)(m_i^*-2)^2\int_{\mathbb{R}^2}e^{v_i}\nonumber\\
&+\int_{\mathbb{R}^2}e^{v_i}|x\cdot\nabla v_i+2|^2+o_\varepsilon(1)\Big]\varepsilon_t^2,\nonumber
\end{align}
for $t\neq t'$
\begin{align}
\sum_{i=1}^n\int_\Omega\sum_{s=1}^N\chi_s e^{v_i(\frac{x-p_{s,\varepsilon}}{\varepsilon_s})+2\ln\frac{1}{\varepsilon_s}}Z_{3,t,i,\varepsilon}^* Z_{3,t',i,\varepsilon}^* =\sum_{i=1}^n&\Big[(N-2)(m_i^*-2)^2\int_{\mathbb{R}^2}e^{v_i}+o_\varepsilon(1)\Big]\varepsilon_t\varepsilon_{t'}\nonumber
\end{align}
\begin{align}
\int_\Omega\sum_{s=1}^N\chi_s e^{v_i(\frac{x-p_{s,\varepsilon}}{\varepsilon_s})+2\ln\frac{1}{\varepsilon_s}}Z_{3,t',i,\varepsilon}^*=\Big[(N-1)(2-m_i^*)\int_{\mathbb{R}^2}e^{v_i}+o_\varepsilon(1)\Big]\nonumber
\end{align}
and
\begin{align}
\int_\Omega\sum_{s=1}^N\chi_s e^{v_i(\frac{x-p_{s,\varepsilon_s}}{\varepsilon_s})+2\ln\frac{1}{\varepsilon_s}}=N\int_{\mathbb{R}^2}e^{v_i}.\nonumber
\end{align}
Now we rewrite (\ref{e:LINEAR01}) and (\ref{e:LINEAR02}) into a linear systems. To do this, we need to introduce several vectors and matrices.
We denote
\begin{align}
X=(\varepsilon_1 r_{3,1},\cdots,\varepsilon_N r_{3N},s_1,\cdots,s_n)^T.\nonumber
\end{align}
Define
$G_1$ to be the $N\times N$ martix such that
\begin{align}
(i,i)\mbox{-element of }G_1=(N-1)\sum_{i=1}^n(m_i^*-2)^2\int_{\mathbb{R}^2}e^{v_i}dx+\sum_{i=1}^n\int_{\mathbb{R}^2}e^{v_i}(\nabla v_i\cdot x+2)^2 dx\nonumber
\end{align}
and
\begin{align}
(i,j)\mbox{-element of }G_1=(N-2)\sum_{i=1}^n(m_i^*-2)^2\int_{\mathbb{R}^2}e^{v_i}dx.\nonumber
\end{align}
for $i\neq j$.
Define the $n\times N$ matrix
\[G_3=
\begin{pmatrix}
  (N-1)(2-m_1^*)\int_{\mathbb{R}^2}e^{v_1} & \cdots &  (N-1)(2-m_1^*) \int_{\mathbb{R}^2}e^{v_1}\\
  \vdots & \cdots & \vdots \\
  (N-1)(2-m_n^*) \int_{\mathbb{R}^2}e^{v_n} & \cdots & (N-1)(2-m_n^*) \int_{\mathbb{R}^2}e^{v_n} \\
\end{pmatrix}.
\]
Define the $n\times n$ diagonal matrix
\begin{align}
G_2=\mbox{diag}\Big(N\int_{\mathbb{R}^2}e^{v_1},\cdots, N\int_{\mathbb{R}^2}e^{v_n}\Big),\nonumber
\end{align}
which is invertible.
Then, the above gives us
\[
\begin{pmatrix}
G_1+o_\varepsilon(1) & G_3^T+o_\varepsilon(1) \\
G_3+o_\varepsilon(1) & G_2+o_\varepsilon(1)
\end{pmatrix}
X=0.
\]

Now we are in the position to prove the invertibility of the matrix
\[
\begin{pmatrix}
G_1 & G_3^T \\
G_3 & G_2
\end{pmatrix}.
\]
To do this, we only need to check the invertibility of $G_1 -G_3^T G_2^{-1} G_3$. By a direct computation, all of the elements of $G_3^T G_2^{-1} G_3$ are equal to
\begin{align}
\frac{(N-1)^2}{N}\sum_{i=1}^n(m_i^*-2)^2\int_{\mathbb{R}^2}e^{v_i}.\nonumber
\end{align}
Then, we get
\begin{align}
(i,i)\mbox{-element of }G_1 -G_3^T G_2^{-1} G_3&=\Big(N-1-\frac{(N-1)^2}{N}\Big)\sum_{i=1}^n(m_i^*-2)^2\int_{\mathbb{R}^2}e^{v_i}dx\nonumber\\
&\quad\quad+\sum_{i=1}^n\int_{\mathbb{R}^2}e^{v_i}(\nabla v_i\cdot x+2)^2 dx\nonumber
\end{align}
and
\begin{align}
(i,j)\mbox{-element of }G_1 -G_3^T G_2^{-1} G_3=\Big(N-2-\frac{(N-1)^2}{N}\Big)\sum_{i=1}^n(m_i^*-2)^2\int_{\mathbb{R}^2}e^{v_i}dx.\nonumber
\end{align}
for $i\neq j$.
This is evidently invertible, which implies that
$r_{3,1}=\cdots=r_{3,N}=s_1=\cdots=s_n=0$. $\Box$

\begin{lemma}\label{l:ERROR}
It holds that
\begin{align}
\sum_{i=1}^n\int_\Omega\Big[L_{i,\varepsilon}w_\varepsilon -g_{i,\varepsilon}(x,w_\varepsilon)\Big]Z_{j,t,i,\varepsilon}^*&=\sum_{i=1}^n\int_\Omega\Bigg[-\Delta \Big(U_\varepsilon^i+\sum_{l=1}^n a^{il}w^I_{i,\varepsilon}\Big) -\rho_{i,\varepsilon}\Bigg(\frac{h_i e^{U_{i,\varepsilon}+w^I_{i,\varepsilon}}}{\int_\Omega h_i e^{U_i,\varepsilon+w^I_{i,\varepsilon}}}-1\Bigg)\Bigg]Z_{i,j,\varepsilon}^*\nonumber\\
&\quad+R_{j,t,\varepsilon}\nonumber
\end{align}
with
$R_{j,t,\varepsilon}=o_\varepsilon(1)\varepsilon O_{m^*}$
for any $t=1,\cdots,N$ and any $j=1,2,3$. Here, $w^I_{i,\varepsilon}$ is the $i$-th component of the 1-frequency part of $w_\varepsilon$ defined in Subsection \ref{susbection:1Frequence}.
\end{lemma}
\noindent{\bf Proof.}
We only consider the cases for $j=1,2$ and the case $j=3$ is similar. Since $w_{i,\varepsilon}=w_{i,\varepsilon}^0+w_{i,\varepsilon}^I+w_{i,\varepsilon}^{II}$ (see Proposition \ref{prop:RecutionFiner}) for any $i=1,\cdots,n$, we get
\begin{align}\label{e:NewExpansionTOTAL}
&L_{i,\varepsilon}w_\varepsilon-g_{i,\varepsilon}(x,w_\varepsilon)\nonumber\\
=&L_{i,\varepsilon}(w_\varepsilon^0+w_\varepsilon^{II})-\Delta(U^i_\varepsilon +\sum_{l=1}^n a^{il}w_{l,\varepsilon}^I)+\rho_{i,\varepsilon}-\rho_{i,\varepsilon}\frac{h_i(x)e^{U_{i,\varepsilon}+w_{i,\varepsilon}^I}}{\int_\Omega h_i(x)e^{U_{i,\varepsilon}+w_{i,\varepsilon}^I}}\nonumber\\
&\quad+\Bigg\{\Big(K_{i,\varepsilon}(x)-\rho_{i,\varepsilon}\frac{h_i(x)e^{U_{i,\varepsilon}+w_{i,\varepsilon}^I}}{\int_\Omega h_i(x)e^{U_{i,\varepsilon}+w_{i,\varepsilon}^I}}\Big)(w_{i,\varepsilon}^0+w_{i,\varepsilon}^{II})\Bigg\}\nonumber\\
&\quad-\Bigg\{\frac{K_{i,\varepsilon}(x)}{\rho_i^*}\int_\Omega K_{i,\varepsilon}(x)(w_{i,\varepsilon}^0+w_{i,\varepsilon}^{II})-\rho_{i,\varepsilon}\frac{h_i(x)e^{U_{i,\varepsilon}+w_{i,\varepsilon}^I}}{(\int_\Omega h_i(x)e^{U_{i,\varepsilon}+w_{i,\varepsilon}^I})^2}\times\nonumber\\
&\quad\quad\quad\times\int_\Omega h_i(x)e^{U_{i,\varepsilon}+w_{i,\varepsilon}^{I}}(w_{i,\varepsilon}^0+w_{i,\varepsilon}^{II})\Bigg\}\nonumber\\
&\quad-\Bigg\{\rho_{i,\varepsilon}\frac{h_i(x)e^{\ U_{i,\varepsilon}+w_{i,\varepsilon}}}{\int_\Omega h_i(x)e^{ U_{i,\varepsilon}+w_{i,\varepsilon}}} - \rho_{i,\varepsilon}\frac{h_i(x)e^{ U_{i,\varepsilon}+w_{i,\varepsilon}^I}}{\int_\Omega h_i(x)e^{ U_{i,\varepsilon}+w_{i,\varepsilon}^I}} - \rho_{i,\varepsilon}\Bigg(\frac{h_i(x)e^{ U_{i,\varepsilon}+w_{i,\varepsilon}^I}}{\int_\Omega h_i(x)e^{ U_{i,\varepsilon}+w_{i,\varepsilon}^I}}\Bigg)(w^0_{i,\varepsilon}+w_{i,\varepsilon}^{II})\nonumber\\
&\quad\quad\quad+\frac{\rho_{i,\varepsilon}h_i(x)e^{ U_{i,\varepsilon}+w_{i,\varepsilon}^I}}{\Big(\int_\Omega h_i(x)e^{ U_{i,\varepsilon}+w_{i,\varepsilon}^I}\Big)^2}\int_\Omega\Big(h_i(x)e^{ U_{i,\varepsilon}+w_{i,\varepsilon}^I}(w^0_{i,\varepsilon}+w_{i,\varepsilon}^{II})\Big)\Bigg\}\nonumber\\
&=:L_{i,\varepsilon}(w_\varepsilon^0+w_\varepsilon^{II})-\Delta(U^i_\varepsilon +\sum_{l=1}^n a^{il}w_{l,\varepsilon}^I)+\rho_{i,\varepsilon}-\rho_{i,\varepsilon}\frac{h_i(x)e^{U_{i,\varepsilon}+w_{i,\varepsilon}^I}}{\int_\Omega h_i(x)e^{U_{i,\varepsilon}+w_{i,\varepsilon}^I}}\nonumber\\
&\quad+A'_{2,i}-A'_{3,i}-N'_{i,\varepsilon}.
\end{align}
Then, for any $t=1,\cdots,N$ and any $j=1,2,3$, we get
\begin{align}
R_{j,t,\varepsilon}=\sum_{i=1}^n\int_\Omega\Big(L_{i,\varepsilon}(w_\varepsilon^0+w_\varepsilon^{II})+A'_{2,i}-A'_{3,i}-N'_{i,\varepsilon}\Big)Z_{i,j,\varepsilon}^*.
\end{align}
Here,
\begin{align}
\sum_{i=1}^n L_{i,\varepsilon}(w_\varepsilon^0+w_\varepsilon^{II})Z_{i,j,\varepsilon}^*= \sum_{i=1}^n L_{i,\varepsilon}Z_{i,j,\varepsilon}^*(w_\varepsilon^0+w_\varepsilon^{II})=0\nonumber
\end{align}
and
\begin{align}
&\sum_{i=1}^n
 \int_\Omega A'_{2,i}Z_{i,j,\varepsilon}^*\nonumber\\
&=\sum_{i=1}^n \Bigg\{\Big(K_{i,\varepsilon}(x)-\rho_{i,\varepsilon}\frac{h_i(x)e^{U_{i,\varepsilon}+w_{i,\varepsilon}^I}}{\int_\Omega h_i(x)e^{U_{i,\varepsilon}+w_{i,\varepsilon}^I}}\Big)(w_{i,\varepsilon}^0+w_{i,\varepsilon}^{II})\Bigg\}Z_{i,j,\varepsilon}^*\nonumber\\
&=O(1)\sum_{i=1}^n\sum_{t=1}^N\int_{B_{5\delta}(p_{t,\varepsilon})} e^{v_i(\frac{x-p_{t,\varepsilon}}{\varepsilon_t})+2\ln\frac{1}{\varepsilon_t}}\Bigg\{1
-\frac{h_i(x)}{h_i(p_{t,\varepsilon})}e^{2\pi m_i^*( {G^*}(x;p_{t,\varepsilon})- {G^*}(p_{t,\varepsilon};p_{t,\varepsilon}))} \Bigg\}\times\nonumber\\
&\quad\times (w_{i,\varepsilon}^0+w^{II}_{i,\varepsilon})Z_{i,j,\varepsilon}^*dx+O\Big[\|w^I_\varepsilon\|_\infty(\|w_\varepsilon^0\|_\infty+\|w_\varepsilon^{II}\|_\infty)\Big]\cdot O(\varepsilon)\nonumber\\
&=\sum_{i=1}^n\sum_{t=1}^N \int_{B_{5\delta}(p_{t,\varepsilon})} e^{v_i(\frac{x-p_{t,\varepsilon}}{\varepsilon_t})+2\ln\frac{1}{\varepsilon_t}}O(|x|^2) (w_{i,\varepsilon}^0+w^{II}_{i,\varepsilon})Z_{i,j,\varepsilon}^*dx +O\Big[\|w^I_\varepsilon\|_\infty(\|w_\varepsilon^0\|_\infty+\|w_\varepsilon^{II}\|_\infty)\Big]\cdot O(\varepsilon)\nonumber\\
&=o_\varepsilon(1)\varepsilon O_{m^*}.\nonumber
\end{align}
Moreover, as in the proof of Proposition \ref{prop:Recution}, we get
\begin{align}
A'_{i,3}(x)=O_{m^*}(\|w_\varepsilon^0\|_\infty+\|w_\varepsilon^{II}\|_\infty)e^{v_i(\frac{x-p_{t,\varepsilon}}{\varepsilon_t})+2\ln\frac{1}{\varepsilon_t}}\nonumber
\end{align}
in $B_{5\delta}(p_{t,\varepsilon})$. This gives that
\begin{align}
\int_\Omega A'_{i,3}Z_{i,j,\varepsilon}^*=o_\varepsilon(1)\varepsilon O_{m^*}.\nonumber
\end{align}
For $N'_{i,\varepsilon}$, it holds that
\begin{align}
N'_{i,\varepsilon}\leq C e^{U_{i,\varepsilon}}(\|w_\varepsilon^0\|_\infty+\|w_\varepsilon^{II}\|_\infty)^2.\nonumber
\end{align}
Then, we get
\begin{align}
\int_\Omega N'_{i,\varepsilon}Z_{i,j,\varepsilon}^*=o_\varepsilon(1)\varepsilon O_{m^*}.\nonumber
\end{align}
Together, they proves the result. $\Box$

\subsection{Projections of $S(U_{1,\varepsilon}+w_{1,\varepsilon}^I,\cdots,U_{n,\varepsilon}+w_{n,\varepsilon}^I)$ on $Z_{1,t,\varepsilon}^*(x)$ and on $Z_{2,t,\varepsilon}^*(x)$}\label{Subsection:Projection1}

In this subsection, we compute the projection of
$S(U_{1,\varepsilon}+w_{1,\varepsilon}^I,\cdots,U_{n,\varepsilon}+w_{n,\varepsilon}^I)$
on $Z_{1,t,\varepsilon}^*(x)$ and on $Z_{2,t,\varepsilon}^*(x)$. Here, $w^I_{i,\varepsilon}$ is the $i$-th component of the 1-frequency part of $w_\varepsilon$ defined in Subsection \ref{susbection:1Frequence}.

\begin{proposition}\label{prop:DeltaUZ12}
It holds that for any $t=1,\cdots,N$, we get
\begin{align}
\sum_{i=1}^n\int_\Omega
& \Bigg[-\Delta \Big(U_\varepsilon^i(x)+\sum_{l=1}^n a^{il}w_{l,\varepsilon}^I\Big)-\rho_{i,\varepsilon}\Bigg(\frac{h_i(x)e^{U_{i,\varepsilon}+w_{i,\varepsilon}^I}}{\int_\Omega h_i(x)e^{U_{i,\varepsilon}+w_{i,\varepsilon}^I}}-1\Bigg)\Bigg] Z_{1,t,i,\varepsilon}^*\nonumber\\
&=\frac{1}{N}\sum_{i=1}^n\Big[\partial_{x_1}\ln h_i(p_{t,\varepsilon})+2\pi m_i^*\partial_{x_1} {G^*}(p_{t,\varepsilon};p_{s,\varepsilon})\Big]\rho_i^*\varepsilon_t+
\varepsilon O_{m^*}.\nonumber
\end{align}
and
\begin{align}
\sum_{i=1}^n\int_\Omega
& \Bigg[-\Delta \Big(U_\varepsilon^i(x)+\sum_{l=1}^n a^{il}w_{l,\varepsilon}^I\Big)-\rho_{i,\varepsilon}\Bigg(\frac{h_i(x)e^{U_{i,\varepsilon}+w_{i,\varepsilon}^I}}{\int_\Omega h_i(x)e^{U_{i,\varepsilon}+w_{i,\varepsilon}^I}}-1\Bigg)\Bigg] Z_{2,t,i,\varepsilon}^*\nonumber\\
&=\frac{1}{N}\sum_{i=1}^n\Big[\partial_{x_2}\ln h_i(p_{t,\varepsilon})+2\pi m_i^*\partial_{x_2} {G^*}(p_{t,\varepsilon};p_{s,\varepsilon})\Big]\rho_i^*\varepsilon_t+
\varepsilon O_{m^*}
\nonumber
\end{align}
Here, as it is defined in (\ref{def:Om}),
\begin{align}
O_{m^*}=\left\{
\begin{aligned}
& O(\varepsilon^{m^*-2})&\mbox{ if } m^*<4;\nonumber\\
& O(\varepsilon^{2}\ln\frac{1}{\varepsilon})&\mbox{ if }m^*= 4.\nonumber
\end{aligned}
\right.
\end{align}
\end{proposition}
\noindent{\bf Proof.}
We only prove the first identity since the second is similar.
By the symmetry, we get
\begin{align}
\sum_{i=1}^n\int_\Omega
& \Bigg[-\Delta \Big(U_\varepsilon^i(x)+\sum_{l=1}^n a^{il}w_{l,\varepsilon}^I\Big)-\rho_{i,\varepsilon}\Bigg(\frac{h_i(x)e^{U_{i,\varepsilon}+w_{i,\varepsilon}^I}}{\int_\Omega h_i(x)e^{U_{i,\varepsilon}+w_{i,\varepsilon}^I}}-1\Bigg)\Bigg] Z_{1,t,i,\varepsilon}^*\nonumber\\
&=\sum_{i=1}^n\int_\Omega
\Bigg[-\Delta \Big(U_\varepsilon^i(x)+\sum_{l=1}^n a^{il}w_{l,\varepsilon}^I\Big)-\rho_{i,\varepsilon} \frac{h_i(x)e^{U_{i,\varepsilon}+w_{i,\varepsilon}^I}}{\int_\Omega h_i(x)e^{U_{i,\varepsilon}+w_{i,\varepsilon}^I}}\Bigg]Z_{1,t,i,\varepsilon}^*\nonumber\\
&= \int_{\Omega\backslash B_{\delta_t}(p_{t,\varepsilon})}+ \int_{B_{\delta_t}(p_{t,\varepsilon})}\sum_{i=1}^n
\Bigg[-\Delta \Big(U_\varepsilon^i(x)+\sum_{l=1}^n a^{il}w_{l,\varepsilon}^I\Big)-\rho_{i,\varepsilon} \frac{h_i(x)e^{U_{i,\varepsilon}+w_{i,\varepsilon}^I}}{\int_\Omega h_i(x)e^{U_{i,\varepsilon}+w_{i,\varepsilon}^I}}\Bigg]Z_{1,t,i,\varepsilon}^*\nonumber\\
&=:I_1+I_2.\nonumber
\end{align}
We compute $I_1$ and $I_2$ separately.

By the symmetry and the fact that $w^I_\varepsilon\in E_\varepsilon$, we get for $m^*\leq 4$
\begin{align}\label{e:I1}
I_1=\varepsilon O_{m^*}.
\end{align}

On the other hand, if $m^*<4$,
first we estimate the integral $\int_\Omega h_i(x)e^{U_{i,\varepsilon}+w_{i,\varepsilon}^I}$.
By a similar method as in Lemma \ref{l:Sexpansion}, we get
\begin{align}
h_i(x)e^{U_{i,\varepsilon}+w^I_{i,\varepsilon}}&=h_i(x)e^{\sum_{s=1}^N W_{i,s,\varepsilon_s}(x)+w^I_{i,\varepsilon}} e^{\sum_{s\neq t}v_i\big(\frac{\delta_s}{\varepsilon_s}\big)+2\ln\frac{1}{\varepsilon_s}}e^{\sum_{s\neq t}2\pi m_i^*[G(x,p_{s,\varepsilon})+\frac{1}{2\pi}\ln\delta_s -\gamma(p_{t,\varepsilon},p_{s,\varepsilon})]}\nonumber\\
&=e^{v_i\big(\frac{x-p_{t,\varepsilon}}{\varepsilon_t}\big)+2\ln\frac{1}{\varepsilon_t}}\Big(\Pi_{s\neq t}\varepsilon_s\Big)^{m_i^*-2} e^{2\pi m_i^*\sum_{s\neq t}[G(p_{t,\varepsilon},p_{s,\varepsilon})-\gamma(p_{s,\varepsilon},p_{s,\varepsilon})]}\times\nonumber\\
&\quad\times\Big[h_i(p_{t,\varepsilon})+\nabla_y\Big(h_i(y)e^{2\pi m^*[ {G^*}(y;p_{t,\varepsilon})- {G^*}(p_{t,\varepsilon};p_{t,\varepsilon})]}\Big)\Big|_{y=p_{t,\varepsilon}}\cdot(x-p_{t,\varepsilon})+O(|x-p_{t,\varepsilon}|)^2\Big]\times\nonumber\\
&\quad\times [1+|(1-\chi_t)(W_{j,t,\varepsilon}^{**}-W_{j,t,\varepsilon}^{*})| +O(|(1-\chi_t)(W_{j,t,\varepsilon}^{**}-W_{j,t,\varepsilon}^{*})|^2)]\nonumber\\
&\quad\times [1+w_{i,\varepsilon}^I+O(|w_{i,\varepsilon}^I|^2)].\nonumber
\end{align}
By the symmetry of $w_{i,\varepsilon}^I$ and Proposition \ref{prop:RecutionFiner}, we get for any $t=1,\cdots,N$
\begin{align}
\int_{B_{\delta_t}(p_{t,\varepsilon})}h_i(x)e^{U_{i,\varepsilon}+w^I_{i,\varepsilon}}dx=\frac{\rho_i^* h_i(p_{t,\varepsilon})}{N}e^{\sum_{l\neq t}2\pi m_i^*[G(p_{t,\varepsilon},p_{l,\varepsilon}) -\gamma(p_{l,\varepsilon},p_{l,\varepsilon})]} \Big(\Pi_{l\neq s}\varepsilon_l\Big)^{m_i^* -2}(1+O_{m^*}).\nonumber
\end{align}

Then, due to the symmetry, we get
\begin{align}
I_2&=-\sum_{i=1}^n\frac{\rho_{i,\varepsilon}}{\rho_i^*}\int_{B_{\delta_t}(p_{t,\varepsilon})} e^{v_i\Big(\frac{x-p_{t,\varepsilon}}{\varepsilon_t}\Big)+2\ln\frac{1}{\varepsilon_t}}\Big[e^{\ln\frac{h_i(x)}{h_i(p_{t,\varepsilon})}+2\pi m_i^*[ {G^*}(x;p_{t,\varepsilon})- {G^*}(p_{t,\varepsilon};p_{t,\varepsilon})]+w^I_{i,\varepsilon} +(1-\chi_t)(W_{j,t,\varepsilon}^{**}-W_{j,t,\varepsilon}^{*})}-1\Big]\times \nonumber\\
&\quad\times \partial_{x_1}v_i\Big(\frac{x-p_{t,\varepsilon_t}}{\varepsilon_t}\Big)dx \times(1+O_{m^*})\nonumber\\
&=-\sum_{i=1}^n\int_{B_{\delta_t}(p_{t,\varepsilon})} e^{v_i\Big(\frac{x-p_{t,\varepsilon}}{\varepsilon_t}\Big)+2\ln\frac{1}{\varepsilon_t}}\Big[\Big(\frac{\nabla h_i(p_{t,\varepsilon})}{h_i(p_{t,\varepsilon})}+2\pi m_i^*\nabla_y {G^*}(y;p_{t,\varepsilon})|_{y=p_{t,\varepsilon}}\Big)\cdot(x-p_{t,\varepsilon})\nonumber\\
&\quad+(x-p_{t,\varepsilon}) D(x-p_{t,\varepsilon})^T+O(|x-p_{t,\varepsilon}|^3) +w_{i,\varepsilon}^I+(1-\chi_t)(W_{j,t,\varepsilon}^{**}-W_{j,t,\varepsilon}^*)+O(|w_{i,\varepsilon}^I|^2)\nonumber\\
&\quad+O(|(1-\chi_t)(W_{j,t,\varepsilon}^{**}-W_{j,t,\varepsilon}^*)|^2)\Big]
\partial_{x_1}v_i\Big(\frac{x-p_{t,\varepsilon}}{\varepsilon_t}\Big)dx\cdot (1+o_\varepsilon(1))\nonumber\\
&=-\sum_{i=1}^n\int_{B_\delta(p_{t,\varepsilon})} e^{v_i\Big(\frac{x-p_{t,\varepsilon}}{\varepsilon_t}\Big)+2\ln\frac{1}{\varepsilon_t}}\Big[\Big(\frac{\partial_{x_1} h_i(p_{t,\varepsilon})}{h_i(p_{t,\varepsilon})}+2\pi m_i^*\partial_{x_1}  {G^*}(p_{t,\varepsilon};p_{t,\varepsilon})\Big)(x-p_{t,\varepsilon})_1\nonumber\\
&\quad+O(|x-p_{t,\varepsilon}|^3) \Big]
\partial_{x_1}v_i\Big(\frac{x-p_{t,\varepsilon}}{\varepsilon_t}\Big)dx\cdot (1+o_\varepsilon(1))\nonumber\\
&\quad+o_\varepsilon(1)\varepsilon O_{m^*}.\nonumber
\end{align}
Here, $(x-p_{t,\varepsilon})_1$ is the first component of $x-p_{t,\varepsilon}$.
In the last step, we use the symmetry and the fact that $w_\varepsilon^I\in E_\varepsilon$.
The matrix $D$ is defined as $$D=D^2\exp\Big\{\ln\frac{h_i(x)}{h_i(p_{t,\varepsilon_t})}+2\pi m_i^*[ {G^*}(x;p_{t,\varepsilon})- {G^*}(p_{t,\varepsilon};p_{t,\varepsilon})] +(1-\chi_t)(W_{j,t,\varepsilon}^{**}-W_{j,t,\varepsilon}^*)\Big\}\Big|_{p_{t,\varepsilon_t}}.$$
By a direct computation, we get
\begin{align}
\sum_{i=1}^n\int_{B_{\delta_t}(0)}e^{v_i(\frac{x}{\varepsilon_t})+2\ln\frac{1}{\varepsilon_t}}|x|^3\Big|\partial_{x_1}v_i\Big(\frac{x-p_{t,\varepsilon}}{\varepsilon_t}\Big)\Big|dx\leq C\varepsilon^3_t\sum_{i=1}^n\int_{B_{\frac{\delta_t}{\varepsilon_t}}(0)}e^{v_i(y)}|y|^3 dy=C\varepsilon_t^{m^*-1}\nonumber
\end{align}
and
\begin{align}
\sum_{i=1}^n\int_{B_{\delta_t}(p_{t,\varepsilon_t})} e^{v_i(\frac{x-p_{t,\varepsilon}}{\varepsilon_t})+2\ln\frac{1}{\varepsilon_t}}(x-p_{t,\varepsilon})_1\partial_{x_1}v_i\Big(\frac{x-p_{t,\varepsilon}}{\varepsilon_t}\Big)dx=\varepsilon_t\sum_{i=1}^n\int_{B_\frac{\delta}{\varepsilon_t}(0)}e^{v_i(x)}x_1\partial_{x_1}v_i(x)dx.\nonumber
\end{align}
Notice that
$\frac{\partial}{\partial x_1}\big(x_1 e^{v_i(x)}\big)=e^{v_i(x)}+x_1\partial_{x_1} v_i(x) e^{v_i(x)}$.
We get
$\int_{\partial B_\frac{\delta_t}{\varepsilon_t}(0)}x_1 e^{v_i(x)} dS=\int_{B_\frac{\delta_t}{\varepsilon_t}(0)}e^{v_i(x)}dx+\int_{B_\frac{\delta_t}{\varepsilon_t}(0)}x_1\partial_{x_1}v_i(x)e^{v_i(x)}dx$.
This implies that
\begin{align}
-\int_{B_\frac{\delta_t}{\varepsilon_t}(0)}x_1\partial_{x_1}v_i(x)e^{v_i(x)}dx=\frac{\rho_i^*}{N}+O(\varepsilon^{m^*-2}).\nonumber
\end{align}
Then, we get for $m^*<4$
\begin{align}\label{e:I2}
I_2=\frac{1}{N}\sum_{i=1}^n\Big[\partial_{x_1}\ln h_i(p_{t,\varepsilon})+2\pi m_i^*\partial_{x_1} {G^*}(p_{t,\varepsilon},p_{s,\varepsilon})\Big]\rho_i^*\varepsilon_t+O(\varepsilon^{m^*-1}).
\end{align}

If $m^*=4$ for $i=1,\cdots,n$, notice that
\begin{align}
\sum_{i=1}^n\int_{B_{\delta_t}(0)}e^{v_i(\frac{x}{\varepsilon_t})+2\ln\frac{1}{\varepsilon+t}}|x|^3\Big|\partial_{x_1}v_i\Big(\frac{x}{\varepsilon_t}\Big)\Big|dx
&\leq C\varepsilon^3\sum_{i=1}^n\int_{B_{\frac{\delta_t}{\varepsilon_t}(0)}} e^{v_i(x)}|x|^2 dx\nonumber\\
&=C\varepsilon^3\sum_{i=1}^n\Bigg\{\int_{B_{\frac{\delta_t}{\varepsilon_t}(0)}\backslash B_{\sqrt{\ln\frac{1}{\varepsilon_t}}}(0)} + \int_{B_{\sqrt{\ln\frac{1}{\varepsilon_t}}}(0)}\Bigg\} e^{v_i(x)}|x|^2 dx\nonumber\\
&=O\Big(\varepsilon^3\ln\frac{1}{\varepsilon}\Big)+C\varepsilon^3\int_1^\frac{\delta_t}{\varepsilon}\frac{dr}{r}=O\Big(\varepsilon^3\ln\frac{1}{\varepsilon}\Big).\nonumber
\end{align}

By a similar approach, it holds that
\begin{align}\label{e:I2'}
I_2=\frac{1}{N}\sum_{i=1}^n\Big[\partial_{x_1}\ln h_i(p_{t,\varepsilon})+2\pi m_i^*\partial_{x_1} {G^*}(p_{t,\varepsilon};p_{t,\varepsilon})\Big]\rho_i^*\varepsilon_t+O\Big(\varepsilon^3\ln\frac{1}{\varepsilon}\Big).
\end{align}

Combining (\ref{e:I1}), (\ref{e:I2}) and (\ref{e:I2'}), we get
\begin{align}
\sum_{i=1}^n\int_\Omega
& \Bigg[-\Delta \Big(U_\varepsilon^i(x)+\sum_{l=1}^n a^{il}w_{l,\varepsilon}^I\Big)-\rho_{i,\varepsilon}\Bigg(\frac{h_i(x)e^{U_{i,\varepsilon}+w_{i,\varepsilon}^I}}{\int_\Omega h_i(x)e^{U_{i,\varepsilon}+w_{i,\varepsilon}^I}}-1\Bigg)\Bigg] Z_{1,t,i,\varepsilon}^*\nonumber\\
&=\frac{1}{N}\sum_{i=1}^n\Big[\partial_{x_1}\ln h_i(p_{t,\varepsilon})+2\pi m_i^*\partial_{x_1} {G^*}(p_{t,\varepsilon};p_{t,\varepsilon})\Big]\rho_i^*\varepsilon_t+\varepsilon O_{m^*}.\nonumber
\end{align}

$\Box$

\subsection{Projection of $S(U_{1,\varepsilon}+w_{1,\varepsilon}^I,\cdots,U_{n,\varepsilon}+w_{n,\varepsilon}^I)$ on $Z_{3,t,\varepsilon}^*(x)$: $m^*<4$ Case}\label{Subsection:Projection2}

In this subsection, we compute the projection of
$S(U_{1,\varepsilon}+w_{1,\varepsilon}^I,\cdots,U_{n,\varepsilon}+w_{n,\varepsilon}^I)$
on
$\mbox{span}\{Z_{3,t,\varepsilon}^*(x)\}$ with $m^*<4$.
Here, $w^I_{i,\varepsilon}$ is the $i$-th component of the 1-frequency part of $w_\varepsilon$ defined in Subsection \ref{susbection:1Frequence}.

\begin{lemma}\label{l:ProjZ3}
Assume that $m^*<4$, for any $t=1,\cdots,N$, we get
\begin{align}
&\sum_{i=1}^n\int_{\Omega}\Bigg[-\Delta (U_\varepsilon^i+\sum_{l=1}^n a^{il}w_{l,\varepsilon}^I)-\rho_{i,\varepsilon}\Big( \frac{h_i e^{U_{i,\varepsilon}+w_{i,\varepsilon}^I}}{\int_\Omega h_i e^{U_{i,\varepsilon}+w_{i,\varepsilon}^I}} -1\Big)\Bigg]Z_{3,i,t,\varepsilon}^*\nonumber\\
&=\sum_{i\in\hat{I}}\frac{m^*-2}{N}\Big(D_{i,t}+\sum_{s\neq t}\frac{e^{ {H_{i,s}(p^*)}}}{e^{ {H_{i,t}(p^*)}}} D_{i,s}+o_\delta(1)\Big)\varepsilon_t^{m^*-1}+2\pi\Lambda_{I,N}(\rho_\varepsilon)\varepsilon_t+o_\varepsilon(1)\varepsilon^{m^*-1}.\nonumber
\end{align}
Here, as it is defined in (\ref{NumberD}),
\[
D_{i,s}:=\Bigg[\int_{\Omega_t}\frac{\frac{h_i(x)}{h_i(p_{s}^*)}e^{2\pi m^*( {G^*}(x;p_{s,\varepsilon})- {G^*}(p_s^*;p_{s}^*))} -1}{|x-p_{s}^*|^{m^*}}dx
-\int_{\Omega_s^c}\frac{dx}{|x-p_{s}^*|^{m^*}}\Bigg]e^{I_i}\nonumber
\]
for $i\in\hat{I}=\{i\in I|m^*_i=m^*\}$ and $s=1,\cdots,N$. The sub-domains $\Omega_s\subset\Omega$ satisfy
\begin{itemize}
    \item $\overline{\cup_{s=1}^N\Omega_s}=\overline{\Omega}$;
    \item $\Omega_{s}\cap\Omega_{s'}=\emptyset$ for any $s,s'=1,\cdots,N$ with $s\neq s'$;
    \item $B_{3\delta_s}(p_{s,\varepsilon_s})\in \Omega_s$ for any $s=1,\cdots,N$.
\end{itemize}
\end{lemma}

\noindent{\bf Proof.}
Consider the integration
for any $t=1,\cdots,N$, we get
\begin{align}
\sum_{i=1}^n&\int_{\Omega}\Bigg[-\Delta (U_\varepsilon^i+\sum_{l=1}^n a^{il} w_{l,\varepsilon}^I)-\rho_{i,\varepsilon}\Big( \frac{h_i e^{U_{i,\varepsilon}+w_{i,\varepsilon}^I}}{\int_\Omega h_i e^{U_{i,\varepsilon}+w_{i,\varepsilon}^I}} -1\Big)\Bigg]Z_{3,i,t,\varepsilon}^*\nonumber\\
&=\sum_{i=1}^n\int_{\Omega}\Bigg[-\Delta U_\varepsilon^i-\rho_{i,\varepsilon}\Big( \frac{h_i e^{U_{i,\varepsilon}+w_{i,\varepsilon}^I}}{\int_\Omega h_i e^{U_{i,\varepsilon}+w_{i,\varepsilon}^I}} -1\Big)\Bigg]Z_{3,i,t,\varepsilon}^*\nonumber\\
&=\sum_{i=1}^n\int_{\Omega}\Bigg[-\Delta U_\varepsilon^i-\rho_{i,\varepsilon}\Big( \frac{h_i e^{U_{i,\varepsilon}+w_{i,\varepsilon}^I}}{\int_\Omega h_i e^{U_{i,\varepsilon}+w_{i,\varepsilon}^I}} -1\Big)\Bigg][Z_{3,i,t,\varepsilon}^*(x)-(2-m_i^*)\varepsilon_t]\nonumber\\
&=\sum_{i=1}^n\int_{B_{\delta_t}(p_{t,\varepsilon_t})}\Bigg[-\Delta U_\varepsilon^i-\rho_{i,\varepsilon}\Big( \frac{h_i e^{U_{i,\varepsilon}+w_{i,\varepsilon}^I}}{\int_\Omega h_i e^{U_{i,\varepsilon}+w_{i,\varepsilon}^I}} -1\Big)\Bigg][Z_{3,i,t,\varepsilon}^*(x)-(2-m_i^*)\varepsilon_t].\nonumber
\end{align}
Here, the first equation is due to symmetry and the second to the fact that
\begin{align}
\int_{\Omega}\Bigg[-\Delta U_\varepsilon^i-\rho_{i,\varepsilon}\Big( \frac{h_i e^{U_{i,\varepsilon}+w_{i,\varepsilon}^I}}{\int_\Omega h_i e^{U_{i,\varepsilon}+w_{i,\varepsilon}^I}} -1\Big)\Bigg]dx=0\nonumber
\end{align}
while the third to that $\mbox{supp}\big(Z_{3,i,t,\varepsilon}^*(x)-(2-m_i^*)\varepsilon_t\Big)\subset \overline{B_{\delta_t}(p_{t,\varepsilon})}$.
Then,
\begin{align}
\sum_{i=1}^n&\int_{\Omega}\Bigg[-\Delta U_\varepsilon^i-\rho_{i,\varepsilon}\Big( \frac{h_i e^{U_{i,\varepsilon}+w_{i,\varepsilon}^I}}{\int_\Omega h_i e^{U_{i,\varepsilon}+w_{i,\varepsilon}^I}} -1\Big)\Bigg]Z_{3,i,t,\varepsilon}^*\nonumber\\
&=\Bigg\{\sum_{i=1}^n\int_{B_{3\delta}(p_{t,\varepsilon_t})}\Bigg[-\Delta U_\varepsilon^i-\rho_{i,\varepsilon}\Big( \frac{h_i e^{U_{i,\varepsilon}+w_{i,\varepsilon}^I}}{\int_\Omega h_i e^{U_{i,\varepsilon}+w_{i,\varepsilon}^I}} -1\Big)\Bigg]Z_{3,i,t,\varepsilon}^*(x)\Bigg\}\nonumber\\
&\quad-\Bigg\{\sum_{i=1}^n\int_{B_{3\delta}(p_{t,\varepsilon_t})}\Bigg[-\Delta U_\varepsilon^i-\rho_{i,\varepsilon}\Big( \frac{h_i e^{U_{i,\varepsilon}+w_{i,\varepsilon}^I}}{\int_\Omega h_i e^{U_{i,\varepsilon}+w_{i,\varepsilon}^I}} -1\Big)\Bigg](2-m_i^*)\varepsilon_t\Bigg\}\nonumber\\
&=:II_1-II_2.\nonumber
\end{align}
Now we estimate $II_1$ and $II_2$ separately. To do so, it is necessary to refine the estimate of $\int_\Omega h_i(x)e^{U_{i,\varepsilon}+w_{i,\varepsilon}^I}$.

~

\noindent{\bf Step 1. Refined estimate of $\int_\Omega h_i(x)e^{U_{i,\varepsilon}+w_{i,\varepsilon}^I}$.}

~


For any $x\in B_{\delta_t}(p_{t,\varepsilon})$, it holds that
\begin{align}
&\quad h_i(x)e^{U_{i,\varepsilon}+w_{i,\varepsilon}^I}\nonumber\\
&=h_i(x)e^{W_{i,t}+w_{i,\varepsilon}^I+\sum_{s\neq t}W_{i,s}}\nonumber\\
&=h_i(x)e^{v_i(\frac{x-p_{t,\varepsilon}}{\varepsilon_t})+2\ln\frac{1}{\varepsilon_t}+w_{i,\varepsilon}^I} e^{\sum_{s\neq t}v_i(\frac{\delta_s}{\varepsilon_s})+2\ln\frac{1}{\varepsilon_s}+m_i^*\ln\delta_s} e^{2\pi m_i^*( {G^*}(x;p_{t,\varepsilon})- {G^*}(p_{t,\varepsilon};p_{t,\varepsilon}))} e^{(1-\chi_t)(W^{**}_{i,t,\varepsilon} -W^*_{i,t,\varepsilon})}\nonumber\\
&=h_i(x)e^{v_i(\frac{x-p_{t,\varepsilon}}{\varepsilon_t})+2\ln\frac{1}{\varepsilon_t}+w_{i,\varepsilon}^I}\times e^{2\pi m_i^*( {G^*}(x;p_{t,\varepsilon})- {G^*}(p_{t,\varepsilon};p_{t,\varepsilon}))}\times\nonumber\\
&\quad\times e^{2\pi m_i^*\sum_{s\neq t}[G(p_{t,\varepsilon},p_{s,\varepsilon})-\gamma(p_{s,\varepsilon},p_{s,\varepsilon})]} e^{(1-\chi_t)(W^{**}_{i,t,\varepsilon} -W^*_{i,t,\varepsilon})} e^{NI_i}\Big[\Pi_{s\neq t}\varepsilon_s^{m_i^*-2}\Big(1+\sum_{j=1}^n a_{ij}e^{I_j}\big(\frac{\varepsilon_s}{\delta_s}\big)^{m_j^*-2}+O_{m^*}^2\Big)\Big].\nonumber
\end{align}
Integrating over $B_{\delta_t}(p_{t,\varepsilon})$, we get
\begin{align}
\int_{B_{\delta_t}(p_{t,\varepsilon})}&h_i(x)e^{U_{i,\varepsilon}+w_{i,\varepsilon}^I}dx=e^{2\pi m_i^*\sum_{s\neq t}[G(p_{t,\varepsilon},p_{s,\varepsilon})-\gamma(p_{s,\varepsilon},p_{s,\varepsilon})] }e^{NI_i}\times\nonumber\\
&\times\Big[\Pi_{s\neq t}\varepsilon_s^{m_i^*-2}\Big(1+\sum_{j=1}^n a_{ij}e^{I_j}\big(\frac{\varepsilon}{\delta_s}\big)^{m_j^*-2}+O_{m^*}^2\Big)\Big]\times\nonumber\\
&\times \int_{B_\delta(p_{t,\varepsilon})} h_i(x)e^{v_i(\frac{x-p_{t,\varepsilon}}{\varepsilon_t})+2\ln\frac{1}{\varepsilon_t}}e^{2\pi m_i^*( {G^*}(x;p_{t,\varepsilon})- {G^*}(p_{t,\varepsilon};p_{t,\varepsilon}))+w_{i,\varepsilon}^I+(1-\chi_t)(W^{**}_{i,t,\varepsilon} -W^*_{i,t,\varepsilon})}dx.\nonumber
\end{align}
Here,
\begin{align}
&\int_{B_\delta(p_{t,\varepsilon_t})} h_i(x)e^{v_i(\frac{x-p_{t,\varepsilon}}{\varepsilon_t})+2\ln\frac{1}{\varepsilon_t}} e^{2\pi m_i^*( {G^*}(x;p_{t,\varepsilon})- {G^*}(p_{t,\varepsilon};p_{t,\varepsilon}))}
 e^{w_{i,\varepsilon}^I+(1-\chi_t)(W^{**}_{i,t,\varepsilon}-W^*_{i,t,\varepsilon})}dx\nonumber\\
&=h_i(p_{t,\varepsilon})\int_{B_{\delta_t}(p_{t,\varepsilon})}e^{v_i(\frac{x-p_{t,\varepsilon}}{\varepsilon_t})+2\ln\frac{1}{\varepsilon_t}}
+h_i(p_{t,\varepsilon})\int_{B_{\delta_t}(p_{t,\varepsilon})}e^{v_i(\frac{x-p_{t,\varepsilon}}{\varepsilon_t})+2\ln\frac{1}{\varepsilon_t}}\times\nonumber\\
&\times\Big\{\frac{h_i(x)}{h_i(p_{t,\varepsilon})}e^{2\pi m_i^*( {G^*}(x;p_{t,\varepsilon})- {G^*}(p_{t,\varepsilon};p_{t,\varepsilon}))+w_{i,\varepsilon}^I +(1-\chi_t)(W^{**}_{i,t,\varepsilon}-W^*_{i,t,\varepsilon})}-1\Big\}+O_{m^*}^2\nonumber\\
&=h_i(p_{t,\varepsilon})\frac{\rho_i^*}{N}+2\pi(2-m_i^*)\big(\frac{\varepsilon_t}{\delta_t}\big)^{m_i^*-2}h_i(p_{t,\varepsilon}) \nonumber\\
&+e^{I_i}\int_{B_{\delta_t}(p_{t,\varepsilon})}\frac{\frac{h_i(x)}{h_i(p_{t,\varepsilon})}e^{2\pi m_i^*( {G^*}(x;p_{t,\varepsilon})- {G^*}(p_{t,\varepsilon};p_{t,\varepsilon}))+w_{i,\varepsilon}^I+ (1-\chi_t)(W^{**}_{i,t,\varepsilon}-W^*_{i,t,\varepsilon})}-1}{|x-p_{t,\varepsilon}|^{m_i^*}}dx\cdot\varepsilon_t^{m_i^*-2}+O_{m^*}^2.\nonumber
\end{align}
Here, the integral converges due to Assertion $(W^I_3)$ of Proposition \ref{prop:RecutionFiner}.
Therefore, we get
\begin{align}\label{e:intBpt}
&\int_{B_{\delta_t}(p_{t,\varepsilon})}h_i(x)e^{U_{i,\varepsilon}+w_{i,\varepsilon}^I}dx=e^{2\pi m_i^*\sum_{s\neq t}[G(p_{t,\varepsilon},p_{s,\varepsilon})-\gamma(p_{s,\varepsilon},p_{s,\varepsilon})]}e^{NI_i}\times\nonumber\\
&\times\Big[\Pi_{s\neq t}\varepsilon_s^{m_i^*-2}\Big(1+\sum_{j=1}^n a_{ij}e^{I_j}\big(\frac{\varepsilon_s}{\delta_s}\big)^{m_j^*-2}+O_{m^*}^2\Big)\Big]\times\nonumber\\
&\times\Bigg\{h_i(p_{t,\varepsilon})\frac{\rho_i^*}{N}+2\pi(2-m_i^*)\big(\frac{\varepsilon_t}{\delta_t}\big)^{m_i^*-2}h_i(p_{t,\varepsilon}) \nonumber\\
&+e^{I_i}\int_{B_{\delta_t}(p_{t,\varepsilon})}\frac{\frac{h_i(x)}{h_i(p_{t,\varepsilon})}e^{2\pi m_i^*( {G^*}(x;p_{t,\varepsilon})- {G^*}(p_{t,\varepsilon};p_{t,\varepsilon}))+w_{i,\varepsilon}^I+(1-\chi_t)(W^{**}_{i,t,\varepsilon}-W^*_{i,t,\varepsilon})} -1}{|x-p_{t,\varepsilon}|^{m_i^*}}dx\cdot\varepsilon_t^{m^*-2}+O_{m^*}^2\Bigg\}.
\end{align}
On the other hand, for any $x\in\Omega\backslash \cup_{s=1}^N B_{\delta_s}(p_{s,\varepsilon})$, we get
\begin{align}
&h_i(x)e^{U_{i,\varepsilon}+w_{i,\varepsilon}^I}=h_i(x)e^{\sum_{s=1}^N v_i(\frac{\delta_s}{\varepsilon_s})+2\ln\frac{1}{\varepsilon_s}+m_i^*\ln\delta} e^{2\pi m_i^*\sum_{s=1}^N [G(x,p_{s,\varepsilon})-\gamma(p_{s,\varepsilon},p_{s,\varepsilon})]+w_{i,\varepsilon}^I}.\nonumber
\end{align}
Then, for any $t=1,\cdots,N$ and any $x\in B_{\delta_t}(p_{t,\varepsilon})$,
\begin{align}
h_i(x)e^{U_{i,\varepsilon}+w_{i,\varepsilon}^I}
&= e^{NI_i}\Big[\Pi_{s=1}^N \varepsilon_s^{m_i^*-2}\Big(1+\sum_{j=1}^n a_{ij} e^{I_j}\big(\frac{\varepsilon_s}{\delta}\big)^{m_j^*-2}+O_{m^*}^2\Big)\Big]
e^{2\pi m_i^*\sum_{s\neq t}[G(p_{t,\varepsilon},p_{s,\varepsilon})-\gamma(p_{s,\varepsilon},p_{s,\varepsilon})]}\times\nonumber\\
&\quad\times\frac{e^{2\pi m_i^*( {G^*}(x;p_{t,\varepsilon})- {G^*}(p_{t,\varepsilon};p_{t,\varepsilon}))+w_{i,\varepsilon}^I+(1-\chi_t)(W^{**}_{i,t,\varepsilon}-W^*_{i,t,\varepsilon})}}{|x-p_{t,\varepsilon}|^{m_i^*}}.\nonumber
\end{align}
Select $N$ sub-domains $\Omega_t\subset\Omega$ such that
\begin{itemize}
    \item $\overline{\cup_{s=1}^N\Omega_s}=\overline{\Omega}$;
    \item $\Omega_{s}\cap\Omega_{s'}=\emptyset$ for any $s,s'=1,\cdots,N$ with $s\neq s'$;
    \item $B_{3\delta_s}(p_{s,\varepsilon_s})\in \Omega_s$ for any $s=1,\cdots,N$.
\end{itemize}
Then, we get
\begin{align}\label{e:IntheUWFiner}
&\int_\Omega h_i(x)e^{U_{i,\varepsilon}+w_{i,\varepsilon}^I}dx\nonumber\\
&=\sum_{t=1}^N\Big[\Pi_{s\neq t}\varepsilon_s^{m_i^*-2}\Big(1+\sum_{j=1}^n a_{ij}\big(\frac{\varepsilon_s}{\delta_s}\big)^{m_j^*-2}+O_{m^*}^2\Big)\Big]\times e^{\sum_{s\neq t}2\pi m_i^*[G(p_{t,\varepsilon},p_{s,\varepsilon})-\gamma(p_{s,\varepsilon},p_{s,\varepsilon})]}\times\nonumber\\
&\times\Bigg\{h_i(p_{t,\varepsilon})\frac{\rho_i^*}{N}+\Bigg[e^{I_i}\int_{\Omega_t}\frac{h_i(x)e^{2\pi m_i^*( {G^*}(x;p_{t,\varepsilon})- {G^*}(p_{t,\varepsilon};p_{t,\varepsilon}))+w_{i,\varepsilon}^I+(1-\chi_t)(W^{**}_{i,t,\varepsilon}-W^*_{i,t,\varepsilon})}-h_i(p_{t,\varepsilon})}{|x-p_{t,\varepsilon}|^{m_i^*}}dx \nonumber\\
&-e^{I_i}\int_{\Omega_t^c}\frac{h_i(p_{t,\varepsilon})dx}{|x-p_{t,\varepsilon}|^{m_i^*}}\Bigg]\varepsilon_t^{m_i^*-2}+O_{m^*}^2\Bigg\}.
\end{align}
In the following, for the sake of simplicity, we  denote
\begin{align}
D^\varepsilon_{i,s}:=&\Bigg[\int_{\Omega_t}\frac{\frac{h_i(x)}{h_i(p_{s,\varepsilon})}e^{2\pi m_i^*( {G^*}(x;p_{s,\varepsilon})- {G^*}(p_{s,\varepsilon};p_{s,\varepsilon}))+w_{i,\varepsilon}^I+(1-\chi_s)(W^{**}_{i,s,\varepsilon}-W^*_{i,s,\varepsilon})}-1}{|x-p_{s,\varepsilon}|^{m_i^*}}dx-\int_{\Omega_s^c}\frac{dx}{|x-p_{s,\varepsilon}|^{m_i^*}}\Bigg]e^{I_i}.\nonumber
\end{align}

Based on the above, we have the following claims.

\begin{claim}\label{c:FinerExpansion}
For any $t=1,\cdots,N$ and any $x\in B_{\delta_t}(p_{t,\varepsilon_t})$, we get
\begin{align}
\frac{h_i(x)e^{U_{i,\varepsilon}+w_{i,\varepsilon}^I}}{\int_\Omega h_i(x)e^{U_{i,\varepsilon}+w_{i,\varepsilon}^I}} =\frac{h_i(x)e^{v_i(\frac{x-p_{t,\varepsilon}}{\varepsilon_t}) +2\ln\frac{1}{\varepsilon_t}+w_{i,\varepsilon}^I+ (1-\chi_t)(W_{i,t,\varepsilon}^{**}-W_{i,t,\varepsilon}^{**})}}{\mbox{DENOMINATOR}}.\nonumber
\end{align}
Here,
\begin{align}
\mbox{DENOMINATOR}&=h_i(p_{t,\varepsilon_t})\frac{\rho_i^*}{N}+h_i(p_{t,\varepsilon})D^\varepsilon_{i,t}\varepsilon_t^{m_i^*-2}\nonumber\\
&+\sum_{s\neq t}\big(\frac{\varepsilon_t}{\varepsilon_s}\big)^{m_i^*-2} \Big(1+\sum_{j=1}^n a_{ij}\Big[(\frac{\varepsilon_t}{\delta_t})^{m_j^*-2} - (\frac{\varepsilon_s}{\delta_s})^{m_j^*-2}\Big]+O_{m^*}^2\Big)\times\nonumber\\
&\times e^{ {H_{i,s}(p_\varepsilon)-H_{i,t}(p_\varepsilon )}}  {h_i(p_{t,\varepsilon_t})}\Big(\frac{\rho_i^*}{N}+D^\varepsilon_{i,s}\varepsilon_s^{m_i^*-2}\Big)+O_{m^*}^2.\nonumber
\end{align}
\end{claim}

Under Hypothesis (\ref{ASSUMPTION}), we get
\begin{align}
(\frac{\varepsilon_t}{\delta_t})^{m_j^*-2} - (\frac{\varepsilon_s}{\delta_s})^{m_j^*-2}=O_{m^*}^2\ln\frac{1}{\varepsilon}.\nonumber
\end{align}

\begin{claim}\label{c:Infinitestimal}
Under Hypothesis (\ref{ASSUMPTION1}), i.e.
\begin{align}
(\varepsilon_{t}/\varepsilon_{s})^{m_i^*-2}=e^{ {H_{i,t}(p_\varepsilon)-H_{i,s}(p_\varepsilon )}}+o_\delta(1)\varepsilon^{m^* -2},\nonumber
\end{align}
we get
\begin{align}\label{def:DENOMINATOR}
DENOMINATOR=&h_i(p_{t,\varepsilon})\rho_i^*+\big(h_i(p_{t,\varepsilon})D_{i,t}^\varepsilon +\sum_{s\neq t}e^{ {H_{i,t}(p_\varepsilon)-H_{i,s}(p_\varepsilon )}} {h_i(p_{t,\varepsilon})}D_{i,s}^\varepsilon\nonumber\\
&\quad+o_\delta(1)\big)\varepsilon_t^{m^*_i-2} +o_\varepsilon(1)\varepsilon O_{m^*}.
\end{align}
\end{claim}

~

\noindent{\bf Step 2. Estimating $II_1$ and $II_2$.}

~

By a direct computation, we get
\begin{align}
II_1&=\sum_{i=1}^n\int_{B_{\delta_t}(p_{t,\varepsilon})}\Bigg[-\Delta  U_\varepsilon^i-\rho_{i,\varepsilon}\Big( \frac{h_i e^{U_{i,\varepsilon}+w_{i,\varepsilon}^I}}{\int_\Omega h_i e^{U_{i,\varepsilon}+w_{i,\varepsilon}^I}} -1\Big)\Bigg]Z_{3,i,t,\varepsilon}^*(x)dx\nonumber\\
&=\sum_{i=1}^n\int_{B_{\delta_t}(p_{t,\varepsilon})}\Bigg[-\Delta  U_\varepsilon^i-\rho_{i,\varepsilon}\Big( \frac{h_i(x)e^{v_i(\frac{x-p_{t,\varepsilon}}{\varepsilon_t})+2\ln\frac{1}{\varepsilon_t}+(1-\chi_t)(W^{**}_{i,t,\varepsilon}-W^*_{i,t,\varepsilon})} }{DENOMINATOR}-1\Big)\Bigg]Z_{3,i,t,\varepsilon}^*(x)dx\nonumber
\end{align}
Here, we use Claim \ref{c:FinerExpansion} and $\mbox{DENOMINATOR}$ is defined as in (\ref{def:DENOMINATOR}). Then, we get
\begin{align}\label{e:II1}
II_1&=\sum_{i=1}^n\int_{B_{\delta_t}(p_{t,\varepsilon})}\Bigg[-\Delta U_\varepsilon^i-\rho_{i,\varepsilon}\Big( \frac{h_i e^{U_{i,\varepsilon}+w_{i,\varepsilon}^I}}{\int_\Omega h_i e^{U_{i,\varepsilon}+w_{i,\varepsilon}^I}} -1\Big)\Bigg]Z_{3,i,t,\varepsilon}^*(x)dx\nonumber\\
&=\sum_{i=1}^n\int_{B_{\delta_t}(p_{t,\varepsilon})} \Big[1-\frac{\rho_{i,\varepsilon}h_i(p_{t,\varepsilon})}{\mbox{DENOMINATOR}}\Big] e^{v_i(\frac{x-p_{t,\varepsilon}}{\varepsilon_t})}Z_{3,t,i,\varepsilon}^*(x)+O(\varepsilon^{2m^*-3})\nonumber\\
&+\sum_{i=1}^n\Big[\rho_{i,\varepsilon}-2\pi\sum_{l=1}^N\sum_{j=1}^n a^{ij} m_j^*\Big]\int_{B_{\delta_t}(p_{t,\varepsilon})}Z_{3,t,i,\varepsilon}^*(x)dx\nonumber\\
&+\sum_{i=1}^n\frac{\rho_{i,\varepsilon}h_i(p_{t,\varepsilon})}{\mbox{DENOMINATOR}}\int_{B_{3\delta}(p_{t,\varepsilon})}e^{v_i(\frac{x-p_{t,\varepsilon}}{\varepsilon_t})+2\ln\frac{1}{\varepsilon_t}}\Big[1\nonumber\\
&\quad-\frac{h_i(x)}{h_i(p_{t,\varepsilon})}e^{2\pi m_i^*( {G^*}(x;p_{t,\varepsilon})- {G^*}(p_{t,\varepsilon};p_{t,\varepsilon}))+w_{i,\varepsilon}^I+(1-\chi_t)(W_{i,t,\varepsilon}^{**}- W_{i,t,\varepsilon}^{*})}\Big]Z_{3,i,t,\varepsilon}^*(x)dx\nonumber\\
&=o_\varepsilon(1)\varepsilon^{m^*-1}+\sum_{i=1}^n\Big[\rho_{i,\varepsilon}-2\pi\sum_{l=1}^N\sum_{j=1}^n a^{ij} m_j^*\Big]\int_{B_{\delta_t}(p_{t,\varepsilon})}[(2-m_i^*)\varepsilon_t]dx\nonumber\\
&\quad+\sum_{i=1}^n\frac{\rho_{i,\varepsilon}h_i(p_{t,\varepsilon})}{\mbox{DENOMINATOR}}\int_{B_{3\delta_t}(p_{t,\varepsilon})}e^{v_i(\frac{x-p_{t,\varepsilon_t}}{\varepsilon_t})+2\ln\frac{1}{\varepsilon_t}}\Big[1\nonumber\\
&\quad-\frac{h_i(x)}{h_i(p_{t,\varepsilon})}e^{2\pi m_i^*( {G^*}(x;p_{t,\varepsilon})- {G^*}(p_{t,\varepsilon};p_{t,\varepsilon}))+w_{i,\varepsilon}^I+(1-\chi_t)(W_{i,t,\varepsilon}^{**}- W_{i,t,\varepsilon}^{*})}\Big](2-m_i^*)\varepsilon_t dx.
\end{align}

By a similar method, we get
\begin{align}\label{e:II2}
II_2&=\sum_{i=1}^n\int_{B_{\delta_t}(p_{t,\varepsilon})}\Bigg[-\Delta (U_\varepsilon^i+\sum_{l=1}^n a^{il}w_{l,\varepsilon}^I)-\rho_{i,\varepsilon}\Big( \frac{h_i e^{U_{i,\varepsilon}+w_{i,\varepsilon}^I}}{\int_\Omega h_i e^{U_{i,\varepsilon}+w_{i,\varepsilon}^I}} -1\Big)\Bigg][(2-m_i^*)\varepsilon_t]dx\nonumber\\
&=\sum_{i=1}^n\int_{B_{\delta_t}(p_{t,\varepsilon})}\Big[1-\frac{\rho_{i,\varepsilon}h_i(p_{t,\varepsilon})}{\mbox{DENOMINATOR}}\Big]e^{v_i(\frac{x-p_{t,\varepsilon}}{\varepsilon_t})}[(2-m_i^*)\varepsilon_t]\nonumber\\
&+o_\varepsilon(1)\varepsilon^{m^*-1}+\sum_{i=1}^n\Big[\rho_{i,\varepsilon}-2\pi\sum_{l=1}^N\sum_{j=1}^n a^{ij} m_j^*\Big]\int_{B_{\delta_t}(p_{t,\varepsilon})}[(2-m_i^*)\varepsilon_t]dx\nonumber\\
&\quad+\sum_{i=1}^n\frac{\rho_{i,\varepsilon}h_i(p_{t,\varepsilon})}{\mbox{DENOMINATOR}}\int_{B_{3\delta}(p_{t,\varepsilon})}e^{v_i(\frac{x-p_{t,\varepsilon}}{\varepsilon_t})+2\ln\frac{1}{\varepsilon_t}}\Big[1\nonumber\\
&\quad-\frac{h_i(x)}{h_i(p_{t,\varepsilon})}e^{2\pi m_i^*( {G^*}(x;p_{t,\varepsilon})- {G^*}(p_{t,\varepsilon};p_{t,\varepsilon}))+w_{i,\varepsilon}^I+(1-\chi_t)(W_{i,t,\varepsilon}^{**}- W_{i,t,\varepsilon}^{*})}\Big]\times(2-m_i^*)\varepsilon_t dx.
\end{align}
Combining (\ref{e:II1}) and (\ref{e:II2}), we get
\begin{align}
\sum_{i=1}^n&\int_{\Omega}\Bigg[-\Delta (U_\varepsilon^i+\sum_{l=1}^n a^{il} w_{l,\varepsilon}^I)-\rho_{i,\varepsilon}\Big( \frac{h_i e^{U_{i,\varepsilon}+w_{i,\varepsilon}^I}}{\int_\Omega h_i e^{U_{i,\varepsilon}+w_{i,\varepsilon}^I}} -1\Big)\Bigg]Z_{3,i,t,\varepsilon}^*\nonumber\\
&=II_1 -II_2\nonumber\\
&=-\sum_{i=1}^n\int_{B_{\delta_t}(p_{t,\varepsilon})}\Big[1-\frac{\rho_{i,\varepsilon}h_i(p_{t,\varepsilon})}{\mbox{DENOMINATOR}}\Big]e^{v_i(\frac{x-p_{t,\varepsilon}}{\varepsilon_t})}[(2-m_i^*)\varepsilon_t]dx+o_\varepsilon(1)\varepsilon^{m^*-1}.\nonumber
\end{align}

~

\noindent{\bf Step 3. Completing the proof.}

~

By a direct computation, we get
\begin{align}
\sum_{i=1}^n&\int_{\Omega}\Bigg[-\Delta (U_\varepsilon^i+\sum_{l=1}^n a^{il} w_{l,\varepsilon}^I)-\rho_{i,\varepsilon}\Big( \frac{h_i e^{U_{i,\varepsilon}+w_{i,\varepsilon}^I}}{\int_\Omega h_i e^{U_{i,\varepsilon}+w_{i,\varepsilon}^I}} -1\Big)\Bigg]Z_{3,i,t,\varepsilon}^*\nonumber\\
&=-\sum_{i=1}^n\int_{B_{\delta_t}(p_{t,\varepsilon})}\Big[1-\frac{\rho_{i,\varepsilon}h_i(p_{t,\varepsilon})}{\mbox{DENOMINATOR}}\Big]e^{v_i(\frac{x-p_{t,\varepsilon}}{\varepsilon_t})}[(2-m_i^*)\varepsilon_t]dx+o_\varepsilon(1)\varepsilon^{m_i^*-1}\nonumber\\
&=-\sum_{i=1}^n(2-m_i^*)\varepsilon_t\Big[\frac{\rho_i^*}{N}-2\pi e^{I_i}\big(\frac{\varepsilon_t}{\delta_t}\big)^{m_i^*-2}\Big]\cdot \Big[1-\frac{\rho_{i,\varepsilon}h_i(p_{t,\varepsilon})}{\mbox{DENOMINATOR}}\Big]+o_\varepsilon(1)\varepsilon^{m^*-1}\nonumber\\
&=\Bigg\{-\sum_{i=1}^n(2-m_i^*)\varepsilon_t\Big[\frac{\rho_i^*}{N}-2\pi e^{I_i}\big(\frac{\varepsilon_t}{\delta_t}\big)^{m_i^*-2}\Big]\cdot \frac{\rho_{i,\varepsilon}}{\rho_i^*}\Big[1-\frac{h_i(p_{t,\varepsilon_t})}{\rho_i^*\cdot \mbox{DENOMINATOR}}\Big]\Bigg\}\nonumber\\
&\quad+\Bigg\{-\sum_{i=1}^n(2-m_i^*)\varepsilon_t\Big[\frac{\rho_i^*}{N}-2\pi e^{I_i}\big(\frac{\varepsilon_t}{\delta_t}\big)^{m_i^*-2}\Big]\cdot \frac{\rho_i^* -\rho_{i,\varepsilon}}{\rho_i^*}\Bigg\}+o_\varepsilon(1)\varepsilon^{m_i^*-1}\nonumber\\
&=:II^*_1+II^*_2+o_\varepsilon(1)\varepsilon^{m_i^*-1}.\nonumber
\end{align}
Here, $\mbox{DENOMINATOR}$ is defined as in (\ref{def:DENOMINATOR}).
Now we compute $II^*_1$ and $II^*_2$ separately. First, for $II_2^*$, we get
\begin{align}\label{ineq:II*2}
II_2^*=2\pi\Lambda_{I,N}(\rho_\varepsilon)\varepsilon_t+o_\varepsilon(1)\varepsilon^{m^*-1}.
\end{align}
This is due to Lemma \ref{l:LambdaRho}.
Now we study $II^*_1$.
Notice that for any infinitesimal $o$, we have
\begin{align}
1-\frac{1}{1+o}=o+O(o^2).\nonumber
\end{align}
By Claim \ref{c:Infinitestimal}, we get
\begin{align}
&1-\frac{h_i(p_{t,\varepsilon})}{\rho_i^*\cdot \mbox{DENOMINATOR}}\nonumber\\
&=(\rho_i^*)^{-1}\Big(D_{i,t}^\varepsilon +\sum_{s\neq t} {\frac{e^{(m_i^*-2)H_{i,s}(p_\varepsilon)}}{e^{(m^*_i-2)H_{i,t}(p_\varepsilon)}}} D_{i,s}^\varepsilon+o_\delta(1)\Big)\varepsilon_t^{m^*_i-2} +o_\varepsilon(1) O_{m^*}.\nonumber
\end{align}
Then, we get
\begin{align}\label{ineq:II*1}
II^*_1&=\sum_{i\in\hat{I}}\frac{m_i^*-2}{N}\Big(D_{i,t}^\varepsilon +\sum_{s\neq t} {\frac{e^{H_{i,s}(p_\varepsilon)}}{e^{H_{i,t}(p_\varepsilon)}}} D_{i,s}^\varepsilon+o_\delta(1)\Big)\varepsilon_t^{m^*-2} +o_\varepsilon(1)\varepsilon O_{m^*}.
\end{align}
Here, $\hat{I}=\{i\in I|m^*_i=m^*\}$.
Combining (\ref{ineq:II*2}) and (\ref{ineq:II*1}), we get
\begin{align}
\sum_{i=1}^n&\int_{\Omega}\Bigg[-\Delta (U_\varepsilon^i+\sum_{l=1}^n a^{il} w_{l,\varepsilon}^I)-\rho_{i,\varepsilon}\Big( \frac{h_i e^{U_{i,\varepsilon}+w_{i,\varepsilon}^I}}{\int_\Omega h_i e^{U_{i,\varepsilon}+w_{i,\varepsilon}^I}} -1\Big)\Bigg]Z_{3,i,t,\varepsilon}^*\nonumber\\
&=\sum_{i\in\hat{I}}\frac{m^*-2}{N}\Big(D_{i,t}^\varepsilon +\sum_{s\neq t} {\frac{e^{H_{i,s}(p_\varepsilon)}}{e^{H_{i,t}(p_\varepsilon)}}} D_{i,s}^\varepsilon+o_\delta(1)\Big)\varepsilon_t^{m^*-1}+2\pi\Lambda_{I,N}(\rho_\varepsilon)\varepsilon_t+o_\varepsilon(1)\varepsilon^{m^*-1}\nonumber\\
&=\sum_{i\in\hat{I}}\frac{m^*-2}{N}\Big(D_{i,t}+\sum_{s\neq t} {\frac{e^{H_{i,s}(p^*)}}{e^{H_{i,t}(p^*)}}} D_{i,s}+o_\delta(1)\Big)\varepsilon_t^{m^*-1}+2\pi\Lambda_{I,N}(\rho_\varepsilon)\varepsilon_t+o_\varepsilon(1)\varepsilon^{m^*-1}.\nonumber
\end{align}
Here,
\[
D^\varepsilon_{i,s}:=\Bigg[\int_{\Omega_s}\frac{\frac{h_i(x)}{h_i(p_{s,\varepsilon})}e^{2\pi m_i^*( {G^*}(x;p_{s,\varepsilon})- {G^*}(p_{s,\varepsilon};p_{s,\varepsilon}))}-1}{|x-p_{s,\varepsilon}|^{m_i^*}}dx
-\int_{\Omega_s^c}\frac{dx}{|x-p_{s,\varepsilon}|^{m_i^*}}\Bigg]e^{I_i}
\]
and
\[
D_{i,s}:=\Bigg[\int_{\Omega_t}\frac{\frac{h_i(x)}{h_i(p_{s}^*)}e^{2\pi m_i^*( {G^*}(x;p_s^*)- {G^*}(p_s^*;p_s^*))}-1}{|x-p_{s}^*|^{m_i^*}}dx
-\int_{\Omega_s^c}\frac{dx}{|x-p_{s}^*|^{m_i^*}}\Bigg]e^{I_i}
\]
for $i\in\hat{I}$ and $s=1,\cdots,N$.
This proves the result.
$\Box$

\subsection{Projection of $S(U_{1,\varepsilon}+w_{1,\varepsilon}^I,\cdots,U_{n,\varepsilon}+w_{n,\varepsilon}^I)$ on $Z_{3,t,\varepsilon}^*(x)$: $m^*= 4$ Case}\label{Subsection:Projection3}

In this subsection, we compute the projection of
$S(U_{1,\varepsilon}+w_{1,\varepsilon}^I,\cdots,U_{n,\varepsilon}+w_{n,\varepsilon}^I)$
on
$\mbox{span}\{Z_{3,t,\varepsilon}^*(x)\}$ with $m^*=4$.
Here, $w^I_{i,\varepsilon}$ is the $i$-th component of the 1-frequency part of $w_\varepsilon$ defined in Subsection \ref{susbection:1Frequence}.

\begin{lemma}\label{l:ProjZ3'}
Assume that $m^*=4$ and that
$\sum_{i=1}^n e^{I_i}\sum_{t=1}^N e^{2H_{i,t}(p^*)}h_i(p_{t}^*)L_{i,t}$
is non-zero.
Then,
we get
\begin{align}
\sum_{i=1}^n&\int_{\Omega}\Bigg[-\Delta (U_\varepsilon^i+\sum_{l=1}^n a^{il} w_{l,\varepsilon}^I)-\rho_{i,\varepsilon}\Big( \frac{h_i e^{U_{i,\varepsilon}+w_{i,\varepsilon}^I}}{\int_\Omega h_i e^{U_{i,\varepsilon}+w_{i,\varepsilon}^I}} -1\Big)\Bigg]Z_{3,i,t,\varepsilon}^*\nonumber\\
&=\sum_{i=1}^n\frac{2}{N}\Big(L_{i,t}+\sum_{s\neq t} {\frac{e^{H_{i,s}(p^*)}}{e^{H_{i,t}(p^*)}} }L_{i,s}+o_\delta(1)\Big)\varepsilon_t^3\ln\frac{1}{\varepsilon_t}+2\pi\Lambda_{I,N}(\rho_\varepsilon)\varepsilon_t
+o_\varepsilon(1)\varepsilon_t^3
\ln\frac{1}{\varepsilon_t}.\nonumber
\end{align}
Here,
\begin{align}
L_{i,s}=\Bigg|\frac{\nabla h_i(p_{s}^*)}{h_i(p_{s}^*)}+8\pi\nabla_1  {G^*}(p_s^*;p_s^*)\Bigg|^2 e^{I_i}+\Bigg(\frac{\Delta h_i(p_{s}^*)}{h_i(p_{s}^*)}+8N\pi\Bigg) e^{I_i}\nonumber
\end{align}
for $i=1,\cdots,n$ and $s=1,\cdots,N$ as defined in (\ref{def:Gt}).
\end{lemma}

\noindent{\bf Proof.}
Since a lot of the computations are similar to Lemma \ref{l:ProjZ3}, we only sketch the proof.

~

\noindent{\bf Step 1. Refined estimate on $\int_\Omega h_i(x)e^{U_{i,\varepsilon}+w_{i,\varepsilon}^I}$.}

~

To begin with, notice that
\begin{align}
\int_{B_{\delta_t}(p_{t,\varepsilon_t})}&h_i(x)e^{U_{i,\varepsilon}+w_{i,\varepsilon}^I}dx=e^{8\pi \sum_{s\neq t}[G(p_{t,\varepsilon},p_{s,\varepsilon})-\gamma(p_{s,\varepsilon},p_{s,\varepsilon})] }e^{NI_i}\times\nonumber\\
&\times\Big[\Pi_{s\neq t}\varepsilon_s^2\Big(1+\sum_{j=1}^n a_{ij}e^{I_j}\big(\frac{\varepsilon_s}{\delta_s}\big)^2+O_{m^*}^2\Big)\Big]
 \int_{B_\delta(p_{t,\varepsilon})} h_i(x)e^{v_i(\frac{x-p_{t,\varepsilon}}{\varepsilon_t})+2\ln\frac{1}{\varepsilon_t}}\times\nonumber\\
&\times e^{8\pi( {G^*}(x;p_{t,\varepsilon})- {G^*}(p_{t,\varepsilon};p_{t,\varepsilon}))+w_{i,\varepsilon}^I+(1-\chi_t)(W^{**}_{i,t,\varepsilon} -W^*_{i,t,\varepsilon})}dx.\nonumber
\end{align}
Here,
\begin{align}
\int_{B_\delta(p_{t,\varepsilon})} & h_i(x)e^{v_i(\frac{x-p_{t,\varepsilon}}{\varepsilon_t})+2\ln\frac{1}{\varepsilon_t}} e^{8\pi( {G^*}(x;p_{t,\varepsilon})- {G^*}(p_{t,\varepsilon};p_{t,\varepsilon}))}e^{w_{i,\varepsilon}^I+(1-\chi_t)(W^{**}_{i,t,\varepsilon}-W^*_{i,t,\varepsilon})}dx\nonumber\\
&=h_i(p_{t,\varepsilon})\int_{B_\delta(p_{t,\varepsilon})}e^{v_i(\frac{x-p_{t,\varepsilon}}{\varepsilon_t})+2\ln\frac{1}{\varepsilon_t}}+h_i(p_{t,\varepsilon})\int_{B_\delta(p_{t,\varepsilon})}e^{v_i(\frac{x-p_{t,\varepsilon}}{\varepsilon_t})+2\ln\frac{1}{\varepsilon_t}}\times\nonumber\\
&\times\Big\{\frac{h_i(x)}{h_i(p_{t,\varepsilon_t})}e^{8\pi ( {G^*}(x;p_{t,\varepsilon})- {G^*}(p_{t,\varepsilon};p_{t,\varepsilon}))+w_{i,\varepsilon}^I +(1-\chi_t)(W^{**}_{i,t,\varepsilon}-W^*_{i,t,\varepsilon})}-1\Big\}+O_{m^*}^2\nonumber\\
&=h_i(p_{t,\varepsilon_t})\frac{\rho_i^*}{N}-4\pi\big(\frac{\varepsilon_t}{\delta_t}\big)^2h_i(p_{t,\varepsilon_t}) \nonumber\\
&+\int_{B_\delta(p_{t,\varepsilon})}\Bigg\{\frac{h_i(x)}{h_i(p_{t,\varepsilon})}e^{8\pi( {G^*}(x;p_{t,\varepsilon})- {G^*}(p_{t,\varepsilon};p_{t,\varepsilon}))+w_{i,\varepsilon}^I+ (1-\chi_t)(W^{**}_{i,t,\varepsilon}-W^*_{i,t,\varepsilon})}-1\Big\} e^{v_i(\frac{x-p_{t,\varepsilon}}{\varepsilon_t})+2\ln\frac{1}{\varepsilon}}dx+O_{m^*}^2.\nonumber
\end{align}
By the symmetry and the obvious identity
\begin{align}
\Delta(e^f)=(\Delta f)e^f+|\nabla f|^2 e^f,\nonumber
\end{align}
we get
\begin{align}
&\int_{B_\delta(p_{t,\varepsilon_t})}\Bigg\{\frac{h_i(x)}{h_i(p_{t,\varepsilon_t})}e^{8\pi ( {G^*}(x;p_{t,\varepsilon})- {G^*}(p_{t,\varepsilon};p_{t,\varepsilon}))+w_{i,\varepsilon}^I+ (1-\chi_t)(W^{**}_{i,t,\varepsilon}-W^*_{i,t,\varepsilon})}-1\Big\} e^{v_i(\frac{x-p_{t,\varepsilon}}{\varepsilon_t})+2\ln\frac{1}{\varepsilon}}dx\nonumber\\
&=\frac{e^{-I_i}}{2}\int_{B_{\delta_t}(p_{t,\varepsilon})}\Big( L^\varepsilon_{i,t}|x|^2 +O(|x|^4)\Big)e^{v_i(\frac{x-p_{t,\varepsilon}}{\varepsilon_t})+2\ln\frac{1}{\varepsilon}}dx\nonumber\\
&=\frac{e^{-I_i}}{2} \int_{B_{\delta_t}(p_{t,\varepsilon})}L^\varepsilon_{i,t}|x|^2  e^{v_i(\frac{x-p_{t,\varepsilon}}{\varepsilon_t})+2\ln\frac{1}{\varepsilon}} dx+O(\varepsilon^2)\nonumber
\end{align}
with
\begin{align}
L^\varepsilon_{i,t}=\Bigg|\frac{\nabla h_i(p_{t,\varepsilon})}{h_i(p_{t,\varepsilon})}+8\pi\nabla_1  {G^*}(p_{t,\varepsilon};p_{t,\varepsilon})+\nabla w_{i,\varepsilon}^I(p_{t,\varepsilon})\Bigg|^2 e^{I_i} +\Bigg(\frac{\Delta h_i(p_{t,\varepsilon})}{h_i(p_{t,\varepsilon})}+8N\pi\Bigg) e^{I_i}.\nonumber
\end{align}
By a direct computation, we get
\begin{align}
\int_{B_{\delta_t}(0)}|x|^2 e^{v_i(\frac{x}{\varepsilon_t})+2\ln\frac{1}{\varepsilon_t}}dx&=\varepsilon_t^2\int_{B_\frac{\delta_t}{\varepsilon_t}(0)}|y|^2 e^{v_i(y)}dy=\varepsilon_t^2\Big(C+2\pi e^{I_i}\int_1^{\frac{\delta_t}{\varepsilon_t}}\frac{dr}{r}\Big)\nonumber\\
&=2\pi e^{I_i}\varepsilon_t^2\ln\frac{1}{\varepsilon_t}+O(\varepsilon^2).\nonumber
\end{align}
This implies that
\begin{align}
&\int_{B_\delta(p_{t,\varepsilon})}\Bigg\{\frac{h_i(x)}{h_i(p_{t,\varepsilon})}e^{8\pi ( {G^*}(x;p_{t,\varepsilon})- {G^*}(p_{t,\varepsilon};p_{t,\varepsilon}))+w_{i,\varepsilon}^I+ (1-\chi_t)(W^{**}_{i,t,\varepsilon}-W^*_{i,t,\varepsilon})}-1\Big\}e^{v_i(\frac{x-p_{t,\varepsilon}}{\varepsilon_t})+2\ln\frac{1}{\varepsilon}}dx\nonumber\\
&=\pi  L_{i,\varepsilon}^*\varepsilon_t^2\ln\frac{1}{\varepsilon_t}+O(\varepsilon^2).\nonumber
\end{align}
Then, we get
\begin{align}
\int_{B_{\delta_t}(p_{t,\varepsilon_t})}h_i(x)e^{U_{i,\varepsilon}+w_{i,\varepsilon}^I}dx&=e^{8\pi \sum_{s\neq t}[G(p_{t,\varepsilon},p_{s,\varepsilon})-\gamma(p_{s,\varepsilon},p_{s,\varepsilon})] }e^{NI_i}\times\nonumber\\
&\times\Big[\Pi_{s\neq t}\varepsilon_s^2\Big(1+\sum_{j=1}^n a_{ij}e^{I_j}\big(\frac{\varepsilon_s}{\delta_s}\big)^2+O_{m^*}^2\Big)\Big]\times\nonumber\\
&\times\Big\{h_i(p_{t,\varepsilon})\frac{\rho_i^*}{N}+\pi L_{i,\varepsilon}^*\varepsilon_t^2\ln\frac{1}{\varepsilon_t}+O(\varepsilon^2)\Big\}\nonumber\\
&=e^{8\pi \sum_{s\neq t}[G(p_{t,\varepsilon},p_{s,\varepsilon})-\gamma(p_{s,\varepsilon},p_{s,\varepsilon})] }e^{NI_i}\times\nonumber\\
&\times\Pi_{s\neq t}\varepsilon_s^2\Big\{h_i(p_{t,\varepsilon})\frac{\rho_i^*}{N}+\pi  L_{i,\varepsilon}^*\varepsilon_t^2\ln\frac{1}{\varepsilon_t}+O(\varepsilon^2)\Big\}.\nonumber
\end{align}
Lemma \ref{l:Sexpansion} implies that
\begin{align}\label{e:IntheUWFiner2}
\int_{\Omega}& h_i(x)e^{U_{i,\varepsilon}+w_{i,\varepsilon}^I} dx=e^{8\pi \sum_{s\neq t}[G(p_{t,\varepsilon},p_{s,\varepsilon})-\gamma(p_{s,\varepsilon_s},p_{s,\varepsilon_s})] }e^{NI_i}\Pi_{s\neq t}\varepsilon_s^2\Big\{h_i(p_{t,\varepsilon})\frac{\rho_i^*}{N}+\pi L_{i,\varepsilon}^*(p_\varepsilon)\varepsilon_t^2\ln\frac{1}{\varepsilon_t}+O(\varepsilon^2)\Big\}.
\end{align}

~

\noindent{\bf Step 2. Proof of Lemma \ref{l:ProjZ3'}.}

~

First, notice that
\begin{align}
\int_{B_\delta(0)} &
 e^{v_i\big(\frac{x}{\varepsilon_t}\big)+2\ln\frac{1}{\varepsilon_t}}x_1^2\Big[|x|v'_i\big(\frac{x}{\varepsilon_t}\big)+2\varepsilon_t\Big]dx\nonumber\\
&=\frac{1}{2}\int_{B_\delta(0)}e^{v_i\big(\frac{x}{\varepsilon_t}\big)+2\ln\frac{1}{\varepsilon_t}}|x|^2\Big[|x|v'_i\big(\frac{x}{\varepsilon_t}\big)+2\varepsilon_t\Big]dx\nonumber\\
&=\frac{1}{2}\int_{B_\delta(0)\backslash B_{\varepsilon_t( \ln\frac{1}{\varepsilon_t})^\frac{1}{100}}(0)}+\int_{B_{\varepsilon_t( \ln\frac{1}{\varepsilon_t})^\frac{1}{100}}(0)} e^{v_i\big(\frac{x}{\varepsilon_t}\big)+2\ln\frac{1}{\varepsilon_t}}|x|^2\Big[|x|v'_i\big(\frac{x}{\varepsilon_t}\big)+2\varepsilon_t\Big]dx\nonumber\\
&=-2\pi e^{I_i}\varepsilon_t^3\int_{(\ln\frac{1}{\varepsilon_t})^\frac{1}{100}}^\frac{\delta}{\varepsilon_t}\frac{r^3dr}{r^4}+O\Big(\varepsilon_t^{-1}\int_0^{\varepsilon_t\big(\ln\frac{1}{\varepsilon_t}\big)^\frac{1}{100}}r^3dr\Big)\nonumber\\
&=-2\pi e^{I_i}\varepsilon_t^3\ln\frac{1}{\varepsilon_t}+O\Big(\varepsilon_t^3\ln\Big(\ln\frac{1}{\varepsilon_t}\Big)^\frac{1}{100}+\varepsilon_t^3\Big(\ln\frac{1}{\varepsilon_t}\Big)^\frac{1}{25}\Big)\nonumber\\
&=-2\pi e^{I_i}\varepsilon_t^3\ln\frac{1}{\varepsilon_t}+o_\varepsilon\Big(\varepsilon^3\Big(\ln\frac{1}{\varepsilon}\Big)^\frac{1}{4}\Big).\nonumber
\end{align}
Then, by a similar computation as in {\bf Step 2} and {\bf Step 3} the proof of Lemma \ref{l:ProjZ3}, we find
\begin{align}
\sum_{i=1}^n&\int_{\Omega}\Bigg[-\Delta (U_\varepsilon^i+\sum_{l=1}^n a^{il}w_{l,\varepsilon}^I)-\rho_{i,\varepsilon}\Big( \frac{h_i e^{U_{i,\varepsilon}+w_{i,\varepsilon}^I}}{\int_\Omega h_i e^{U_{i,\varepsilon}+w_{i,\varepsilon}^I}} -1\Big)\Bigg]Z_{3,i,t,\varepsilon}^*\nonumber\\
=&\Bigg\{2\sum_{i=1}^n\varepsilon_t\Big[\frac{\rho_i^*}{N}-2\pi e^{I_i}\big(\frac{\varepsilon_t}{\delta_t}\big)^2\Big]\cdot \frac{\rho_{i,\varepsilon}}{\rho_i^*}\Big[1-\frac{h_i(p_{t,\varepsilon_t})}{\rho_i^*\cdot \mbox{DENOMINATOR}}\Big]\Bigg\}\nonumber\\
&\quad+\Bigg\{2\sum_{i=1}^n\varepsilon_t\Big[\frac{\rho_i^*}{N}-2\pi e^{I_i}\big(\frac{\varepsilon_t}{\delta_t}\big)^2\Big]\cdot \frac{\rho_i^* -\rho_{i,\varepsilon}}{\rho_i^*}\Bigg\}+o_\varepsilon(1)\varepsilon^3\ln\frac{1}{\varepsilon}\nonumber\\
\end{align}
Here,
\begin{align}\label{def:DENOMINATOR1}
\mbox{DENOMINATOR}&=h_i(p_{t,\varepsilon_t})\rho_i^*+\Big(h_i(p_{t,\varepsilon_t})L^\varepsilon_{i,t}+\sum_{s\neq t}  {e^{H_{i,s}(p_\varepsilon)-H_{i,t}(p_\varepsilon )} h_i(p_{t,\varepsilon})}L^\varepsilon_{i,s}\nonumber\\
&\quad+o_\delta(1)\Big)\varepsilon_t^2\ln\frac{1}{\varepsilon_t}+o_\varepsilon(1)\varepsilon_t^2\ln\frac{1}{\varepsilon_t}
\end{align}
and
\begin{align}
L^\varepsilon_{i,s}=\Bigg|\frac{\nabla h_i(p_{s,\varepsilon})}{h_i(p_{s,\varepsilon})}+8\pi\nabla_1  {G^*}(p_{s,\varepsilon};p_{s,\varepsilon})+\nabla w_{i,\varepsilon}^I(p_{s,\varepsilon})\Bigg|^2 e^{I_i}+\Bigg(\frac{\Delta h_i(p_{s,\varepsilon})}{h_i(p_{s,\varepsilon})}+8N\pi\Bigg) e^{I_i}\nonumber
\end{align}
for $s=1,\cdots,N$.
By a similar computation as in Lemma \ref{l:ProjZ3}, we get
\begin{align}
\sum_{i=1}^n&\int_{\Omega}\Bigg[-\Delta (U_\varepsilon^i+\sum_{l=1}^n a^{il} w_{l,\varepsilon}^I)-\rho_{i,\varepsilon}\Big( \frac{h_i e^{U_{i,\varepsilon}+w_{i,\varepsilon}^I}}{\int_\Omega h_i e^{U_{i,\varepsilon}+w_{i,\varepsilon}^I}} -1\Big)\Bigg]Z_{3,i,t,\varepsilon}^*\nonumber\\
&=\sum_{i=1}^n\frac{2}{N}\Big(L^\varepsilon_{i,t}+\sum_{s\neq t} {\frac{e^{H_{i,s}(p_\varepsilon)}}{e^{H_{i,t}(p_\varepsilon)}}} L^\varepsilon_{i,s}+o_\delta(1)\Big)\varepsilon_t^3\ln\frac{1}{\varepsilon_t}+2\pi\Lambda_{I,N}(\rho_\varepsilon)\varepsilon_t+o_\varepsilon(1)\varepsilon_t^3\ln\frac{1}{\varepsilon_t}\nonumber\\
&=\sum_{i=1}^n\frac{2}{N}\Big(L_{i,t}+\sum_{s\neq t} {\frac{e^{H_{i,s}(p^*)}}{e^{H_{i,t}(p^*)}}} L_{i,s}+o_\delta(1)\Big)\varepsilon_t^3\ln\frac{1}{\varepsilon_t}+2\pi\Lambda_{I,N}(\rho_\varepsilon)\varepsilon_t+o_\varepsilon(1)\varepsilon_t^3\ln\frac{1}{\varepsilon_t}.\nonumber
\end{align}
Here,
\begin{align}
L^\varepsilon_{i,s}=\Bigg|\frac{\nabla h_i(p_{s,\varepsilon})}{h_i(p_{s,\varepsilon})}+8\pi\nabla_1  {G^*}(p_{s,\varepsilon};p_{s,\varepsilon})+\nabla w_{i,\varepsilon}^I(p_{s,\varepsilon})\Bigg|^2 e^{I_i}+\Bigg(\frac{\Delta h_i(p_{s,\varepsilon})}{h_i(p_{s,\varepsilon})}+8N\pi\Bigg) e^{I_i}\nonumber
\end{align}
and
\begin{align}\label{def:Lis}
L_{i,s}=\Bigg|\frac{\nabla h_i(p_{s}^*)}{h_i(p_{s}^*)}+8\pi\nabla_1  {G^*}(p_s^*;p_s^*)\Bigg|^2 e^{I_i}+\Bigg(\frac{\Delta h_i(p_{s}^*)}{h_i(p_{s}^*)}+8N\pi\Bigg)e^{I_i}
\end{align}
for $i=1,\cdots,n$ and $s=1,\cdots,N$.
$\Box$

\begin{lemma}\label{l:ProjZ3''}
Assume that $m^*= 4$ and that
$\sum_{i=1}^n \sum_{t}^N  {e^{H_{i,t}(p^*)}} L_{i,t}=0$.
Here, $L_{i,s}$ is defined as in (\ref{def:Lis}). Then,
\begin{align}
&\sum_{i=1}^n\int_{\Omega}\Bigg[-\Delta \Big(U_\varepsilon^i+\sum_{l=1}^n a^{il}w_{l,\varepsilon}^I\Big)-\rho_{i,\varepsilon}\Bigg( \frac{h_i e^{U_{i,\varepsilon}+w_{i,\varepsilon}^I}}{\int_\Omega h_i e^{U_{i,\varepsilon}+w_{i,\varepsilon}^I}} -1\Bigg)\Bigg]Z_{3,i,t,\varepsilon}^*\nonumber\\
&=\sum_{i=1}^n\frac{2}{N}\Big(D_{i,t}+\sum_{s\neq t} {\frac{e^{H_{i,s}(p_\varepsilon)}}{e^{H_{i,t}(p^*)}}} D_{i,s}+o_\delta(1)\Big)\varepsilon_t^{3}+2\pi\Lambda_{I,N}(\rho_\varepsilon)\varepsilon_t+o_\varepsilon(1)\varepsilon^{3}.\nonumber
\end{align}
Here,
\[
D_{i,s}:=\Bigg[\int_{\Omega_t}\frac{\frac{h_i(x)}{h_i(p_{t}^*)}e^{8\pi( {G^*}(x;p_{t,\varepsilon})- {G^*}(p_{t,\varepsilon};p_{t,\varepsilon}))}-1}{|x-p_{t}^*|^{4}}dx
-\int_{\Omega_t^c}\frac{dx}{|x-p_{t}^*|^{4}}\Bigg]e^{I_i}
\]
for $i=1,\cdots,n$ and $s=1,\cdots,N$. The sub-domains $\Omega_s\subset\Omega$ satisfy
\begin{itemize}
    \item $\overline{\cup_{s=1}^N\Omega_s}=\overline{\Omega}$;
    \item $\Omega_{s}\cap\Omega_{s'}=\emptyset$ for any $s,s'=1,\cdots,N$ with $s\neq s'$;
    \item $B_{3\delta_s}(p_{s,\varepsilon_s})\in \Omega_s$ for any $s=1,\cdots,N$.
\end{itemize}
\end{lemma}

\noindent{\bf Proof.}
Since a lot of the computations are similar to  Lemmas \ref{l:ProjZ3} and \ref{l:ProjZ3'}, we only sketch the proof.
It is only need to prove that
\begin{align}
\sum_{i=1}^n&\int_{B_{\delta_t}(p_{t,\varepsilon_t})}\Bigg\{\frac{h_i(x)}{h_i(p_{t,\varepsilon_t})}e^{8\pi ( {G^*}(x;p_{t,\varepsilon})- {G^*}(p_{t,\varepsilon};p_{t,\varepsilon}))+w_{i,\varepsilon}^I+(1-\chi_t)(W^{**}_{i,\varepsilon}-W^*_{i,\varepsilon})}-1\Bigg\} e^{v_i(\frac{x-p_{t,\varepsilon_t}}{\varepsilon_t})+2\ln\frac{1}{\varepsilon_t}}dx\nonumber\\
&=\varepsilon^2\Bigg(e^{I_i}\sum_{i=1}^n\int_{B_{\delta_t}(p_{t,\varepsilon_t})}\frac{\frac{h_i(x)}{h_i(p_{t,\varepsilon_t})}e^{8\pi ( {G^*}(x;p_{t,\varepsilon})- {G^*}(p_{t,\varepsilon};p_{t,\varepsilon}))}-1}{|x-p_t^*|^4}dx+o_\varepsilon(1)\Bigg).\nonumber
\end{align}

Notice that
\begin{align}
\sum_{i=1}^n&\int_{B_{\delta_t}(p_{t,\varepsilon_t})}\Bigg\{\frac{h_i(x)}{h_i(p_{t,\varepsilon})}e^{8\pi ( {G^*}(x;p_{t,\varepsilon})- {G^*}(p_{t,\varepsilon};p_{t,\varepsilon}))+w_{i,\varepsilon}^I+(1-\chi_t)(W^{**}_{i,\varepsilon}-W^*_{i,\varepsilon})}-1\Bigg\} e^{v_i(\frac{x-p_{t,\varepsilon}}{\varepsilon_t})+2\ln\frac{1}{\varepsilon_t}}dx\nonumber\\
&=\sum_{i=1}^n\int_{B_{\delta_t}(p_{t,\varepsilon})\backslash B_{\varepsilon_t^\frac{4}{5}(p_{t,\varepsilon})}}+ \sum_{i=1}^n\int_{B_{\varepsilon_t^\frac{4}{5}}}\Bigg\{\frac{h_i(x)}{h_i(p_{t,\varepsilon})}e^{8\pi ( {G^*}(x;p_{t,\varepsilon})- {G^*}(p_{t,\varepsilon};p_{t,\varepsilon}))} e^{w_{i,\varepsilon}^I+(1-\chi_t)(W^{**}_{i,\varepsilon}-W^*_{i,\varepsilon})}-1\Bigg\} \times\nonumber\\
&\quad\times e^{v_i(\frac{x-p_{t,\varepsilon}}{\varepsilon_t})+2\ln\frac{1}{\varepsilon_t}}dx\nonumber\\
&=:II'_1+II'_2.\nonumber
\end{align}
Here,
\begin{align}
II'_2=O\Big(\varepsilon^{-2}\int_{B_{\varepsilon^\frac{4}{5}}(0)}\varepsilon^2|x|^2 dx\Big)=O(\varepsilon^\frac{16}{5})=o_\varepsilon(1)\varepsilon^3\nonumber
\end{align}
and
\begin{align}
II'_1&=\sum_{i=1}^n\int_{B_{\delta_t}(p_{t,\varepsilon})\backslash B_{\varepsilon^\frac{4}{5}}(p_{t,\varepsilon})}\Bigg\{\frac{h_i(x)}{h_i(p_{t,\varepsilon})}e^{8\pi ( {G^*}(x;p_{t,\varepsilon})- {G^*}(p_{t,\varepsilon};p_{t,\varepsilon}))+w_{i,\varepsilon}^I+(1-\chi_t)(W^{**}_{i,\varepsilon}-W^*_{i,\varepsilon})}-1\Bigg\}e^{v_i(\frac{x-p_{t,\varepsilon}}{\varepsilon_t})+2\ln\frac{1}{\varepsilon_t}}dx\nonumber\\
&=\sum_{i=1}^n\int_{B_{\delta_t}(p_{t,\varepsilon})\backslash B_{\varepsilon^\frac{4}{5}}(p_{t,\varepsilon})}\Bigg\{\frac{h_i(x)}{h_i(p_{t,\varepsilon})}e^{8\pi ( {G^*}(x;p_{t,\varepsilon})- {G^*}(p_{t,\varepsilon};p_{t,\varepsilon}))}-1\Bigg\} e^{v_i(\frac{x-p_{t,\varepsilon}}{\varepsilon_t})+2\ln\frac{1}{\varepsilon_t}}dx\nonumber\\
&\quad\quad+O\Bigg(\varepsilon^2\int_{B_{\delta_t}(0)\backslash B_{\varepsilon^\frac{4}{5}(0)}}\frac{\varepsilon\ln\frac{1}{\varepsilon}}{|x|^2}dx\Bigg)\nonumber\\
&=\sum_{i=1}^n\int_{B_{\delta_t}(p_{t,\varepsilon})}\Bigg\{\frac{h_i(x)}{h_i(p_{t,\varepsilon_t})}e^{8\pi ( {G^*}(x;p_{t,\varepsilon})- {G^*}(p_{t,\varepsilon};p_{t,\varepsilon}))}-1\Bigg\} e^{v_i(\frac{x-p^*_{t}}{\varepsilon_t})+2\ln\frac{1}{\varepsilon_t}}dx\nonumber\\
&\quad\quad+O\Bigg(\int_{B_{\delta_t}(0)\backslash B_{\varepsilon^\frac{4}{5}(0)}}\frac{\varepsilon^2\ln\frac{1}{\varepsilon}}{|x|^3}dx\Bigg)\varepsilon^2+O\Bigg(\varepsilon^2\int_{B_{\delta_t}(0)\backslash B_{\varepsilon^\frac{4}{5}(0)}}\frac{\varepsilon\ln\frac{1}{\varepsilon}}{|x|^2}dx\Bigg)\nonumber\\
&=\sum_{i=1}^n\int_{B_{\delta_t}(p_{t,\varepsilon})\backslash B_{\varepsilon^\frac{4}{5}}(p_{t,\varepsilon_t})}\Bigg\{\frac{h_i(x)}{h_i(p_{t,\varepsilon_t})}e^{8\pi ( {G^*}(x;p_{t,\varepsilon})- {G^*}(p_{t,\varepsilon};p_{t,\varepsilon}))} -1\Big\} e^{v_i(\frac{x-p^*_{t}}{\varepsilon_t})+2\ln\frac{1}{\varepsilon_t}}dx\nonumber\\
&\quad\quad+O\Bigg(\int_{B_{\delta_t}(0)\backslash B_{\varepsilon^\frac{4}{5}(0)}}\frac{\varepsilon^2\ln\frac{1}{\varepsilon}}{|x|^3}dx\Bigg)\varepsilon^2+o_\varepsilon(1)\varepsilon^2\nonumber\\
&=\varepsilon^2\Bigg(e^{I_i}\sum_{i=1}^n\int_{B_{\delta_t}(p_{t,\varepsilon})}\frac{\frac{h_i(x)}{h_i(p_{t,\varepsilon_t})}e^{8\pi ( {G^*}(x;p_{t,\varepsilon})- {G^*}(p_{t,\varepsilon};p_{t,\varepsilon}))}-1}{|x-p_t^*|^4}dx+o_\varepsilon(1)\Bigg).\nonumber
\end{align}
Notice that in the third equation, the first error is due to Hypothesis (\ref{ASSUMPTION}) and the difference between
\[
\Bigg\{\frac{h_i(x)}{h_i(p_{t,\varepsilon_t})}e^{8\pi ( {G^*}(x;p_{t,\varepsilon})- {G^*}(p_{t,\varepsilon};p_{t,\varepsilon}))}-1\Bigg\}
e^{v_i(\frac{x-p_{t,\varepsilon}}{\varepsilon_t})+2\ln\frac{1}{\varepsilon_t}}
\]
and
\[
\Bigg\{\frac{h_i(x)}{h_i(p_{t,\varepsilon_t})}e^{8\pi ( {G^*}(x;p_{t,\varepsilon})- {G^*}(p_{t,\varepsilon};p_{t,\varepsilon}))w_{i,\varepsilon}^I+(1-\chi_t)(W^{**}_{i,\varepsilon}-W^*_{i,\varepsilon})}-1\Bigg\}
e^{v_i(\frac{x-p^*_{t}}{\varepsilon_t})+2\ln\frac{1}{\varepsilon_t}}.\nonumber
\]
Notice the number
\begin{align}
\sum_{i\in\hat{I}}\frac{m^*-2}{N}\Big(D_{i,t}+\sum_{s\neq t} {\frac{e^{H_{i,s}(p^*)}}{e^{H_{i,t}(p^*)}}} D_{i,s}\Big)
\end{align}
is finite
due to
the assumption
$\sum_{i=1}^n \sum_{t=1}^N  {e^{H_{i,t}(p^*)}} L_{i,t}=0$.
Here,
\[
D_{i,s}:=\Bigg[\int_{\Omega_t}\frac{\frac{h_i(x)}{h_i(p_{s}^*)}e^{8\pi( {G^*}(x;p_{s}^*)- {G^*}(p_s^*;p_s^*))}-1}{|x-p_{s}^*|^{4}}dx
-\int_{\Omega_t^c}\frac{dx}{|x-p_{s}^*|^{4}}\Bigg]e^{I_i}
\]
and
\[
L_{i,t}:=\Bigg|\frac{\nabla h_i(p_{t}^*)}{h_i(p_{t}^*)}+8\pi\nabla_1  {G^*}(p_t^*;p_t^*)\Bigg|^2 e^{I_i}+\Bigg(\frac{\Delta h_i(p_{t}^*)}{h_i(p_{t}^*)}+8N\pi\Bigg)e^{I_i}.
\]
Together, they prove the result.
$\Box$

\subsection{Proof of Theorem \ref{t:MAIN}}\label{Subsection:proof}
Now we prove Theorem \ref{t:MAIN}.
We only verify the Case $(1)$ of Theorem \ref{t:MAIN} since Cases $(2)$ and $(3)$ are similar.
To do this, we will
\begin{itemize}
    \item Find a genuine solution to Problem (\ref{e:001});
    \item Verify Hypotheses (\ref{ASSUMPTION0}), (\ref{ASSUMPTION}) and (\ref{ASSUMPTION1});
    \item Verify (\ref{SE1}), (\ref{SE2}) and (\ref{SE3}).
\end{itemize}
By Lemma \ref{l:ERROR} and Proposition \ref{prop:DeltaUZ12}, for $j=1,2$ and $t=1,\cdots,N$, we get
\begin{align}
\sum_{i=1}^n\int_\Omega
&\Bigg(-\Delta (U_\varepsilon^i+w_\varepsilon^i)-\rho_{i,\varepsilon}\Bigg(\frac{h_i(x)e^{U_{i,\varepsilon}+w_{i,\varepsilon}}}{\int_\Omega h_i(x)e^{U_{i,\varepsilon}+w_{i,\varepsilon}}}-1\Bigg)\Bigg)Z_{j,t,i,\varepsilon}^*\nonumber\\
&=\sum_{i=1}^n\int_\Omega
\Bigg(-\Delta \Big(U_\varepsilon^i+\sum_{j=1}^n a^{ij} w^I_{j,\varepsilon}\Big)-\rho_{i,\varepsilon}\Bigg(\frac{h_i(x)e^{U_{i,\varepsilon}+w^i_{i,\varepsilon}}}{\int_\Omega h_i(x)e^{U_{i,\varepsilon}+w^I_{i,\varepsilon}}}-1\Bigg)\Bigg)Z_{j,t,i,\varepsilon}^*+R_{j,t,\varepsilon}
\nonumber\\
&=
\frac{1}{N}\sum_{i=1}^n\Big[\partial_{x_j}\ln h_i(p_{t,\varepsilon})+2\pi m_i^*\partial_{x_j} {G^*}(p_{t,\varepsilon};p_{t,\varepsilon})\Big]\rho_i^*\varepsilon_t+\varepsilon O_{m^*}.
\nonumber
\end{align}
By the non-degeneracy of the critical point $(p_1^*,\cdots,p_N^*)$ to $\sum_{t=1}^N\Big\{\sum_{i=1}^n\Big[\ln h_i(x_t)+2\pi m_i^* {G^*}(x;p_t^*)\Big]\rho_i^*\Big\}$
such that $p_{t,\varepsilon_t}\to p_t^*$ as $\varepsilon\to0$ and for any $t=1,\cdots,N$
\begin{align}\label{e:NablaGOm} \sum_{i=1}^n\Big(\nabla\ln h_i(p_{t,\varepsilon})+2\pi m_i^*\nabla  {G^*}(p_{t,\varepsilon};p_{t,\varepsilon})\Big)=O_{m^*}.
\end{align}
Therefore, $(p_{1,\varepsilon},\cdots,p_{N,\varepsilon})$ satisfies Hypothesis (\ref{ASSUMPTION}).

By a similar idea, we get
\begin{align}
\sum_{i=1}^n\int_{\Omega}&\Bigg[-\Delta (U_\varepsilon^i+w_\varepsilon^i)-\rho_{i,\varepsilon}\Big( \frac{h_i e^{U_{i,\varepsilon}+w_\varepsilon^i}}{\int_\Omega h_i e^{U_{i,\varepsilon}+w_\varepsilon^i}}-1\Big)\Bigg]Z_{3,t,i,\varepsilon}^*\nonumber\\
&=\sum_{i=1}^n\int_{\Omega}\Bigg[-\Delta (U_\varepsilon^i+\sum_{j=1}^n a^{ij}w_{j,\varepsilon}^I) -\rho_{i,\varepsilon}\Big( \frac{h_i e^{U_{i,\varepsilon}+w_{i,\varepsilon}^I}}{\int_\Omega h_i e^{U_{i,\varepsilon}+w_{i,\varepsilon}^I}} -1\Big)\Bigg]Z_{3,i,\varepsilon}^*+R_{3,t,\varepsilon}\nonumber\\
&=\sum_{i=1}^n\frac{m_i^*-2}{N}\Big(D_{i,t}+\sum_{s\neq t} {\frac{e^{H_{i,s}(p^*)}}{e^{H_{i,t}(p^*)}}} D_{i,s}+o_\delta(1)\Big)\varepsilon_t^{m^*_i-1}+2\pi\Lambda_{I,N}(\rho_\varepsilon)\varepsilon_t+o_\varepsilon(1)\varepsilon^{m^*-1}.\nonumber
\end{align}
with
\begin{align}
D_{i,s}:=\Bigg[\int_{\Omega_t}\frac{\frac{h_i(x)}{h_i(p_{t}^*)}e^{2\pi m_i^*( {G^*}(x;p_t^*)- {G^*}(p_t^*;p_t^*))}-1}{|x-p_{t}^*|^{m_i^*}}dx-\int_{\Omega_t^c}\frac{dx}{|x-p_{t}^*|^{m_i^*}}\Bigg]e^{I_i}\nonumber
\end{align}
for $s=1,\cdots,N$.
Due to the geometry of the hypersurface
\begin{align}
\Gamma_{N}=\{\rho=(\rho_1,\cdots,\rho_N)\in(\mathbb{R}_+)^N|\Lambda_{I,N}(\rho)=0\},\nonumber
\end{align}
we can find a sequence $\rho_\varepsilon=(\rho_{1,\varepsilon},\cdots,\rho_{n,\varepsilon})$ such that
\begin{itemize}
    \item $\rho_{i,\varepsilon}\to\rho_i^*$ as $\varepsilon\to0$;
    \item $\frac{\rho_{j,\varepsilon} -\rho_j^*}{\rho_{i,\varepsilon}-\rho_i^*}\sim 1$.
\end{itemize}
The existence of such a sequence is due to Sard's theorem \cite{HirschBook1976}.
Therefore, we can find a sequence $(\varepsilon_1,\cdots,\varepsilon_N)$ such that
\begin{align}\label{e:LambdaD}
\Lambda(\rho_\varepsilon)=-\sum_{i=1}^n\frac{m_i^*-2}{2\pi N}\Big(D_{i,t}+\sum_{s\neq t} {\frac{e^{H_{i,s}(p^*)}}{e^{H_{i,t}(p^*)}}} D_{i,s}+o_\delta(1)\Big)\varepsilon_t^{m^*_i-2}.
\end{align}
This implies that Hypothesis (\ref{ASSUMPTION}) is satisfied.

Now we are in the position to verify Hypothesis (\ref{ASSUMPTION1}).
The idea to verify Hypothesis (\ref{ASSUMPTION1}) is to proceed the blow-up analysis for the constructed solution $(U_{1,\varepsilon}+w_{1,\varepsilon},\cdots,U_{n,\varepsilon}+w_{n,\varepsilon})$. Denote for any $t=1,\cdots,N$, choose
$p_{t,\varepsilon}$ to be the blow-up point. By the construction of the solution $(U_{1,\varepsilon}+w_{1,\varepsilon},\cdots,U_{n,\varepsilon}+w_{n,\varepsilon})$, we find that
\begin{align}
M_{t,\varepsilon}^i:=U_{i,\varepsilon}(p_{t,\varepsilon})+w_{i,\varepsilon}(p_{t,\varepsilon})=U_i(0)+\ln\frac{1}{2\varepsilon_t}+O(|w^0_{i,\varepsilon}(p_{t,\varepsilon_t})|)=U_i(0)+\ln\frac{1}{2\varepsilon'_t}.\nonumber
\end{align}
Here, $\varepsilon'_t$, the scale in the view point of blow up.
Then, we get
$\varepsilon'_t=\varepsilon_t e^{O(|w^0_{i,\varepsilon}(p_{t,\varepsilon_t})|)}$.
Repeating the computation as in Subsection \ref{subsection:RefindeReduction1} (See also \cite[pp. 2608-2609]{LinZhang2013})
, we get $O(|w^0_{i,\varepsilon}(p_{t,\varepsilon_t})|)=O(\delta^2)\varepsilon^2$.
Therefore, we get for any $t,s=1,\cdots,N$ with $t\neq s$
\begin{align}
\big(\frac{\varepsilon'_t}{\varepsilon'_s}\big)^{m_i^*-2}=\big(\frac{\varepsilon_t}{\varepsilon_s}\big)^{m_i^*-2}\Big(1+O(\delta^2)\varepsilon^2\Big)\nonumber
\end{align}
for suitable small $\theta>0$.
By a standard approach as in \cite[Theorem 2.6]{GuZhangArxiv}, we find that
\begin{align}\label{e:varepsilonH1}
\big(\frac{\varepsilon'_t}{\varepsilon'_s}\big)^{m_i^*-2}=\frac{e^{(m_i^*-2)H_{i,t}(p_\varepsilon)}}{e^{(m_i^*-2)H_{i,s}(p_\varepsilon)}}+O(\delta^{4-m^*})(\varepsilon'_t)^{m^*_i-2}.
\end{align}
For the case $m^*=4$, a similar argument implies that
\begin{align}\label{e:varepsilonH2}
\big(\frac{\varepsilon'_t}{\varepsilon'_s}\big)^2=\frac{e^{2H_{i,t}(p_\varepsilon)}}{e^{2H_{i,s}(p_\varepsilon)}}+O(\delta^2)(\varepsilon'_t)^2.
\end{align}
This is sufficient to verify Hypothesis (\ref{ASSUMPTION1}).

Now we verify (\ref{SE1}), (\ref{SE2}) and (\ref{SE3}). To do this, we only need to check that
\begin{align}\label{Ineq:pp}
|p_{t,\varepsilon}-p_{t,max}|=
O(\varepsilon^2).
\end{align}
Here, $p_{t,max}$ attends
$\max_{i=1,\cdots,n}\max_{x\in B_{5\delta_t}(p_{t,\varepsilon})}\big(U_{i,\varepsilon}(x)+w_{i,\varepsilon(x)}\big)$.
Indeed, this follows Assertion ($W_3^I$) of Proposition \ref{prop:RecutionFiner} and a similar argument as \cite[Lemma 4.1]{LinZhang2013}. Consequently,
(\ref{e:NablaGOm}), (\ref{e:LambdaD}), (\ref{e:varepsilonH1}) and (\ref{e:varepsilonH2}) imply (\ref{SE1}), (\ref{SE2}) and (\ref{SE3}).
The proof is completed.
    $\Box$

\appendix

\section{The Invertibility of a Linear Operator}\label{App:Invertbility}

In this part, we prove the invertibility of the operator $Q_\varepsilon L_\varepsilon$. A standard method can be found in \cite{LinYan2013,Huang2019} and the references therein. Recall that in (\ref{ProjectionQ}) and (\ref{OperatorL}) we define
\begin{align}\label{ProjectionQ}
Q_\varepsilon(u)=u&-\sum_{j=1}^2\sum_{t=1}^N r_{j,t}\Big(\chi_t\partial_{x_j}v_1(\frac{x-p_{t,\varepsilon_t}}{\varepsilon_t})\sum_{s=1}^Ne^{v_1(\frac{x-p_{s,\varepsilon}}{\varepsilon_s})+2\ln\frac{1}{\varepsilon_s}},\cdots,\nonumber\\
&\quad\quad\chi_t\partial_{x_j}v_n(\frac{x-p_{t,\varepsilon}}{\varepsilon_t})\sum_{s=1}^Ne^{v_n(\frac{x-p_{s,\varepsilon}}{\varepsilon_s})+2\ln\frac{1}{\varepsilon_s}}\Big)\nonumber\\
&-\sum_{t=1}^N r_{3,t}\Big(Z_{3,t,1,\varepsilon }^*\sum_{s=1}^N\chi_s e^{v_1(\frac{x-p_{s,\varepsilon}}{\varepsilon_s})+2\ln\frac{1}{\varepsilon_s}},\cdots,Z_{3,t,n,\varepsilon }^*\sum_{s=1}^N\chi_s e^{v_n(\frac{x-p_{s,\varepsilon}}{\varepsilon_s})+2\ln\frac{1}{\varepsilon_s}}\Big)\nonumber\\
&-\Big(s_1\sum_{t=1}^N\chi_s e^{v_1(\frac{x-p_{s,\varepsilon}}{\varepsilon_s})+2\ln\frac{1}{\varepsilon_s}},\cdots,s_n\sum_{s=1}^N\chi_t e^{v_n(\frac{x-p_{s,\varepsilon}}{\varepsilon_s})+2\ln\frac{1}{\varepsilon_s}}\Big)\nonumber
\end{align}
and
\begin{align}
L_{i,\varepsilon}(w_\varepsilon^1,\cdots,w_\varepsilon^n) = -\Delta w_\varepsilon^i-K_{i,\varepsilon}(x)w_{i,\varepsilon} +\frac{K_{i,\varepsilon}(x)}{\rho_i^*}\int_\Omega K_{i,\varepsilon}(x) w_{i,\varepsilon}.\nonumber
\end{align}
Moreover, let
\begin{align}
L_\varepsilon=(L_{1,\varepsilon},\cdots,L_{n,\varepsilon}).\nonumber
\end{align}

This appendix aims to prove the following theorem.
\begin{theorem}\label{t:invertibility}
Under Hypothesis (\ref{ASSUMPTION1}), there exists a positive constant $C$ independence of $\varepsilon$ such that
\begin{equation}\label{Ineq:Invert}
\|w\|_{\mathbb{X}_\varepsilon}+\|w\|_{L^\infty(\Omega)}\leq C\ln\frac{1}{\varepsilon} \|Q_\varepsilon L_\varepsilon w\|_{\mathbb{Y}_\varepsilon}.
\end{equation}
for any $w\in E_\varepsilon$. Furthermore, $Q_\varepsilon L_\varepsilon$ is an isomorphism from $E_\varepsilon$ to $F_\varepsilon$.
\end{theorem}

\begin{remark}
Under Hypothesis (\ref{ASSUMPTION1}), $\varepsilon_t$ are equivalence to each other for $t=1,\cdots,N$.
\end{remark}

For any $w_\varepsilon\in E_\varepsilon$, there exist constants $r_{j,t}$ and $s_i$ for
$j=1,2,3$, $t=1,\cdots,N$ and $i=1,\cdots,n$ such that
\begin{align}\label{e:expansion}
L_{i,\varepsilon} w_\varepsilon &=\phi_{i,\varepsilon}+ \sum_{j=1}^2 \sum_{t=1}^N r_{j,t}\partial_{x_j} v_i( \frac{x-p_{t,\varepsilon}}{\varepsilon_t})\cdot\sum_{s=1}^N\chi_s e^{v_i(\frac{x-p_{s,\varepsilon}}{\varepsilon_s})+2\ln\frac{1}{\varepsilon_s}}\nonumber\\
&\quad+\sum_{t=1}^N r_{3,t}Z_{3,t,i,\varepsilon}^*\cdot\sum_{s=1}^N\chi_s e^{v_i(\frac{x-p_{s,\varepsilon}}{\varepsilon_s})+2\ln\frac{1}{\varepsilon_s}}+s_i \sum_{s=1}^N\chi_s e^{v_i(\frac{x-p_{s,\varepsilon}}{\varepsilon_s})+2\ln\frac{1}{\varepsilon_s}}
\end{align}
for $i=1,\cdots,n$.

We begin by estimating the constants $r_{j,t}$ and $s_i$ for
$j=1,2,3$, $t=1,\cdots,N$ and $i=1,\cdots,n$.

\begin{lemma}\label{l:constantsrrs}
Assume that $m^*\leq4$.
For any $j=1,2$ and any $t=1,\cdots,N$, we have
\begin{align}\label{Ineq:r12}
r_{j,t}=O(\|\phi_\varepsilon\|_{\mathbb{Y}_\varepsilon})+O(\varepsilon^{m^*-2})\max_{i=1,\cdots,n}\sup_{y\in B_{2\delta}(p_{t,\varepsilon_t})}|w_{i,\varepsilon}(y)|,
\end{align}
For $t=1,\cdots, N$ and any $i=1,\cdots,n$, we get
\begin{align}\label{Ineq:r3}
r_{3,t}=O(\varepsilon^{m^*-3})\sup|w_\varepsilon|+\frac{1}{\varepsilon}O(\|\phi_\varepsilon\|_{\mathbb{Y}_\varepsilon})
\end{align}
and
\begin{align}\label{Ineq:si}
s_i=O(\varepsilon^{m^*-2})\sup|w_\varepsilon|+O(\|\phi_\varepsilon\|_{\mathbb{Y}_\varepsilon}).
\end{align}
\end{lemma}
\noindent{\bf Proof.}
First, we prove (\ref{Ineq:r12}). Multiplying by  $\chi_t\partial_{x_j}v_i(\frac{x-p_t}{\varepsilon_t})$ on both sides of (\ref{e:expansion}), integrating over $\Omega$ and summing in $i=1,\cdots,n$, we get

\begin{align}\label{e:Expansion1}
\sum_{i=1}^n & \int_\Omega (L_{i,\varepsilon}w_\varepsilon)\chi_t \partial_{x_j} v_i(\frac{x-p_{t,\varepsilon}}{\varepsilon_t})\nonumber\\
&=\sum_{i=1}^n\int_\Omega \phi_{i,\varepsilon}\chi_t \partial_{x_j}v_i(\frac{x-p_{t,\varepsilon}}{\varepsilon_t})+ r_{j,t}\sum_{i=1}^n\int_\Omega |\partial_{x_j} v_i(\frac{x-p_{t,\varepsilon}}{\varepsilon_t})|^2 \chi_t \Big[ \sum_{s=1}^N\chi_s e^{v_i(\frac{x-p_{\varepsilon_s}}{\varepsilon_s})+2\ln\frac{1}{\varepsilon_s}}\Big].
\end{align}
By a direct computation, we get
\begin{align}
\mbox{LHS of (\ref{e:Expansion1})}=
\sum_{i=1}^n \int_\Omega (L_{i,\varepsilon}w_\varepsilon)\chi_t \partial_{x_j} v_i(\frac{x-p_{t,\varepsilon}}{\varepsilon_t})=\sum_{i=1}^n  \int_\Omega (L_{i,\varepsilon}Z_{j,t\varepsilon}^*)w_{i,\varepsilon}.\nonumber
\end{align}
Here,
\begin{align}
L_{i,\varepsilon}Z_{j,t,\varepsilon}^*
=-\Delta\Big[\sum_{l=1}^n a^{il}\chi_t\partial_{x_j} v_l(\frac{x-p_{t,\varepsilon}}{\varepsilon_t})\Big]-K_{i,\varepsilon}(x)\chi_t\partial_{x_j}v_i(\frac{x-p_{t,\varepsilon}}{\varepsilon_t})\nonumber
\end{align}
and we use the symmetry. Moreover, by a pointwise computation, we get
\begin{align}
L_{i,\varepsilon}Z_{j,t,\varepsilon}^*=-\Delta\Big[\sum_{l=1}^n a^{il}\chi_t\partial_{x_j} v_l(\frac{x-p_{t,\varepsilon}}{\varepsilon_t})\Big]-K_{i,\varepsilon}(x)\chi_t\partial_{x_j}v_i(\frac{x-p_{t,\varepsilon}}{\varepsilon_t})=O(\varepsilon^{m^*-2}).\nonumber
\end{align}
Therefore,
\begin{align}
\mbox{LHS of (\ref{e:Expansion1})}=O(\varepsilon^{m^*-2})\max_{i=1,\cdots,n}\sup_{y\in B_{2\delta}(p_{t,\varepsilon_t})}|w_{i,\varepsilon}(y)|.\nonumber
\end{align}
On the other hand, we have
\begin{align}
&\mbox{RHS of (\ref{e:Expansion1})}\nonumber\\
&\quad\quad =\sum_{i=1}^n\int_\Omega \phi_{i,\varepsilon}\chi_t \partial_{x_j}v_i(\frac{x-p_{t,\varepsilon}}{\varepsilon_t})+ r_{j,t}\sum_{i=1}^n\int_\Omega |\partial_{x_j} v_i(\frac{x-p_{t,\varepsilon}}{\varepsilon_t})|^2 \chi_t \Big[ \sum_{s=1}^N\chi_s e^{v_i(\frac{x-p_{s,\varepsilon}}{\varepsilon_s})+2\ln\frac{1}{\varepsilon_s}}\Big].\nonumber
\end{align}
Here,
\begin{align}
\sum_{i=1}^n  \int_\Omega \phi_{i,\varepsilon}\chi_t \partial_{x_j}v_i(\frac{x-p_{t,\varepsilon}}{\varepsilon_t})
&=O\Big(\varepsilon_t^2\sum_{i=1}^n\int_{B_\frac{2\delta_t}{\varepsilon_t}}|h_{i,\varepsilon}(p_{t,\varepsilon}+\varepsilon_t y)|\cdot|\partial_{x_j}v_i(y)|dy\Big)=O(\|\phi_\varepsilon\|_{\mathbb{Y}_\varepsilon})\nonumber
\end{align}
and
\begin{align}
\sum_{i=1}^n\int_\Omega|\partial_{x_j}v_i(\frac{x-p_{t,\varepsilon}}{\varepsilon_t})|^2\chi_t^3 e^{v_i(\frac{x-p_{t,\varepsilon}}{\varepsilon_t})+2\ln\frac{1}{\varepsilon_t}}=\sum_{i=1}^n\int_{\mathbb{R}^2}|\partial_{x_1}v_i|^2e^{v_i}dx+o_\varepsilon(1).\nonumber
\end{align}
Therefore,
\begin{align}
\mbox{RHS of (\ref{e:Expansion1})}=O(\|\phi_\varepsilon\|_{\mathbb{Y}_\varepsilon})+\Big(\sum_{i=1}^n\int_{\mathbb{R}^2}|\partial_{x_1}v_i|^2e^{v_i}dx+o_\varepsilon(1)\Big)r_{j,t}.\nonumber
\end{align}
Plugging these into (\ref{e:Expansion1}), we get
\begin{align}
r_{j,t}=O(\|\phi_\varepsilon\|_{\mathbb{Y}_\varepsilon})+O(\varepsilon^{m^*-2})\max_{i=1,\cdots,n}\sup_{y\in B_{2\delta}(p_{t,\varepsilon_t})}|w_{i,\varepsilon}(y)|\nonumber
\end{align}
for any $j=1,2$ and any $t=1,\cdots,N$.
This proves (\ref{Ineq:r12}).

Now we prove (\ref{Ineq:r3}) and (\ref{Ineq:si}).
Multiplying by  $Z_{3,t,i,\varepsilon}^*$ on both sides of (\ref{e:expansion}), integrating over $\Omega$ and summing in $i=1,\cdots,n$, we
get
\begin{align}\label{e:Expansion2}
\sum_{i=1}^n & \int_\Omega (L_{i,\varepsilon}w_\varepsilon) Z_{3,t,i,\varepsilon}^*dx\nonumber\\
&=\sum_{i=1}^n\int_\Omega \phi_{i,\varepsilon}
Z_{3,t,i,\varepsilon}^* +r_{3,t}\Big[\sum_{i=1}^n\int_{\Omega}\chi_t(Z_{3,t,i,\varepsilon}^*)^2e^{v_i(\frac{x-p_{t,\varepsilon}}{\varepsilon_t})+2\ln\frac{1}{\varepsilon_t}}dx\nonumber\\
&+ \sum_{s\neq t}\sum_{i=1}^n\int_\Omega\chi_s(2-m_i^*)^2\varepsilon_t^2 e^{v_i(\frac{x-p_{s,\varepsilon}}{\varepsilon_s})+2\ln\frac{1}{\varepsilon_s}}dx\Big]\nonumber\\
&+\sum_{s\neq t}r_{3,s}\sum_{i=1}^n\Big[\sum_{l\neq s,l\neq t}\int_\Omega\chi_l e^{v_i(\frac{x-p_{l,\varepsilon}}{\varepsilon_l})+2\ln\frac{1}{\varepsilon_l}}(2-m_i^*)^2\varepsilon_s\varepsilon_t dx\nonumber\\
& +\int_\Omega \chi_t e^{v_i(\frac{x-p_{t,\varepsilon}}{\varepsilon_t})+2\ln\frac{1}{\varepsilon_t}}(2-m_i^*)\varepsilon_s Z_{3,i,t,\varepsilon}^* dx+ \int_\Omega \chi_s e^{v_i(\frac{x-p_{s,\varepsilon}}{\varepsilon_s})+2\ln\frac{1}{\varepsilon_s}}(2-m_i^*)\varepsilon_t Z_{3,i,s,\varepsilon}^* dx\Big]\nonumber\\
&+\sum_{i=1}^n s_i\Big
[\int_\Omega \chi_t Z_{3,t,i,\varepsilon}^* e^{v_i(\frac{x-p_{t,\varepsilon}}{\varepsilon_t})+2\ln\frac{1}{\varepsilon_t}}dx+\sum_{s\neq t}\int_\Omega \chi_s(2-m_i^*)\varepsilon_t e^{v_i(\frac{x-p_{s,\varepsilon}}{\varepsilon_s})+2\ln\frac{1}{\varepsilon_s}}dx\Big].
\end{align}
Here,
\begin{align}
\mbox{LHS of (\ref{e:Expansion2})}= \sum_{i=1}^n  \int_\Omega (L_{i,\varepsilon}w_\varepsilon) Z_{3,t,i,\varepsilon}^*dx = \sum_{i=1}^n  \int_\Omega (L_{i,\varepsilon} Z_{3,t,\varepsilon}^*)w_{i,\varepsilon} dx \nonumber
\end{align}
Notice that
\begin{align}
L_{i,\varepsilon} Z_{3,t,\varepsilon}^* &=-\Delta\Big[\sum_{l=1}^n a^{il}Z_{3,i,t,\varepsilon}^*\Big]-K_{i,\varepsilon}(x)Z_{3,i,t,\varepsilon}^*+\frac{K_{i,\varepsilon}(x)}{\rho_i^*}\int_\Omega K_{i,\varepsilon}(x)Z_{3,i,t,\varepsilon}^* dx=O(\varepsilon^{m^*-1}).\nonumber
\end{align}
This implies that
\begin{align}
|\mbox{LHS of (\ref{e:Expansion2})}|=O(\varepsilon^{m^*-1})\max_{i=1,\cdots,n}\sup_{y\in B_{2\delta}(p_{t,\varepsilon_t})}|w_{i,\varepsilon}(y)|\nonumber
\end{align}
On the other hand, we get
\begin{align}
\mbox{RHS of (\ref{e:Expansion2})}&=\Bigg\{\sum_{i=1}^n\int_\Omega \phi_{i,\varepsilon}
Z_{3,t,i,\varepsilon}^* \Bigg\}+\Bigg\{r_{3,t}\Big[\sum_{i=1}^n\int_{\Omega}\chi_t(Z_{3,t,i,\varepsilon}^*)^2e^{v_i(\frac{x-p_{t,\varepsilon}}{\varepsilon_t})+2\ln\frac{1}{\varepsilon_t}}dx\nonumber\\
&+ \sum_{s\neq t}\sum_{i=1}^n\int_\Omega\chi_s(2-m_i^*)^2\varepsilon_t^2 e^{v_i(\frac{x-p_{s,\varepsilon}}{\varepsilon_s})+2\ln\frac{1}{\varepsilon_s}}dx\Big]\nonumber\\
&+\sum_{s\neq t}r_{3,s}\sum_{i=1}^n\Big[\sum_{l\neq s,l\neq t}\int_\Omega\chi_l e^{v_i(\frac{x-p_{l,\varepsilon}}{\varepsilon_l})+2\ln\frac{1}{\varepsilon_l}}(2-m_i^*)^2\varepsilon_s\varepsilon_t dx\nonumber\\
& +\int_\Omega \chi_t e^{v_i(\frac{x-p_{t,\varepsilon}}{\varepsilon_t})+2\ln\frac{1}{\varepsilon_t}}(2-m_i^*)\varepsilon_s Z_{3,i,t,\varepsilon}^* dx+ \int_\Omega \chi_s e^{v_i(\frac{x-p_{s,\varepsilon}}{\varepsilon_s})+2\ln\frac{1}{\varepsilon_s}}(2-m_i^*)\varepsilon_t Z_{3,i,s,\varepsilon}^* dx\Big]
\Bigg\}\nonumber\\
&+\Bigg\{\sum_{i=1}^n s_i\Big
[\int_\Omega \chi_t Z_{3,t,i,\varepsilon}^* e^{v_i(\frac{x-p_{t,\varepsilon}}{\varepsilon_t})+2\ln\frac{1}{\varepsilon_t}}dx+\sum_{s\neq t}\int_\Omega \chi_s(2-m_i^*)\varepsilon_t e^{v_i(\frac{x-p_{s,\varepsilon}}{\varepsilon_s})+2\ln\frac{1}{\varepsilon_s}}dx\Big]\Bigg\}\nonumber\\
&=:A_1+A_2+A_3.
\end{align}
Here,
\begin{align}\label{Ineq:hY}
A_1&=\sum_{i=1}^n\int_\Omega \phi_{i,\varepsilon}
Z_{3,t,i,\varepsilon}^* \leq C\varepsilon\sum_{i=1}^n\int_\Omega |\phi_{i,\varepsilon} |\nonumber\\
&\leq C\varepsilon_t^3\sum_{i=1}^n \sum_{t=1}^N\int_{B_{\frac{2\delta_t}{\varepsilon_t}(p_{t,\varepsilon})}}|\phi_{i,\varepsilon}(y)|dx+
C\varepsilon\sum_{i=1}^n\Big(\int_{\Omega\backslash\cup_{t=1}^N B_{2\delta_t}(p_{t,\varepsilon_t})}|\phi_{i,\varepsilon}|^2\Big)^\frac{1}{2}\nonumber\\
&=O(\varepsilon\|\phi_{i,\varepsilon}\|_{\mathbb{Y}_\varepsilon})
\end{align}
and
\begin{align}
A_2=&r_{3,t}\Big[\sum_{i=1}^n\int_{\Omega}\chi_t(Z_{3,t,i,\varepsilon}^*)^2e^{v_i(\frac{x-p_{t,\varepsilon}}{\varepsilon_t})+2\ln\frac{1}{\varepsilon_t}}dx+ \sum_{s\neq t}\sum_{i=1}^n\int_\Omega\chi_s(2-m_i^*)^2\varepsilon_t^2 e^{v_i(\frac{x-p_{s,\varepsilon}}{\varepsilon_s})+2\ln\frac{1}{\varepsilon_s}}dx\Big]\nonumber\\
&+\sum_{s\neq t}r_{3,s}\sum_{i=1}^n\Big[\sum_{l\neq s,l\neq t}\int_\Omega\chi_l e^{v_i(\frac{x-p_{l,\varepsilon}}{\varepsilon_l})+2\ln\frac{1}{\varepsilon_l}}(2-m_i^*)\varepsilon_s\varepsilon_t dx+\int_\Omega \chi_t e^{v_i(\frac{x-p_{t,\varepsilon}}{\varepsilon_t})+2\ln\frac{1}{\varepsilon_t}}(2-m_i^*)^2\varepsilon_s Z_{3,i,t,\varepsilon}^* dx\nonumber\\
&+ \int_\Omega \chi_s e^{v_i(\frac{x-p_{s,\varepsilon}}{\varepsilon_s})+2\ln\frac{1}{\varepsilon_s}}(2-m_i^*)\varepsilon_t Z_{3,i,s,\varepsilon}^* dx\Big].\nonumber
\end{align}
Notice that
\begin{align}
\sum_{i=1}^n \int_\Omega\chi_t(Z_{3,t,i,\varepsilon}^*)^2 e^{v_i(\frac{x-p_{t,\varepsilon_t}}{\varepsilon_t})+2\ln\frac{1}{\varepsilon_t}}dx=\Big[\sum_{i=1}^n\int_{\mathbb{R}^2}(\nabla v_i\cdot x+2)^2 e^{v_i}dx+o_\varepsilon(1)\Big]\varepsilon_t^2,\nonumber
\end{align}
\begin{align}
\sum_{i=1}^n\int_\Omega \chi_s(2-\frac{m_i^*}{2\pi})^2\varepsilon_t^2 e^{v_i(\frac{x-p_{s,\varepsilon}}{\varepsilon_s})+2\ln\frac{1}{\varepsilon_s}}dx =\Big[\sum_{i=1}^n(2-m_i^*)^2\int_{\mathbb{R}^2}e^{v_i}dx +o_\varepsilon(1)\Big]\varepsilon_t^2,\nonumber
\end{align}
\begin{align}
\sum_{i=1}^n\int_\Omega \chi_l(2-\frac{m_i^*}{2\pi})^2\varepsilon_s\varepsilon_t e^{v_i(\frac{x-p_{l,\varepsilon}}{\varepsilon_l})+2\ln\frac{1}{\varepsilon_l}}dx =\Big[\sum_{i=1}^n(2-m_i^*)^2\int_{\mathbb{R}^2}e^{v_i}dx +o_\varepsilon(1)\Big]\varepsilon_t\varepsilon_s\nonumber
\end{align}
and
\begin{align}
\sum_{i=1}^n\int_\Omega \chi_t e^{v_i(\frac{x-p_{t,\varepsilon}}{\varepsilon_t})+2\ln\frac{1}{\varepsilon_t}}(2-m_i^*)\varepsilon_s Z_{3,i,t,\varepsilon}^* dx&\sim \sum_{i=1}^n\int_\Omega \chi_s e^{v_i(\frac{x-p_{s,\varepsilon}}{\varepsilon_s})+2\ln\frac{1}{\varepsilon_s}}(2-m_i^*)\varepsilon_s Z_{3,i,s,\varepsilon}^* dx\nonumber\\
&=O(\varepsilon^{m^*}).\nonumber
\end{align}
Therefore, we get
\begin{align}
A_2&=r_{3,t}\Big[\Big(\sum_{i=1}^n\int_{\mathbb{R}^2}(\nabla v_i\cdot x+2)^2 e^{v_i}dx\Big)\varepsilon_t^2+\sum_{s\neq t}\Big(\sum_{i=1}^n(2-m^*_i)^2\int_{\mathbb{R}^2}e^{v_i}dx\Big)\varepsilon_t^2 +o_\varepsilon(1)\varepsilon^2\Big]\nonumber\\
&\quad+\sum_{s\neq t}r_{3,s}\Big[\sum_{l\neq s,l\neq t}\Big(\sum_{i=1}^n(2-m_i^*)^2\int_{\mathbb{R}^2}e^{v_i}dx\Big)\varepsilon_s\varepsilon_t+o_\varepsilon(1)\varepsilon^2 \Big].\nonumber
\end{align}
On the other hand, we get
\begin{align}
A_3=\sum_{i=1}^ns_i\Big[(N-1)(2-m_i^*)\int_{\mathbb{R}^2}e^{v_i}dx \varepsilon_t+o_\varepsilon(1)\varepsilon \Big].\nonumber
\end{align}
Combining them together, we get
\begin{align}
\mbox{RHS of (\ref{e:Expansion2})}&= O(\varepsilon\|h_\varepsilon\|_{\mathbb{Y}_\varepsilon})\nonumber\\
&\quad+r_{3,t}\Big[\Big(\sum_{i=1}^n\int_{\mathbb{R}^2}(\nabla v_i\cdot x+2)^2 e^{v_i}dx\Big)\varepsilon_t^2+\sum_{s\neq t}\Big(\sum_{i=1}^n(2-m^*_i)^2\int_{\mathbb{R}^2}e^{v_i}dx\Big)\varepsilon_t^2+o_\varepsilon(1)\varepsilon^2 \Big]\nonumber\\
&\quad+\sum_{s\neq t}r_{3,s}\Big[\sum_{l\neq s,l\neq t}\Big(\sum_{i=1}^n(2-m_i^*)^2\int_{\mathbb{R}^2}e^{v_i}dx\Big)\varepsilon_s\varepsilon_t+o_\varepsilon(1)\varepsilon^2 \Big]\nonumber\\
&\quad+\sum_{i=1}^ns_i\Big[(N-1)(2-m_i^*)\int_{\mathbb{R}^2}e^{v_i}dx \varepsilon_t+o_\varepsilon(1)\varepsilon \Big].\nonumber
\end{align}
Plugging this into (\ref{e:Expansion2}), we get for any $t=1,\cdots,N$
\begin{align}
r_{3,t}&\Big[\Big(\sum_{i=1}^n\int_{\mathbb{R}^2}(\nabla v_i\cdot x+2)^2 e^{v_i}dx\Big)\varepsilon_t^2+\sum_{s\neq t}\Big(\sum_{i=1}^n(2-m^*_i)^2\int_{\mathbb{R}^2}e^{v_i}dx\Big)\varepsilon_s\varepsilon_t+o_\varepsilon(1)\varepsilon^2 \Big]\nonumber\\
&\quad+\sum_{s\neq t}r_{3,s}\Big[\sum_{l\neq s,l\neq t}\Big(\sum_{i=1}^n(2-m_i^*)^2\int_{\mathbb{R}^2}e^{v_i}dx\Big)\varepsilon_s\varepsilon_t+o_\varepsilon(1)\varepsilon^2 \Big]\nonumber\\
&\quad+\sum_{i=1}^ns_i\Big[(N-1)(2-m_i^*)\int_{\mathbb{R}^2}e^{v_i}dx \varepsilon_t+o_\varepsilon(1)\varepsilon \Big]=O(\varepsilon\|\phi_\varepsilon\|_{\mathbb{Y}_\varepsilon})+ O(\varepsilon^{m^*-1})\max_{i=1,\cdots,n}\sup_{y\in B_{2\delta}(p_{t,\varepsilon_t})}|w_{i,\varepsilon}(y)|.\nonumber
\end{align}
In other words, we get for any $t=1,\cdots,N$
\begin{align}\label{e:CoeffientEquation1}
\varepsilon_t r_{3,t}&\Big[\Big(\sum_{i=1}^n\int_{\mathbb{R}^2}(\nabla v_i\cdot x+2)^2 e^{v_i}dx\Big)+\sum_{s\neq t}\Big(\sum_{i=1}^n(2-m^*_i)^2\int_{\mathbb{R}^2}e^{v_i}dx\Big)+o_\varepsilon(1)\Big]\nonumber\\
&\quad+\sum_{s\neq t}\varepsilon_s r_{3,s}\Big[\sum_{l\neq s,l\neq t}\Big(\sum_{i=1}^n(2-m_i^*)^2\int_{\mathbb{R}^2}e^{v_i}dx\Big)+o_\varepsilon(1)\Big]\nonumber\\
&\quad+\sum_{i=1}^ns_i\Big[(N-1)(2-m_i^*)\int_{\mathbb{R}^2}e^{v_i}dx +o_\varepsilon(1)\Big]=O(\|\phi_\varepsilon\|_{\mathbb{Y}_\varepsilon})+ O(\varepsilon^{m^*-2}) \max_{i=1,\cdots,n}\sup_{y\in B_{2\delta}(p_{t,\varepsilon_t})}|w_{i,\varepsilon}(y)|.
\end{align}

Denote
$e_i=(0,\cdots,0,1,0,\cdots,0)$.
Here, $e_i$ is a $n$-dimensional vector with $i$-th element equals to $1$ and the rest are $0$. Denote $e_{i,j}$ by the $j$-th element of $e_i$. Integrating (\ref{e:expansion}) over $\Omega$, we get
\begin{align}\label{e:Expansion3}
\int_\Omega L_{i,\varepsilon} w_\varepsilon dx &=\int_\Omega \phi_{i,\varepsilon} dx+\sum_{t=1}^N r_{3,t}\int_\Omega Z_{3,t,i,\varepsilon}^*\cdot\sum_{s=1}^N\chi_s e^{v_i(\frac{x-p_{s,\varepsilon}}{\varepsilon_s})+2\ln\frac{1}{\varepsilon_s}}dx+s_i \sum_{s=1}^N\int_\Omega \chi_s e^{v_i(\frac{x-p_{s,\varepsilon}}{\varepsilon_s})+2\ln\frac{1}{\varepsilon_s}}dx.
\end{align}
Here,
\begin{align}
\mbox{LHS of (\ref{e:Expansion3})} =\int_\Omega L_{i,\varepsilon} w_\varepsilon dx=\sum_{j=1}^n\int_\Omega (L_{j,\varepsilon}e_i)w_{j,\varepsilon}dx.\nonumber
\end{align}
Notice that
$L_{j,\varepsilon}e_i=0$
for $j=1,\cdots,n$ and $j\neq i$. And
\begin{align}
L_{i,\varepsilon} e_i=- K_{i,\varepsilon}(x)+\frac{\int_\Omega K_{i,\varepsilon}(x) dx}{\rho_i^*}h_{i,\varepsilon}(x).\nonumber
\end{align}
This implies that
\begin{align}
\mbox{LHS of (\ref{e:Expansion3})} =O(\varepsilon^{m^*-2})\max_{i=1,\cdots,n}\sup_{x\in B}|w_{i,\varepsilon}|.\nonumber
\end{align}
On the other hand, by a direct computation, we get
\begin{align}\label{Ineq:phiY}
\int_\Omega \phi_{i,\varepsilon}dx=O(\|\phi_\varepsilon\|_{\mathbb{Y}_\varepsilon}),
\end{align}
\begin{align}
\int_\Omega Z_{3,t,i,\varepsilon}^*\Big[\sum_{s=1}^N\chi_s e^{v_i(\frac{x-p_{s,\varepsilon}}{\varepsilon_s})+2\ln\frac{1}{\varepsilon_s}}\Big]dx=&\int_\Omega Z_{3,t,i,\varepsilon}^* e^{v_i(\frac{x-p_{t,\varepsilon}}{\varepsilon_t})+2\ln\frac{1}{\varepsilon_t}}\chi_tdx\nonumber\\
&+\sum_{s\neq t}(2-m_i^*)\varepsilon_t\int_\Omega\chi_s e^{v_i(\frac{x-p_{s,\varepsilon}}{\varepsilon_s})+2\ln\frac{1}{\varepsilon_s}}dx\nonumber\\
&=O(\varepsilon^{m^*-1})+(2-m_i^*)\varepsilon_t(N-1)\int_{\mathbb{R}^2}e^{v_i}dx,\nonumber
\end{align}
\begin{align}
\sum_{s=1}^N\int_\Omega \chi_s e^{v_i(\frac{x-p_{s,\varepsilon}}{\varepsilon_s})+2\ln\frac{1}{\varepsilon_s}}dx=N\int_{\mathbb{R}^2}e^{v_i}dx+O(\varepsilon^{m^*-2}).\nonumber
\end{align}
Here, the first inequality is obtained in a similar way as in (\ref{Ineq:hY}).
Plugging these into (\ref{e:Expansion3}), we get
\begin{align}
\sum_{t=1}^N r_{3,t}\Big[\varepsilon_t(N-1)(2-m_i^*) \int_{\mathbb{R}^2}e^{v_i}+o_\varepsilon(1)\varepsilon\Big]+s_i\Big(N\int_{\mathbb{R}^2}e^{v_i}dx+o_\varepsilon(1)\Big) =O(\varepsilon^{m^*-2})\sup|w_\varepsilon|+O(\|\phi_\varepsilon\|_{\mathbb{Y}_\varepsilon}).\nonumber
\end{align}
In other words,
\begin{align}\label{e:CoeffientEquation2}
\sum_{t=1}^N (\varepsilon_t r_{3,t})\Big[(N-1)(2-m_i^*) \int_{\mathbb{R}^2}e^{v_i}+o_\varepsilon(1)\Big]+s_i\Big(N\int_{\mathbb{R}^2}e^{v_i}dx+o_\varepsilon(1)\Big) =O(\varepsilon^{m^*-2})\sup|w_\varepsilon|+O(\|\phi_\varepsilon\|_{\mathbb{Y}_\varepsilon})
\end{align}
for $i=1,\cdots,n$.

Rewriting (\ref{e:CoeffientEquation1}) and (\ref{e:CoeffientEquation2}) into matrix form, we denote
\begin{align}
X=(\varepsilon_1 r_{3,1},\cdots,\varepsilon_N r_{3N},s_1,\cdots,s_n)^T\nonumber
\end{align}
and
\begin{align}
\xi=(1,\cdots,1)^T,\nonumber
\end{align}
a $(N+n)$-dimensional vector. Define
$G_1$ to be the $N\times N$ martix such that
\begin{align}
(i,i)\mbox{-element of }G_1=(N-1)\sum_{i=1}^n(m_i^*-2)^2\int_{\mathbb{R}^2}e^{v_i}dx+\sum_{i=1}^n\int_{\mathbb{R}^2}e^{v_i}(\nabla v_i\cdot x+2)^2 dx\nonumber
\end{align}
and
\begin{align}
(i,j)\mbox{-element of }G_1=(N-2)\sum_{i=1}^n(m_i^*-2)^2\int_{\mathbb{R}^2}e^{v_i}dx.\nonumber
\end{align}
for $i\neq j$.
Define the $n\times N$ matrix
\[G_3=
\begin{pmatrix}
  (N-1)(2-m_1^*)\int_{\mathbb{R}^2}e^{v_1} & \cdots &  (N-1)(2-m_1^*) \int_{\mathbb{R}^2}e^{v_1}\\
  \vdots & \cdots & \vdots \\
  (N-1)(2-m_n^*) \int_{\mathbb{R}^2}e^{v_n} & \cdots & (N-1)(2-m_n^*) \int_{\mathbb{R}^2}e^{v_n} \\
\end{pmatrix}.
\]
Define the $n\times n$ diagonal matrix
\begin{align}
G_2=\mbox{diag}\Big(N\int_{\mathbb{R}^2}e^{v_1},\cdots, N\int_{\mathbb{R}^2}e^{v_n}\Big),\nonumber
\end{align}
which is invertible.
Then, (\ref{e:CoeffientEquation1}) and (\ref{e:CoeffientEquation2}) together give us
\[
\begin{pmatrix}
G_1+o_\varepsilon(1) & G_3^T+o_\varepsilon(1) \\
G_3+o_\varepsilon(1) & G_2+o_\varepsilon(1)
\end{pmatrix}
X=O(\varepsilon^{m^*-2}\sup|w_\varepsilon|+\|\phi_\varepsilon\|_{\mathbb{Y}_\varepsilon})\xi.
\]

By a similar computation as in (\ref{l:Coefficients}), the matrix
\[
\begin{pmatrix}
G_1 & G_3^T \\
G_3 & G_2
\end{pmatrix}.
\]
is invertible.
Then, for $t=1,\cdots, N$ and any $i=1,\cdots,n$, we get
\begin{align}
r_{3,t}=O(\varepsilon^{m^*-3})\sup|w_\varepsilon|+\frac{1}{\varepsilon}O(\|\phi_\varepsilon\|_{\mathbb{Y}_\varepsilon})\nonumber
\end{align}
and
\begin{align}
s_i=O(\varepsilon^{m^*-2})\sup|w_\varepsilon|+O(\|\phi_\varepsilon\|_{\mathbb{Y}_\varepsilon}).\nonumber
\end{align}
This proves (\ref{Ineq:r3}) and (\ref{Ineq:si}).
$\Box$

Before proving Theorem \ref{t:invertibility}, let us recall Theorem \ref{t:KernelLocal} and Corollary \ref{coro:KernelLocal}.

\begin{theorem}
Suppose $\phi=(\phi_1,\cdots,\phi_n)$ satisfies
\begin{equation}
\left\{
\begin{array}{lr}
-\Delta(\sum_{l=1}^n a^{il}\phi_l)=e^{v_i(y)}\phi_i\mbox{ in }\mathbb{R}^2,\nonumber\\
|\phi_i(x)|\leq C(1+|x|)^\tau\mbox{ for }i=1,\cdots,n\nonumber
\end{array}
\right.
\end{equation}
for small $\tau>0$. Then, $\phi$ is a linear combination of
\begin{align}
(\partial_{x_j}v_1,\cdots,\partial_{x_j}v_n)\mbox{ for }j=1,2\nonumber
\end{align}
and
\begin{align}
(|x|\dot{v}_1(|x|)+2,\cdots,|x|\dot{v}_n(|x|)+2).\nonumber
\end{align}
\end{theorem}

\begin{corollary}
Under the assumptions  {in} the previous theorem, if
\begin{align}
\sum_{i=1}^n\int_{\mathbb{R}^2} e^{v_i}\partial_{x_j}v_i\phi_idx=0\mbox{ for }j=1,2\nonumber
\end{align}
and
\begin{align}
\sum_{i=1}^n\int_{\mathbb{R}^n}e^{v_i}(|x|\dot{v}_i(|x|)+2)\phi_idx=0,\nonumber
\end{align}
then $\phi_i\equiv0$ for $i=1,\cdots,n$.
\end{corollary}

~

~

\noindent{\bf Proof of Theorem \ref{t:invertibility}.}
We prove Theorem \ref{t:invertibility} by contradiction. To this end, we assume there exists a sequence of positive numbers $\varepsilon^{(m)}\in\mathbb{R}^+$, two sequence of functions $w_{\varepsilon^{(m)}}\in E_{\varepsilon^{(m)}}$ and $\phi_{\varepsilon^{(m)}}\in F_{\varepsilon^{(m)}}$ such that
\begin{equation}\label{Assumption:invert}
    \left\{
   \begin{array}{lr}
\varepsilon^{(m)}\to0+,\\
Q_{\varepsilon^{(m)}}L_{\varepsilon^{(m)}} w_{\varepsilon^{(m)}}=\phi_{\varepsilon^{(m)}}:=(\phi_{1,\varepsilon^{(m)}},\cdots, \phi_{n,\varepsilon^{(m)}}),\\
\|w_{\varepsilon^{(m)}}\|_{\mathbb{X}_{\varepsilon^{(m)}}}+\|w_{\varepsilon^{(m)}}\|_{(L^\infty(\Omega))^n}=1,\\
\|\phi_{\varepsilon^{(m)}}\|_{\mathbb{Y}_{\varepsilon^{(m)}}}=o\Big(\ln\frac{1}{\varepsilon^{(m)}}\Big)
   \end{array}
   \right.
\end{equation}
as $m\to\infty$. For convenience, throughout the rest this proof, we suppress the $m$ in $\varepsilon^{(m)}$ and only write $\varepsilon$. We begin by rescaling (\ref{e:expansion}). For any $t,s=1,\cdots,N$, any $i=1,\cdots,n$ and any $j=1,2,3$, we denote
\begin{align}
\widetilde{w}_{t,i,\varepsilon}(y)&:=w_{i,\varepsilon}(\varepsilon_t y+p_{t,\varepsilon}),\nonumber\\
\widetilde{w}^i_{t,\varepsilon}(y)&:=\sum_{j=1}^n a^{ij}w_{j,\varepsilon}(\varepsilon_t y+p_{t,\varepsilon}),\nonumber\\
\widetilde{w}_{t,\varepsilon}(y)&:=(\widetilde{w}_{t,1,\varepsilon}(y), \cdots,\widetilde{w}_{t,n,\varepsilon}(y)),\nonumber\\
\widetilde{\phi}_{t,i,\varepsilon}(y)&:=\phi_{i,\varepsilon}(\varepsilon_t y+p_{t,\varepsilon}),\nonumber\\
\widetilde{K}_{t,i,\varepsilon}(y)&:=K_{i,\varepsilon}(\varepsilon_t y+p_{t,\varepsilon}),\nonumber\\
\widetilde{\chi}_{t,s}(y)&:=\chi_s(\varepsilon_t y+p_{t,\varepsilon}),\nonumber\\
\widetilde{Z}_{t,j,s,\varepsilon}^*&:= Z^*_{j,s,\varepsilon}(\varepsilon_t y+p_{t,\varepsilon}),\nonumber\\
\Omega_{t,\varepsilon}&:=\{y|\varepsilon_t y+p_{t,\varepsilon}\in\Omega\}\nonumber
\end{align}
and
\begin{align}
\widetilde{L}_{t,i,\varepsilon}\widetilde{w}_{t,\varepsilon}:=-\Delta_y\Big(\sum_{l=1}^n a^{il}\widetilde{w}_{t,l,\varepsilon}\Big)-\varepsilon_t^2\widetilde{K}_{t,i,\varepsilon}\widetilde{w}_{t,i,\varepsilon}+\frac{\varepsilon_t^2\widetilde{K}_{t,i,\varepsilon}}{\rho_i^*}\int_{\Omega_{t,\varepsilon}} \widetilde{K}_{t,i,\varepsilon}(y)\widetilde{w}_{t,i,\varepsilon}(y)dy.\nonumber
\end{align}
Based on these notations, we can rewrite
(\ref{e:expansion}) into
\begin{align}\label{e:expansionNew}
\widetilde{L}_{t,i,\varepsilon} \widetilde{w}_{t,\varepsilon} &=\varepsilon_t^2\widetilde{\phi}_{t,i,\varepsilon}+ \sum_{j=1}^3 \sum_{s=1}^N r_{j,t}\widetilde{Z}_{t,j,s,i,\varepsilon}^*\cdot\sum_{l=1}^N\widetilde{\chi}_{t,l} e^{v_i(y+\frac{p_{t,\varepsilon} -p_{l,\varepsilon}}{\varepsilon_l})+2\ln\frac{\varepsilon_t}{\varepsilon_l}}+s_i \sum_{l=1}^N\widetilde{\chi}_{t,l} e^{v_i(y+\frac{p_{t,\varepsilon} -p_{l,\varepsilon}}{\varepsilon_l})+2\ln\frac{\varepsilon_t}{\varepsilon_l}}.
\end{align}
The rest of this proof is proceeded in the following steps.

~

\noindent{\bf Step 1. We show that for any $R>0$, $\max_{i=1,\cdots,n,t=1,\cdots,N}\sup_{x\in B_R(0)}\widetilde{w}_{t,i,\varepsilon}(x)\to 0$ as $\varepsilon\to0$.}

~

By (\ref{Assumption:invert}), we know that $|\widetilde{w}_{t,i,\varepsilon}|=|w_{t,i,\varepsilon}|\leq 1$ for any $t=1,\cdots,N$ and any $i=1,\cdots,n$.
On the other hand, by Lemma \ref{Assumption:invert} and (\ref{Assumption:invert}), for $p\in(1,2)$ and any $R>0$, we get
\begin{align}
\varepsilon^2_t\|\widetilde{\phi}_{t,i,\varepsilon}\|_{L^p(B_R(0))}&=\varepsilon^2_t\Big(\int_{B_R}|\widetilde{\phi}_{t,i,\varepsilon}|^p \rho^p\frac{dy}{\rho^p}\Big)^\frac{1}{p}=\Big(\varepsilon^4_t\int_{B_R(0)}|\widetilde{\phi}_{t,i,\varepsilon}|^2\rho^2 dy\Big)^\frac{1}{2}\Big(\int_{\mathbb{R}^2}\frac{dy}{\rho^\frac{2p}{2-p}}\Big)^\frac{2-p}{2p}\nonumber\\
&\leq C\|\phi_\varepsilon\|_{\mathbb{Y}_\varepsilon}=o(\frac{1}{\ln\frac{1}{\varepsilon}}).\nonumber
\end{align}
For any $j=1,2$ and any $t=1,\cdots,N$,
\begin{align}
&\|r_{j,t}\widetilde{Z}_{t,j,s,i,\varepsilon}^*\cdot\sum_{l=1}^N\widetilde{\chi}_{t,l} e^{v_i(y+\frac{p_{t,\varepsilon} -p_{l,\varepsilon}}{\varepsilon_l})+2\ln\frac{\varepsilon_t}{\varepsilon_l}}\|_{L^P(B_R(0))}\nonumber\\
&\quad\leq Cr_{j,t}=O(\|\phi_\varepsilon\|_{\mathbb{Y}_\varepsilon})+O(\varepsilon^{m^*-2})\cdot\max_{i=1,\cdots,n}\max_{t=1,\cdots,N}\sup_{y\in B_{2\delta}(p_{t,\varepsilon_t})}|w_{i,\varepsilon}(y)|=o(\frac{1}{\ln\frac{1}{\varepsilon}}),\nonumber
\end{align}
\begin{align}
&\|r_{3,t}\widetilde{Z}_{t,3,s,i,\varepsilon}^*\cdot\sum_{l=1}^N\widetilde{\chi}_{t,l} e^{v_i(y+\frac{p_{t,\varepsilon} -p_{l,\varepsilon}}{\varepsilon_l})+2\ln\frac{\varepsilon_t}{\varepsilon_l}}\|_{L^P(B_R(0))}\nonumber\\
&\quad\leq C\varepsilon r_{3,t}=O(\|\phi_\varepsilon\|_{\mathbb{Y}_\varepsilon})+O(\varepsilon^{m^*-2})\max_{i=1,\cdots,n}\max_{t=1,\cdots,N}\sup_{y\in B_{2\delta}(p_{t,\varepsilon_t})}|w_{i,\varepsilon}(y)|=o(\frac{1}{\ln\frac{1}{\varepsilon}})\nonumber
\end{align}
and
\begin{align}
\|s_i & \sum_{l=1}^N\widetilde{\chi}_{t,l} e^{v_i(y+\frac{p_{t,\varepsilon} -p_{l,\varepsilon}}{\varepsilon_l})+2\ln\frac{\varepsilon_t}{\varepsilon_l}}\|_{L^p(B_R(0))}=O(s_i)\nonumber\\
&=O(\|\phi_\varepsilon\|_{\mathbb{Y}_\varepsilon})+O(\varepsilon^{m^*-2})\cdot\max_{i=1,\cdots,n}\max_{t=1,\cdots,N}\sup_{y\in B_{2\delta}(p_{t,\varepsilon_t})}|w_{i,\varepsilon}(y)|=o(\frac{1}{\ln\frac{1}{\varepsilon}}).\nonumber
\end{align}
By the classical elliptic regularity (see \cite{GilbargTrudingerBook2001}), there exists bounded continuous functions $\widetilde{w}_{t,i}$ for $t=1,\cdots,N$ and $i=1,\cdots,n$ and a positive constant $\alpha\in(0,1)$ such that
$\widetilde{w}_{t,i,\varepsilon}\to\widetilde{w}_{t,i}\mbox{ in }C^\alpha_{loc}(\mathbb{R}^2)$.
Moreover, $\widetilde{w}_{t,i}$ satisfies
\begin{align}
-\Delta\Big(\sum_{l=1}^n a^{il}\widetilde{w}_{t,l}\Big)=e^{v_i(y)}\widetilde{w}_{t,i}+\frac{e^{v_i(y)}}{\rho_i^*}\int_{\mathbb{R}^2}\sum_{s=1}^N\Bigg[\frac{e^{m_1\sum_{l=1}^N  {G^*}(p_{t}^*,p_{l}^*)}h_i(p_{t}^*)}{e^{m_1\sum_{l=1}^N  {G^*}(p_{s}^*,p_{l}^*)}h_i(p_{s}^*)}\Bigg]^\frac{2}{m_i^*-2} e^{v_i(x)}\widetilde{w}_{s,i}dx.\nonumber
\end{align}
for $i=1,\cdots,n$ and $t=1,\cdots,N$.
Recalling that $w_\varepsilon\in E_\varepsilon$, we get
\begin{align}
\int_{\mathbb{R}^2}\sum_{s=1}^N\Bigg[\frac{e^{m_1\sum_{l=1}^N  {G^*}(p_{t}^*,p_{l}^*)}h_i(p_{t}^*)}{e^{m_1\sum_{l=1}^N  {G^*}(p_{s}^*,p_{l}^*)}h_i(p_{s}^*)}\Bigg]^\frac{2}{m_i^*-2} e^{v_i(x)}\widetilde{w}_{s,i}dx=0\nonumber
\end{align}
for any $i=1,\cdots,n$.
Then, $\widetilde{w}_{t,i}\equiv0$ for any $t=1,\cdots,N$ and any $i=1,\cdots,n$ due to Corollary \ref{coro:KernelLocal}. Then,
$\widetilde{w}_{t,i,\varepsilon}\to0$ in $C_{loc}^\alpha(\mathbb{R}^2)$
as $\varepsilon\to0$. The result follows immediately.

~

~

\noindent{\bf Step 2. For any $\theta>0$, we show that there exists a small constant $\beta_\theta>0$ such that $\max_{i=1,\cdots,n}\max_{t=1,\cdots,N}\|\widetilde{w}_{t,i,\varepsilon}\|_{L^\infty(B_\frac{\beta_\theta}{\varepsilon_t}(0))}\leq\theta$.}

~

We argue by contradiction. Suppose
there exists a constant $\theta_0>0$, a sequence of $\widetilde{w}_\varepsilon$ and $\beta_\varepsilon\to0$ such that
\begin{align}\label{ASSUMPTION02}
\max_{i=1,\cdots,n}\max_{t=1,\cdots,N}\|\widetilde{w}_{t,i,\varepsilon}\|_{L^\infty(B_\frac{\beta_\varepsilon}{\varepsilon}(0))}\geq\theta_0.
\end{align}
Let $x_\varepsilon\in B_{\frac{\beta_\varepsilon}{\varepsilon}}(0)$ be the point where $\widetilde{w}_{t_0,i}$ achieve $\max_{i=1,\cdots,n}\max_{t=1,\cdots,N}\|\widetilde{w}_{t,i,\varepsilon}\|_{L^\infty(B_\frac{\beta_\varepsilon}{\varepsilon}(0))}$ for some $t_0=1,\cdots,N$. Without loss of generality, we assume that $t_0=1$. This implies that
\begin{align}\label{Ineq:Sizeofx}
|x_\varepsilon|\leq\frac{\beta_\varepsilon}{\varepsilon}=o_\varepsilon(1).
\end{align}
By {\bf Step 1}, we know that
$|x_\varepsilon|\to\infty$.
We define the rescaled functions as follows. For $t=1,\cdots,N$ and $i=1,\cdots,n$,
\begin{align}
\overline{w}_{t,i,\varepsilon}(y):&
=\widetilde{w}_{t,i,\varepsilon}(|x_\varepsilon|y),\nonumber\\
\overline{w}_{t,\varepsilon}(y)&:=(\overline{w}_{t,1,\varepsilon}(y), \cdots,\overline{w}_{t,n,\varepsilon}(y)),\nonumber\\
\overline{\phi}_{t,i,\varepsilon}(y)&:=\widetilde{\phi}_{t,i,\varepsilon}(|x_\varepsilon|y),\nonumber\\
\overline{K}_{t,i,\varepsilon}(y)&:=\widetilde{K}_{t,i,\varepsilon}(|x_\varepsilon|y),\nonumber\\
\overline{\chi}_{t,s}(y)&:=\widetilde{\chi}_{t,s}(|x_\varepsilon|y),\nonumber\\
\overline{Z}_{t,j,s,\varepsilon}^*(y)&:=\widetilde{Z}_{t,j,s,\varepsilon}^*(|x_\varepsilon|y)\nonumber
\end{align}
and
\begin{align}
\overline{L}_{t,i,\varepsilon}\overline{w}_{t,\varepsilon}:=-\Delta_y\Big(\sum_{l=1}^n a^{il}\overline{w}_{t,l,\varepsilon}\Big)-|x_\varepsilon|^2\overline{K}_{t,i,\varepsilon}\overline{w}_{t,i,\varepsilon}+\frac{|x_\varepsilon|^2\overline{K}_{t,i,\varepsilon}}{\rho_i^*}\int_{\frac{1}{|x_\varepsilon|^2}\Omega_{t,\varepsilon}} \overline{K}_{t,i,\varepsilon}(y)\overline{w}_{t,i,\varepsilon}(y)dy.\nonumber
\end{align}
Then, we get
\begin{align}\label{e:expansionNew01}
\overline{L}_{t,i,\varepsilon} \overline{w}_{t,\varepsilon} &=\varepsilon_t^2|x_\varepsilon|^2\overline{\phi}_{t,i,\varepsilon}+ \sum_{j=1}^3 \sum_{s=1}^N r_{j,t}|x_\varepsilon|^2\overline{Z}_{t,j,s,i,\varepsilon}^*\cdot\sum_{l=1}^N\overline{\chi}_{t,l} e^{v_i(|x_\varepsilon|y+\frac{p_{t,\varepsilon} -p_{l,\varepsilon}}{\varepsilon_l})+2\ln\frac{\varepsilon_t}{\varepsilon_l}}\nonumber\\
&\quad\quad+s_i|x_\varepsilon|^2 \sum_{l=1}^N\overline{\chi}_{t,l} e^{v_i(|x_\varepsilon|y+\frac{p_{t,\varepsilon} -p_{l,\varepsilon}}{\varepsilon_l})+2\ln\frac{\varepsilon_t}{\varepsilon_l}}.
\end{align}
For any $R_2>R_2>0$, as $\varepsilon\to0+$, we notice that
\begin{align}
|x_\varepsilon|^2\overline{K}_{t,i,\varepsilon}\overline{w}_{t,i,\varepsilon}=O(|x_\varepsilon|^{2-m^*})\mbox{ in }L^\infty(B_{R_2}(0)\backslash B_{R_1}(0))\mbox{ sense}\nonumber
\end{align}
and
\begin{align}
\frac{|x_\varepsilon|^2\overline{K}_{t,i,\varepsilon}}{\rho_i^*}\int_{\frac{1}{|x_\varepsilon|^2}\Omega_{t,\varepsilon}} \overline{K}_{t,i,\varepsilon}(y)\overline{w}_{t,i,\varepsilon}(y)dy=O(|x_\varepsilon|^{2-m^*})\mbox{ in }L^\infty(B_{R_2}(0)\backslash B_{R_1}(0))\mbox{ sense.}\nonumber
\end{align}
For $p\in(1,2)$, we get
\begin{align}
\|\varepsilon_t^2|x_\varepsilon|^2\overline{\phi}_{t,i,\varepsilon}\|_{L^p(B_{R_2}(0)\backslash B_{R_1}(0))}&=\varepsilon_t^2|x_\varepsilon|^2\Big(\int_{B_{R_2}(0)\backslash B_{R_1}(0)}|\phi_{i,\varepsilon}(\varepsilon_t|x_\varepsilon|y+p_{t,\varepsilon})|^p dy\Big)^\frac{1}{p}\nonumber\\
&\leq \varepsilon_t^2|x_\varepsilon|^{2-\frac{2}{p}}\Big(\int_{\frac{1}{|x_\varepsilon|}\cdot B_{R_2}(0)\backslash B_{R_1}(0)} |\phi_{i,\varepsilon}(\varepsilon y+p_{t,\varepsilon})|^2\rho(y)^2 dy\Big)^\frac{1}{2}\times\nonumber\\
&\quad\quad\times\Big(\int_{\frac{1}{|x_\varepsilon|}\cdot B_{R_2}(0)\backslash B_{R_1}(0)}\frac{dy}{\rho(y)^{\frac{2p}{2-p}}}\Big)^\frac{2-p}{2p}\nonumber
\end{align}
By (\ref{Ineq:Sizeofx}), we get
\begin{align}
\|\varepsilon_t^2|x_\varepsilon|^2\overline{\phi}_{t,i,\varepsilon}\|_{L^p(B_{R_2}(0)\backslash B_{R_1}(0))}&\leq \varepsilon_t^2|x_\varepsilon|^{2-\frac{2}{p}}\Big(\int_{B_{\frac{2\delta}{\varepsilon}}(0)} |\phi_{i,\varepsilon}(\varepsilon y+p_{t,\varepsilon})|^2\rho(y)^2 dy\Big)^\frac{1}{2}\times\nonumber\\
&\quad\quad\times\Big(\int_{\frac{1}{|x_\varepsilon|}\cdot B_{R_2}(0)\backslash B_{R_1}(0)}\frac{dy}{\rho(y)^{\frac{2p}{2-p}}}\Big)^\frac{2-p}{2p}\nonumber\\
&\leq C |x_\varepsilon|^\frac{\alpha}{2}\|\phi_{i}\|_{\mathbb{Y}_\varepsilon}=o_\varepsilon(1).\nonumber
\end{align}
Moreover, for $j=1,2$ we get
\begin{align}
\|&\sum_{s=1}^N r_{j,t}|x_\varepsilon|^2\overline{Z}_{t,j,s,i,\varepsilon}^*\cdot\sum_{l=1}^N\overline{\chi}_{t,l} e^{v_i(|x_\varepsilon|y+\frac{p_{t,\varepsilon} -p_{l,\varepsilon}}{\varepsilon_l})+2\ln\frac{\varepsilon_t}{\varepsilon_l}}\|_{L^p(B_{R_2}(0)\backslash B_{R_1}(0))}\nonumber\\
&\leq C\sum_{s=1}^N r_{j,t}|x_\varepsilon|^{1-m^*}=O(\|\phi_\varepsilon\|_{\mathbb{Y}_\varepsilon}+\varepsilon^{m^*-2}\max_{i=1,\cdots,n}\sup_{y\in B_{2\delta}(p_{t,\varepsilon_t})}|w_{i,\varepsilon}(y)|)|x_\varepsilon|^{1-m^*}=o_\varepsilon(1),\nonumber
\end{align}
\begin{align}
\|&\sum_{s=1}^N r_{j,t}|x_\varepsilon|^2\overline{Z}_{t,3,s,i,\varepsilon}^*\cdot\sum_{l=1}^N\overline{\chi}_{t,l} e^{v_i(|x_\varepsilon|y+\frac{p_{t,\varepsilon} -p_{l,\varepsilon}}{\varepsilon_l})+2\ln\frac{\varepsilon_t}{\varepsilon_l}}\|_{L^p(B_{R_2}(0)\backslash B_{R_1}(0))}\nonumber\\
&\leq C\varepsilon\sum_{s=1}^N r_{3,t}|x_\varepsilon|^{2-m^*}=O(\|\phi_\varepsilon\|_{\mathbb{Y}_\varepsilon}+\varepsilon^{m^*-2}\max_{i=1,\cdots,n}\sup_{y\in B_{2\delta}(p_{t,\varepsilon_t})}|w_{i,\varepsilon}(y)|)|x_\varepsilon|^{2-m^*}=o_\varepsilon(1)\nonumber
\end{align}
and
\begin{align}
\|&\sum_{s=1}^N s_i|x_\varepsilon|^2 \sum_{l=1}^N\overline{\chi}_{t,l} e^{v_i(|x_\varepsilon|y+\frac{p_{t,\varepsilon} -p_{l,\varepsilon}}{\varepsilon_l})+2\ln\frac{\varepsilon_t}{\varepsilon_l}}\|_{L^p(B_{R_2}(0)\backslash B_{R_1}(0))}\nonumber\\
&\leq C\sum_{s=1}^N s_1|x_\varepsilon|^{2-m^*}=O(\|\phi_\varepsilon\|_{\mathbb{Y}_\varepsilon}+\varepsilon^{m^*-2}\max_{i=1,\cdots,n}\sup_{y\in B_{2\delta}(p_{t,\varepsilon_t})}|w_{i,\varepsilon}(y)|)|x_\varepsilon|^{2-m^*}=o_\varepsilon(1).\nonumber
\end{align}

Then, there exists a bounded harmonic function $\overline{w}_{t,i,\infty}$ such that
\begin{align}
\overline{w}_{t,i,\varepsilon}\to \overline{w}_{t,i,\infty}\mbox{ in }L_{loc}^\infty(\mathbb{R}^2\backslash \{0\}).\nonumber
\end{align}
Moreover, by Liouville theorem, $\overline{w}_{t,i,\infty}$ is a constant. Denote it by $c_{t,i}$.
By (\ref{ASSUMPTION02}), we get
\begin{align}\label{Ineq:cLowerBdd}
\max_{t=1,\cdots,N}\max_{i=1,\cdots,n}c_{t,i}\geq\theta_0>0.
\end{align}
In other words, for any $t=1,\cdots,N$ and any $i=1,\cdots,n$, we get
\begin{align}\label{Ineq:wconstant}
\widetilde{w}_{t,i,\varepsilon}(x)=c_{t,i}+o_\varepsilon(1)\mbox{ for }R_1|x_\varepsilon|\leq |y|\leq R_2|x_\varepsilon|.
\end{align}
Define
\begin{align}
\hat{w}_{t,i,\varepsilon}(s)=\frac{1}{2\pi}\int_0^{2\pi}\widetilde{w}_{t,i,\varepsilon}(s(\cos\theta,\sin\theta))d\theta.\nonumber
\end{align}
We get
\begin{align}
\hat{w}_{t,i,\varepsilon}(s)=c_{t,i}+o_\varepsilon(1)\mbox{ for }R_1|x_\varepsilon|\leq s\leq R_2|x_\varepsilon|.\nonumber
\end{align}
immediately.

On the other hand, (\ref{e:expansionNew}) implies that
\begin{align}
-\Delta & (\sum_{l=1}^n a^{ij} \hat{w}_{t,j,\varepsilon}) -e^{v_i(y)}\hat{w}_{t,i,\varepsilon}+\frac{e^{v_i(y)}}{\rho_i^*}\int_{\Omega_{t,\varepsilon}} \widetilde{K}_{t,i,\varepsilon}(y)\widetilde{w}_{t,i,\varepsilon}(y)dy\nonumber\\
&=\frac{\varepsilon_t^2}{2\pi}\int_0^{2\pi}\widetilde{\phi}_{\varepsilon,i}(s\cos\theta,s\sin\theta)d\theta+\frac{r_{3,t}e^{v_i(y)}}{2\pi}\int_{0}^{2\pi}\chi_{t}Z_{3,i,\varepsilon}(s\cos\theta,s\sin\theta) +\frac{s_i}{2\pi}\chi_t e^{v_i(y)}.\nonumber
\end{align}
Notice that $w_\varepsilon\in E_\varepsilon$, we get
$\int_{\Omega_{t,\varepsilon}} \widetilde{K}_{t,i,\varepsilon}(y)\widetilde{w}_{t,i,\varepsilon}(y)dy=0$.
Therefore,
\begin{align}\label{e:ExpansionNew02}
-\Delta & (\sum_{l=1}^n a^{ij} \hat{w}_{t,j,\varepsilon}) -e^{v_i(y)}\hat{w}_{t,i,\varepsilon}\nonumber\\
&=\frac{\varepsilon_t^2}{2\pi}\int_0^{2\pi}\widetilde{\phi}_{\varepsilon,i}(s\cos\theta,s\sin\theta)d\theta+\frac{r_{3,t}e^{v_i(y)}}{2\pi}\int_{0}^{2\pi}\chi_{t}Z_{3,i,\varepsilon}(s\cos\theta,s\sin\theta) +\frac{s_i}{2\pi}\chi_t e^{v_i(y)}.
\end{align}
Multiplying the both sides of (\ref{e:ExpansionNew02}) by $s$ and integrating from $0$ to $t$, we get
\begin{align}
-t\frac{d}{dt} & (\sum_{j=1}^n a^{ij}\hat{w}_{t,i,\varepsilon}(t)) -\int_{B_t(0)}e^{v_i(y)}\hat{w}_{t,i,\varepsilon}(|y|)dy\nonumber\\
&=\int_{B_t(0)}\widetilde{\phi}_{i,\varepsilon}dy+r_{3,t}\int_{B_t(0)}\chi_t e^{v_i(y)}\widetilde{Z}_{3,i,t,\varepsilon}(y)dy +s_i\int_{B_t(0)}\chi_t e^{v_i(y)}dy.\nonumber
\end{align}
In other words,
\begin{align}\label{e:ExpansionNew031}
-t\frac{d}{dt} & \hat{w}_{t,i,\varepsilon}(t) -\int_{B_t(0)}\sum_{j=1}^n a_{ij}e^{v_j(y)}\hat{w}_{t,j,\varepsilon}(|y|)dy\nonumber\\
&=\int_{B_t(0)}\sum_{j=1}^n a_{ij}\widetilde{\phi}_{j,\varepsilon}dy+r_{3,t}\int_{B_t(0)}\chi_t \sum_{j=1}^n a_{ij}e^{v_i(y)}\widetilde{Z}_{3,j,t,\varepsilon}(y)dy +\sum_{j=1}^n a_{ij}s_j\int_{B_t(0)}\chi_t e^{v_j(y)}dy.
\end{align}
In order to proceed with our proof, we need the following estimate.

\begin{claim}\label{c:wlndecay}
For any small $\tau>0$, there exsits a positive constant $C_\tau$ such that $$\max_{t=1,\cdots,N,i=1,\cdots,n}|w_{t,i,\varepsilon}(y)|\leq\frac{o_{\varepsilon}(1)\cdot
C_\tau(1+|y|)^\tau}{\ln\frac{1}{\varepsilon}}.$$
\end{claim}

We will sketch the proof of Claim \ref{c:wlndecay} at the end of this appendix since this is similar to Proposition \ref{prop:w0Pointwise}.

By a similar computation as in {\bf Step 1}, we get
\begin{align}
\mbox{RHS of (\ref{e:ExpansionNew031})}&=O(\|\phi_\varepsilon\|_{\mathbb{Y}_\varepsilon})+\varepsilon^{-1} O_{m^*}\cdot\max_{i=1,\cdots,n}\max_{t=1,\cdots,N}\sup_{y\in B_{2\delta}(p_{t,\varepsilon_t})}|w_{i,\varepsilon}(y)|=o_\varepsilon\Big(\frac{1}{\ln\frac{1}{\varepsilon}}\Big
)+\varepsilon^{-1}O_{m^*}.\nonumber
\end{align}
On the other hand, for $t$ sufficiently large, we get
\begin{align}
\int_{B_t(0)}\sum_{j=1}^n a_{ij}e^{v_j(y)}\hat{w}_{t,j,\varepsilon}(|y|)dy& \leq \frac{o_\varepsilon(1)}{\ln\frac{1}{\varepsilon}}
\int_{B_t(0)}\frac{dy}{(1+|y|)^{m^*-\tau}}=\frac{o_\varepsilon(1)}{\ln\frac{1}{\varepsilon}}\nonumber
\end{align}
for sufficiently small $\tau$.
Then, we get
$\frac{d}{dt}\hat{w}_{t,j,\varepsilon}(t)=\frac{o_\varepsilon(1)}{\ln\frac{1}{\varepsilon}}\frac{1}{t}$.
Integrating the both side from $R$ to $|x_\varepsilon|$, we get
\begin{align}
|\hat{w}_{t,j,\varepsilon}(|x_\varepsilon|)|\leq|\hat{w}_{t,j,\varepsilon}(R)|+\frac{o_\varepsilon(1)\ln\frac{|x_\varepsilon|}{R}}{\ln\frac{1}{\varepsilon}}.\nonumber
\end{align}
Recall that $|x_\varepsilon|=\frac{o_\varepsilon(1)}{\varepsilon}$, as we assumed. For sufficiently large $R$ and small $\varepsilon$, we get
$|\hat{w}_{t,j,\varepsilon}(|x_\varepsilon|)|\leq\frac{\theta_0}{2}$
for any $t=1,\cdots,N$ and any $i=1,\cdots,n$. Here, $\theta_0$ is the positive constant appears in (\ref{ASSUMPTION02}).
This contradicts with (\ref{Ineq:cLowerBdd}) and (\ref{Ineq:wconstant}).

~

\noindent{\bf Step 3. Estimating $\|w_\varepsilon^i\|_{L^\infty(\Omega)}+\|w_\varepsilon^i\|_{\mathbb{X}_\varepsilon}$.}

~

Now we come back to (\ref{e:expansion}), i.e.
\begin{align}
L_{i,\varepsilon} w_\varepsilon &=\phi_{i,\varepsilon}+ \sum_{j=1}^2 \sum_{t=1}^N r_{j,t}\partial_{x_j} v_i( \frac{x-p_{t,\varepsilon}}{\varepsilon_t})\cdot\sum_{s=1}^N\chi_s e^{v_i(\frac{x-p_{s,\varepsilon}}{\varepsilon_s})+2\ln\frac{1}{\varepsilon_s}}\nonumber\\
&\quad+\sum_{t=1}^N r_{3,t}Z_{3,t,i,\varepsilon}^*\cdot\sum_{s=1}^N\chi_s e^{v_i(\frac{x-p_{s,\varepsilon}}{\varepsilon_s})+2\ln\frac{1}{\varepsilon_s}}+s_i \sum_{s=1}^N\chi_s e^{v_i(\frac{x-p_{s,\varepsilon}}{\varepsilon_s})+2\ln\frac{1}{\varepsilon_s}}\nonumber
\end{align}
for $i=1,\cdots,n$. Here,
\begin{align}
L_{i,\varepsilon}(w_\varepsilon^1,\cdots,w_\varepsilon^n) = -\Delta w_\varepsilon^i-K_{i,\varepsilon}(x)w_{i,\varepsilon} +\frac{K_{i,\varepsilon}(x)}{\rho_i^*}\int_\Omega K_{i,\varepsilon}(x) w_{i,\varepsilon}.\nonumber
\end{align}
For the sake of simplicity, we write $B:=\cup_{t=1}^N B_{\delta_t}(p_{t,\varepsilon})$.
By the classical elliptic estimate (see, for instance, \cite{GilbargTrudingerBook2001}), we get
\begin{align}
\|w^i_\varepsilon\|_{L^\infty(\Omega\backslash B)}&\leq C\|w^i_\varepsilon\|_{L^\infty(\partial B)}+C\|\phi_{i,\varepsilon}\|_{L^p(\Omega\backslash B)}\nonumber\\
&\quad+C\Big\|\sum_{j=1}^2 \sum_{t=1}^N r_{j,t}\partial_{x_j} v_i( \frac{x-p_{t,\varepsilon}}{\varepsilon_t})\cdot\sum_{s=1}^N\chi_s e^{v_i(\frac{x-p_{s,\varepsilon}}{\varepsilon_s})+2\ln\frac{1}{\varepsilon_s}}\Big\|_{L^p(\Omega\backslash B)}\nonumber\\
&\quad+C\Big\|\sum_{t=1}^N r_{3,t}Z_{3,t,i,\varepsilon}^*\cdot\sum_{s=1}^N\chi_s e^{v_i(\frac{x-p_{s,\varepsilon}}{\varepsilon_s})+2\ln\frac{1}{\varepsilon_s}}\Big\|_{L^p(\Omega\backslash B)}\nonumber\\
&\quad+C\Big\|\sum_{t=1}^N r_{3,t}Z_{3,t,i,\varepsilon}^*\cdot\sum_{s=1}^N\chi_s e^{v_i(\frac{x-p_{s,\varepsilon}}{\varepsilon_s})+2\ln\frac{1}{\varepsilon_s}}\Big\|_{L^p(\Omega\backslash B)}\nonumber\\
&\quad+C\Big\|s_i \sum_{s=1}^N\chi_s e^{v_i(\frac{x-p_{s,\varepsilon}}{\varepsilon_s})+2\ln\frac{1}{\varepsilon_s}}\Big\|_{L^p(\Omega\backslash B)} +C\Big\|\frac{K_{i,\varepsilon}(x)}{\rho_i^*}\int_\Omega K_{i,\varepsilon}(x) w_{i,\varepsilon}\Big\|_{L^p(\Omega\backslash B)}.\nonumber
\end{align}
Here, $\|w^i_\varepsilon\|_{L^\infty(\partial B)}$ is sufficiently small due to {\bf Step 2}. Denote it by $\theta$. Moreover, we get
\begin{align}
\|\phi_{\varepsilon,i}\|_{L^p(\Omega\backslash B)}\leq C\|\phi_{\varepsilon,i}\|_{\mathbb{Y}_\varepsilon}=o_\varepsilon\Big(\frac{1}{\ln\frac{1}{\varepsilon}}\Big),\nonumber
\end{align}
\begin{align}
\Big\|r_{j,t}\partial_{x_j} v_i( \frac{x-p_{t,\varepsilon}}{\varepsilon_t})\cdot\sum_{s=1}^N\chi_s e^{v_i(\frac{x-p_{s,\varepsilon}}{\varepsilon_s})+2\ln\frac{1}{\varepsilon_s}}\Big\|_{L^p(\Omega\backslash B)}\leq C(\|\phi_{i,\varepsilon}\|_{\mathbb{Y}_\varepsilon}+\varepsilon^{m^*-2})\cdot \varepsilon^{\frac{2}{p}-2+m^*-2}=o_\varepsilon(1),\nonumber
\end{align}
\begin{align}
\Big\|&r_{3,t}Z_{3,t,i,\varepsilon}^*\cdot\sum_{s=1}^N\chi_s e^{v_i(\frac{x-p_{s,\varepsilon}}{\varepsilon_s})+2\ln\frac{1}{\varepsilon_s}}\Big\|_{L^p(\Omega\backslash B)}\nonumber\\
&\leq \Big[O(\varepsilon^{-1}\|\phi_\varepsilon\|_{\mathbb{Y}_\varepsilon})+\varepsilon^{m^*-2} \cdot\max_{i=1,\cdots,n}\max_{t=1,\cdots,N}\sup_{y\in B_{2\delta}(p_{t,\varepsilon_t})}|w_{i,\varepsilon}(y)|\Big]\varepsilon^{\frac{2}{p}+m^*-3}=o_\varepsilon(1),\nonumber
\end{align}
\begin{align}
\Big\|s_i & \sum_{s=1}^N\chi_s e^{v_i(\frac{x-p_{s,\varepsilon}}{\varepsilon_s})+2\ln\frac{1}{\varepsilon_s}}\Big\|_{L^p(\Omega\backslash B)}\nonumber\\
&\leq \Big[O(\|\phi_\varepsilon\|_{\mathbb{Y}_\varepsilon})+\varepsilon^{-1} O_{m^*}\cdot\max_{i=1,\cdots,n}\max_{t=1,\cdots,N}\sup_{y\in B_{2\delta}(p_{t,\varepsilon_t})}|w_{i,\varepsilon}(y)|\Big]\varepsilon^{\frac{2}{p}+m^*-4}=o_\varepsilon(1)\nonumber
\end{align}
and
$\int_\Omega K_{i,\varepsilon}(x) w_{i,\varepsilon}=0$
since $w_\varepsilon\in E_\varepsilon$.
Here, we used Lemma \ref{l:constantsrrs}. Therefore,
\begin{align}\label{Ineq:wLinfty}
\|w_\varepsilon^i\|_{L^\infty(\Omega)}\leq \theta+o_\varepsilon(1).
\end{align}

Now we estimate $\|w_\varepsilon^i\|_{\mathbb{X}_\varepsilon}$.
By (\ref{e:expansionNew}), for any $t=1,\cdots,N$ we get
\begin{align}
\|\Delta_y\widetilde{w}^i_{t,\varepsilon} {\widetilde{\rho}_\beta}\|_{L^2(B_\frac{2\delta}{\varepsilon_t}(0))}&\leq \int_{B_\frac{2\delta}{\varepsilon}}| {\varepsilon_t^2\widetilde{K}_{t,i,\varepsilon}\widetilde{w}_{t,i,\varepsilon}}\cdot(1+|y|)^{1+ {\frac{\beta}{2}}}|^2dy+C\varepsilon_t^2\int_{B_\frac{2\delta}{\varepsilon}}\sum_{j=1}^n|\phi_{t,j,\varepsilon}(\varepsilon_t y+p_{t,\varepsilon})|^2dy\nonumber\\
&\quad+\|\sum_{j=1}^3 \sum_{s=1}^N r_{j,t}\widetilde{Z}_{t,j,s,i,\varepsilon}^*\cdot\sum_{l=1}^N\widetilde{\chi}_{t,l} e^{v_i(y+\frac{p_{t,\varepsilon} -p_{l,\varepsilon}}{\varepsilon_l})+2\ln\frac{\varepsilon_t}{\varepsilon_l}}\|^2_{L^2({B_\frac{2\delta}{\varepsilon}}(0))}\nonumber\\
&\quad+\|s_i \sum_{l=1}^N\widetilde{\chi}_{t,l} e^{v_i(y+\frac{p_{t,\varepsilon} -p_{l,\varepsilon}}{\varepsilon_l})+2\ln\frac{\varepsilon_t}{\varepsilon_l}}\|^2_{L^2(B_\frac{2\delta}{\varepsilon}(0))}.\nonumber
\end{align}
Notice that
\begin{align}
\int_{B_\frac{2\delta}{\varepsilon}}| {\varepsilon_t^2\widetilde{K}_{t,i,\varepsilon}\widetilde{w}_{t,i,\varepsilon}}\cdot(1+|y|)^{1+ {\frac{\beta}{2}}}|^2dy&\leq C\sum_{j=1}^n\int_{B_\frac{2\delta}{\varepsilon}}\frac{(1+|y|)^{ {2+\beta}}}{(1+|y|)^{2m^*}}|\widetilde{w}^i_{t,\varepsilon}(\varepsilon_t y+p_{t,\varepsilon})|^2 dy\nonumber\\
&\leq C {|w_\varepsilon|^2_{L^\infty(B_{2\delta})}},\nonumber
\end{align}
\begin{align}
\varepsilon_t^2\int_{B_\frac{2\delta}{\varepsilon}}\sum_{j=1}^n|\phi_{t,j,\varepsilon}(\varepsilon_t y+p_{t,\varepsilon})|^2dy\leq  {\|\phi_\varepsilon\|^2_{\mathbb{Y}_\varepsilon}}\nonumber
\end{align}
and
\begin{align}
&\|\sum_{j=1}^3 \sum_{s=1}^N r_{j,t}\widetilde{Z}_{t,j,s,i,\varepsilon}^*\cdot\sum_{l=1}^N\widetilde{\chi}_{t,l} e^{v_i(y+\frac{p_{t,\varepsilon} -p_{l,\varepsilon}}{\varepsilon_l})+2\ln\frac{\varepsilon_t}{\varepsilon_l}}\|^2_{L^2({B_\frac{2\delta}{\varepsilon}}(0))}\nonumber\\
&\quad+\|s_i \sum_{l=1}^N\widetilde{\chi}_{t,l} e^{v_i(y+\frac{p_{t,\varepsilon} -p_{l,\varepsilon}}{\varepsilon_l})+2\ln\frac{\varepsilon_t}{\varepsilon_l}}\|^2_{L^2(B_\frac{2\delta_t}{\varepsilon_t}(0))}\leq [\varepsilon^{-1}O_{m^*}|w_\varepsilon|_{L^\infty(\Omega)}+O(\|\phi_\varepsilon\|_{\mathbb{Y}_\varepsilon})]^2.\nonumber
\end{align}
On the other hand, we get
\begin{align}
\|\widetilde{w}_{t,\varepsilon}^i\hat{\rho}\|^2_{L^2(B_\frac{2\delta}{\varepsilon}(0))}&\leq\int_{B_\frac{2\delta}{\varepsilon}}|\widetilde{w}^i_{t,\varepsilon}(\varepsilon_t y+p_{t,\varepsilon})|^2\cdot\frac{dy}{(1+|y|)^2\big(\ln(2+|y|)\big)^{2+\alpha}}\leq C {|w_{t,\varepsilon}^i|_{L^\infty(B_{2\delta}(0))}^2}.\nonumber
\end{align}
By a standard elliptic estimate (\cite{GilbargTrudingerBook2001}) on $\Omega\backslash\cup_{t=1}^N B_{2\delta_t}(p_{t,\varepsilon})$, we get
\begin{align}
\|w_\varepsilon^i\|_{\mathbb{X}_\varepsilon}\leq C(\|w_\varepsilon^i\|_{L^\infty(\Omega)}+\|\phi_\varepsilon\|_{\mathbb{Y}_\varepsilon})\leq C\Big(\theta+o\Big(\frac{1}{|\ln\varepsilon|}\Big)\Big).\nonumber
\end{align}
Combining (\ref{Ineq:wLinfty}), we obtain a contradiction with (\ref{Assumption:invert}).
Now (\ref{Ineq:Invert}) is proved.

~~

\noindent{\bf Step 4. $Q_\varepsilon$ is an isomorphism.}

The result follows from the Fredholm alternative.
$\Box$

~~

Now we are in the position to prove Claim \ref{c:wlndecay}.

~

\noindent{\bf Proof of Claim \ref{c:wlndecay}.}
We only sketch the proof since this is similar to Proposition \ref{prop:w0Pointwise}.
To prove this, inspired by \cite{LinZhang2013,BartolucciYangZhang2024}, we argue by contradiction as follows. Denote
\begin{align}
\Lambda_{\varepsilon}:=\max_{t=1,\cdots,n,i=1,\cdots,N}\frac{\ln\frac{1}{\varepsilon}\cdot\hat{w}_{t,i,\varepsilon}(y)}{(1+|y|)^\gamma}.\nonumber
\end{align}
If $\Lambda_\varepsilon=o_\varepsilon(1)$, then the proof is completed.
Suppose that there exists a constant $c_0>0$ such that $\Lambda_\varepsilon\geq c_0$ and $y_\varepsilon\in B_{\frac{2\delta_t}{\varepsilon_t}}(0)$ attends $\Lambda_\varepsilon$. Denote
\begin{align}
\hat{\hat{w}}_{t,i,\varepsilon}(y):=\frac{\ln\frac{1}{\varepsilon}\cdot\hat{w}_{t,i,\varepsilon}(y)}{\Lambda_\varepsilon (1+|y_\varepsilon|)^\tau}.\nonumber
\end{align}
An immediate observation is that
\begin{align}
|\hat{\hat{w}}_{t,i,\varepsilon}(y)|=\Big|\frac{\ln\frac{1}{\varepsilon}\cdot\hat{w}_{t,i,\varepsilon}(y)}{\Lambda_\varepsilon(1+|y|)^\tau}\cdot\frac{(1+|y|)^\tau}{(1+|y_\varepsilon|)^\tau}\Big|\leq\frac{(1+|y|)^\tau}{(1+|y_\varepsilon|)^\tau}.\nonumber
\end{align}
Then,
\begin{align}\label{e:ExpansionNew03}
-\Delta_y\Big(\sum_{j=1}^n a^{ij}\hat{\hat{w}}_{t,j,\varepsilon}\Big)-e^{v_i(y)}\hat{\hat{w}}_{t,j,\varepsilon}=\frac{\varepsilon_t^2\ln\frac{1}{\varepsilon}\cdot\hat{\phi}_{\varepsilon,i}}{\Lambda_\varepsilon (1+|y_\varepsilon|)^\tau}+\frac{\ln\frac{1}{\varepsilon}\cdot r_{3,t}e^{v_i(y)}\widetilde{Z}_{3,t,i,\varepsilon}^*}{\Lambda_\varepsilon (1+|y_\varepsilon|)^\tau}+\frac{\ln\frac{1}{\varepsilon}\cdot s_i e^{v_i(y)}}{\Lambda_\varepsilon (1+|y_\varepsilon|)^\tau}.
\end{align}
Here,
\begin{align}
\hat{\phi}_{\varepsilon,i}(s)=\frac{1}{2\pi}\int_0^{2\pi}\widetilde{\phi}_{\varepsilon,i}(s\cdot\cos\theta,s\cdot\sin\theta)d\theta.\nonumber
\end{align}

Our argument is divided into three parts with respect to the growth of $|y_\varepsilon|$. For all of the three cases, we will draw contradictions.

~

\noindent{\bf Case 1. $|y_\varepsilon|$ is bounded.}

~

By a similar approach as in {\bf Step 1} of the proof of Theorem \ref{t:invertibility}, we get for any $R>0$
\begin{align}
\Big\|\frac{\varepsilon_t^2\ln\frac{1}{\varepsilon}\hat{\phi}_{\varepsilon,i}}{\Lambda_\varepsilon (1+|y_\varepsilon|)^\tau}\Big\|_{L^p(B_R(0))}=\Big\|\frac{\ln\frac{1}{\varepsilon}\cdot r_{3,t}e^{v_i(y)}\widetilde{Z}_{3,t,i,\varepsilon}^*}{\Lambda_\varepsilon (1+|y_\varepsilon|)^\tau}\Big\|_{L^p(B_R(0))}=\Big\|\frac{\ln\frac{1}{\varepsilon}\cdot s_i e^{v_i(y)}}{\Lambda_\varepsilon (1+|y_\varepsilon|)^\tau}\Big\|_{L^p(B_R(0))}=o_\varepsilon(1),\nonumber
\end{align}
This implies that
$\hat{\hat{w}}_{t,i,\varepsilon}\to w^*_{t,i}\mbox{ in }C_{loc}^\beta(\mathbb{R}^2)$
with $\beta\in(0,1)$ and
\begin{equation}
\left\{
\begin{array}{lr}
-\Delta\Big(\sum_{j=1}^n a^{ij}w_{t,j}^*\Big)-e^{v_i(y)}w_{t,i}^*=0\mbox{ in }\mathbb{R}^2,\\
|w_{t,i}^*(y)|\leq C(1+|y|)^\tau\mbox{ and }w_{t,i}^*\mbox{ are radial};\\
\sum_{i=1}^n\int_{\mathbb{R}^2}e^{v_i(y)}(y\cdot\nabla v_i(y)+2)w_{t,i}^*(y)dy=0\nonumber\\
w_{t,i}^*(y_0)=1\mbox{ with }y_0=\lim_{\varepsilon\to0}y_\varepsilon.
\end{array}
\right.
\end{equation}
Here, the third assertion is due to the fact that $w_\varepsilon\in E_\varepsilon$. This contradicts with Corollary \ref{coro:KernelLocal}.

~

\noindent{\bf Case 2. $|y_\varepsilon|$ is unbounded and $|y_\varepsilon|=o_\varepsilon(1)\varepsilon^{-1}$.}

~

By a similar computation as in {\bf Case 1}, we know that $\hat{\hat{w}}_{t,i,\varepsilon}(0)\to0$ as $\varepsilon\to0$. On the other hand, we always have
$\hat{\hat{w}}_{t,i,\varepsilon}(y_\varepsilon)=\pm 1$.
Then, by Green's formula, we get
\begin{align}
\frac{1}{2}&\leq|\hat{\hat{w}}_{t,i,\varepsilon}(y_\varepsilon)-\hat{\hat{w}}_{t,i,\varepsilon}(0)|\nonumber\\
&\leq\Bigg\{\int_{B_{\frac{2\delta}{\varepsilon}}(0)}\Big|G_\varepsilon(y_\varepsilon,\eta) -G_\varepsilon(0,\eta)\Big|\cdot \Big|\sum_{j=1}^n a_{ij}e^{v_j(\eta)}\hat{\hat{w}}_{t,j,\varepsilon}(\eta)\Big|d\eta\Bigg\}\nonumber\\
&\quad+\Bigg\{\int_{B_{\frac{2\delta}{\varepsilon}}(0)}\Big|G_\varepsilon(y_\varepsilon,\eta) -G_\varepsilon(0,\eta)\Big|\cdot\Big|\frac{\varepsilon_t^2\ln\frac{1}{\varepsilon}\cdot\sum_{j=1}^n a_{ij}\hat{\phi}_{\varepsilon,j}(\eta)}{\Lambda_\varepsilon (1+|y_\varepsilon|)^\tau}+\frac{\ln\frac{1}{\varepsilon}\cdot r_{3,t}\sum_{j=1}^n a_{ij}e^{v_j(\eta)}\widetilde{Z}_{3,t,j,\varepsilon}^*(\eta)}{\Lambda_\varepsilon (1+|y_\varepsilon|)^\tau}\nonumber\\
&\quad\quad\quad+\frac{\ln\frac{1}{\varepsilon}\cdot \sum_{j=1}^n s_j e^{v_j(y)}}{\Lambda_\varepsilon (1+|y_\varepsilon|)^\tau}\Big|\Bigg\}d\eta=:A+B.
\end{align}
Here, $G_\varepsilon(y,\eta)$ denotes the Green's function with Dirichlet boundary condition on $B_{\frac{2\delta_t}{\varepsilon_t}}(0)$. The boundary terms are canceled out since $\widetilde{w}^*_{t,i,\varepsilon}$ is radial.
We will draw a contradiction by proving that $A+B=o_\varepsilon(1)$.

Notice that we have the pointwise estimate
\begin{align}\label{Ineq:Pointwise}
\Big|e^{v_j(\eta)}\hat{\hat{w}}_{t,j,\varepsilon}(\eta)\Big|\leq\frac{C}{(1+|y_\varepsilon|)^\tau(1+|\eta|)^{m^*-\tau}}.
\end{align}
On the other hand, by similar computations as in (\ref{ineq:gLp1}), we get for $p\in(1,2)$
\begin{align}\label{Ineq:Lp01}
\Bigg(\int_{B_{\frac{2\delta}{\varepsilon}}(0)}\Bigg|\frac{\varepsilon_t^2\ln\frac{1}{\varepsilon}\cdot\hat{\phi}_{\varepsilon,j}(\eta)}{\Lambda_\varepsilon (1+|y_\varepsilon|)^\tau}\Bigg|^pd\eta
\Bigg)^\frac{1}{p}
=\frac{o_\varepsilon(1)}{\Lambda_\varepsilon(1+|y_\varepsilon|)^\tau}.
\end{align}
Here, we use (\ref{Assumption:invert}). As in (\ref{ineq:Z3Lp1}) and (\ref{ineq:esLp1}) we get for $p\in(1,2)$
\begin{align}\label{Ineq:Lp02}
\Bigg(\int_{B_\frac{2\delta}{\varepsilon}(0)}\Bigg|\frac{\ln\frac{1}{\varepsilon}\cdot r_{j,t}e^{v_i(\eta)}\widetilde{Z}_{j,t,i,\varepsilon}^*(\eta)}{\Lambda_\varepsilon (1+|y_\varepsilon|)^\tau} \Bigg|^p\Bigg)^\frac{1}{p}=\frac{o_\varepsilon(1)}{\Lambda_\varepsilon(1+|y_\varepsilon|)^\tau}
\end{align}
for $j=1,2,3$
and
\begin{align}\label{Ineq:Lp03}
\Bigg(\int_{B_\frac{2\delta}{\varepsilon}(0)}\Bigg|\frac{\ln\frac{1}{\varepsilon}\cdot s_i e^{v_i(\eta)}}{\Lambda_\varepsilon (1+|y_\varepsilon|)^\tau}\Bigg|^p d\eta\Bigg)^\frac{1}{p}=\frac{o_\varepsilon(1)}{\Lambda_\varepsilon(1+|y_\varepsilon|)^\tau}.
\end{align}
Here, we use (\ref{Assumption:invert}) and Lemma \ref{l:constantsrrs}.
To estimate $A+B$, we apply \cite[Theorem 3.2]{LinZhang2013} (see also \cite[Proposition 3.1]{BartolucciYangZhang2024}). By \cite[Lemma 3.2]{LinZhang2013}, we get
\begin{equation}\label{Ineq:GreenEstimates}
    |G_\varepsilon(y_\varepsilon,
    \eta)-G_\varepsilon(0,\eta)| \leq \left\{
    \begin{aligned}
    & C(\ln|y_\varepsilon| +|\ln|\eta||)&\mbox{ if }\eta\in\Sigma_1:=\{\eta\in B_{\frac{2\delta}{\varepsilon}}(0)||\eta|<\frac{|y_\varepsilon|}{2}\}; \\
    &C(\ln|y_\varepsilon|+|\ln|y-\eta||)&\mbox{ if }\eta\in\Sigma_2:=\{\eta\in B_{\frac{2\delta}{\varepsilon}}(0)||y_\varepsilon-\eta|<\frac{|y_\varepsilon|}{2}\};\\
    &\frac{C|y_\varepsilon|}{|\eta|}&\mbox{ if }\eta\in\Sigma_3:= B_{\frac{2\delta}{\varepsilon}(0)}\backslash(\Sigma_1\cup\Sigma_2).
    \end{aligned}
    \right.
\end{equation}
In follows, we estimate $A+B$ with the help of (\ref{Ineq:Pointwise}), (\ref{Ineq:Lp01}), (\ref{Ineq:Lp02}), (\ref{Ineq:Lp03}) and (\ref{Ineq:GreenEstimates}).

As for $A$, by a similar method as in (\ref{ineq:ATOTAL1}), (\ref{ineq:A11}), (\ref{ineq:A21}), (\ref{ineq:A31}) and (\ref{Ineq:A}), we get
\begin{align}
A&\leq C\int_{B_\frac{2\delta}{\varepsilon}(0)}\frac{|G_\varepsilon(y_\varepsilon,\eta)-G_\varepsilon(0,\eta)|d\eta}{(1+|y_\varepsilon|)^\tau(1+|\eta|)^{m^*-\tau}}=C\Big\{ \int_{\Sigma_1}+\int_{\Sigma_2}+\int_{\Sigma_3}\Big\}\frac{|G_\varepsilon(y_\varepsilon,\eta)-G_\varepsilon(0,\eta)|d\eta}{(1+|y_\varepsilon|)^\tau(1+|\eta|)^{m^*-\tau}}=O_\varepsilon(1).\nonumber
\end{align}

On the other hand, by a similar method as in (\ref{ineq:BTOTAL1}), (\ref{ineq:B11}), (\ref{ineq:B21}) and (\ref{ineq:B31}), we get
\begin{align}
B&\leq \frac{o_\varepsilon(1)}{\Lambda_\varepsilon(1+|y_\varepsilon|)^\tau}\Big(\int_{B_\frac{2\delta}{\varepsilon}}|G_\varepsilon(y_\varepsilon,\eta)-G_\varepsilon(0,\eta)|^{p'}d\eta\Big)^\frac{1}{p'}\nonumber\\
&\leq \frac{o_\varepsilon(1)}{\Lambda_\varepsilon(1+|y_\varepsilon|)^\tau}\Big[\Big(\int_{\Sigma_1}|G_\varepsilon(y_\varepsilon,\eta)-G_\varepsilon(0,\eta)|^{p'}d\eta\Big)^\frac{1}{p'}+\Big(\int_{\Sigma_2}|G_\varepsilon(y_\varepsilon,\eta)-G_\varepsilon(0,\eta)|^{p'}d\eta\Big)^\frac{1}{p'}\nonumber\\
&\quad\quad\quad+\Big(\int_{\Sigma_3}|G_\varepsilon(y_\varepsilon,\eta)-G_\varepsilon(0,\eta)|^{p'}d\eta\Big)^\frac{1}{p'}\Big]=:o_\varepsilon(1)\nonumber
\end{align}
Hence, we get
$A+B=o_\varepsilon(1)$, which is a contradiction.

~

\noindent{\bf Case 3. $|y_\varepsilon|$ is unbounded and $|y_\varepsilon|\simeq\varepsilon^{-1}$.}

~
Combining \cite[Proposition 3.1]{BartolucciYangZhang2024},
we obtain a contradiction by a similar approach as in {\bf Case 3} of Proposition \ref{prop:w0Pointwise}. Claim
\ref{c:wlndecay} is proved.
$\Box$

~~

~~

~~

\noindent{\bf  {Acknowledgment}} The authors would like to express their  sincere gratitude to the anonymous reviewers.

\noindent{\bf  {{\bf Author Contributions.}}}  {The authors are listed in alphabetical order of their surnames, and all authors have made equal contributions.}

\noindent{\bf  { Competing Interests.}} The authors have no relevant financial or non-financial competing interests.

\noindent{\bf  { Data availability statement.}} No datasets were generated or analysed during  the current study.

\noindent{\bf  { Fundings.}}
Haoyu Li was supported by FAPESP Proc 2022/15812-0.   {Lei Zhang is partially supported by Simon's Foundation grant SFI-MPS-TSM-00013752}

\end{document}